\documentclass[a4paper,11pt]{amsart}
\usepackage{amssymb,graphicx}
\newtheorem{lemma}{Lemma}
\newtheorem{prop}{Proposition}
\newtheorem{thm}{Theorem}

\theoremstyle{definition}
\newtheorem{defn}{Definition}
\theoremstyle{remark}
\newtheorem{rem}{Remark}
\newtheorem{ex}{Example}
\newcounter{numl}

\newcommand{\labelnuml}{\textup{(\roman{numl})}}
\newenvironment{numlist}{\begin{list}{\labelnuml}%
{\usecounter{numl}\setlength{\leftmargin}{0pt}%
\setlength{\itemindent}{2\parindent}%
\setlength{\itemsep}{\smallskipamount}\def
\makelabel ##1{\hss \llap {\upshape ##1}}}}{\end{list}}
\newenvironment{bulletlist}{\begin{list}{\labelitemi}%
{\setlength{\leftmargin}{\parindent}\def
\makelabel ##1{\hss \llap {\upshape ##1}}}}{\end{list}}
\newenvironment{numlproof}{\begin{proof}\begin{numlist}}
{\qed \end{numlist}\begingroup\let\qed\relax\end{proof}\endgroup}
\newcommand{\R}{{\mathbb R}}
\newcommand{\C}{{\mathbb C}}
\newcommand{\Z}{{\mathbb Z}}

\newcommand{\Q}{{\mathbb Q}}
\newcommand{\T}{{\mathbb T}}
\newcommand{\J}{{\mathcal J}}
\newcommand{\cH}{{\mathcal H}}
\newcommand{\cB}{{\mathcal B}}
\newcommand{\cC}{{\mathcal C}}
\newcommand{\cF}{{\mathcal F}}
\newcommand{\cL}{{\mathcal L}}

\newcommand{\cK}{{\mathcal K}}
\newcommand{\cO}{{\mathcal O}}
\newcommand{\cV}{{\mathcal V}}
\newcommand{\cX}{{\mathcal X}}

\newcommand{\Scal}{\mathit{Scal}}
\newcommand{\Id}{\mathit{Id}}
\newcommand{\ka}{K{\"a}hler }

\newcommand{\trace}{\mathop{\mathrm{tr}}\nolimits}

\newcommand{\pfaff}{\mathop{\mathrm{pf}}\nolimits}
\newcommand{\grad}{\mathop{\mathrm{grad}}\nolimits}
\newcommand{\g}[1]{g(#1)}
\newcommand{\Hess}{\mathrm{Hess}}
\newcommand{\symprod}{\mathbin{\raise1pt\hbox{$\scriptstyle\bigcirc$}}}

\newcommand{\restr}[1]{|_{#1}^{\vphantom x}}
\newcommand{\Fa}{F}
\newcommand{\tFa}{F'}
\newcommand{\pF}{F}
\newcommand{\oF}{\Theta}
\newcommand{\Mp}{p}
\newcommand{\Mpc}{p_{\mathrm{c}}}
\newcommand{\Mpn}{p_{\mathrm{nc}}}

\newcommand{\cc}{c}

\newcommand{\ang}{t}
\begin{document}
\title[Hamiltonian 2-forms in K{\smash{\"a}}hler geometry, II]
{Hamiltonian 2-forms in K{\smash{\"a}}hler geometry,\\
II Global Classification}
\author[V. Apostolov]{Vestislav Apostolov}
\address{Vestislav Apostolov \\ D{\'e}partement de Math{\'e}matiques\\
UQAM\\ C.P. 8888 \\ Succ. Centre-ville \\ Montr{\'e}al (Qu{\'e}bec) \\
H3C 3P8 \\ Canada}
\email{apostolo@math.uqam.ca}
\author[D. Calderbank]{David M.~J.~Calderbank}
\address{David M. J. Calderbank \\ School of Mathematics\\
University of Edinburgh\\ King's Buildings\\ Mayfield Road\\
Edinburgh EH9 3JZ\\ Scotland}
\email{davidmjc@maths.ed.ac.uk}
\author[P. Gauduchon]{Paul Gauduchon}
\address{Paul Gauduchon \\ Centre de Math\'ematiques\\
{\'E}cole Polytechnique \\ UMR 7640 du CNRS
\\ 91128 Palaiseau \\ France}
\email{pg@math.polytechnique.fr}
\author[C. T\o nnesen-Friedman]{Christina W.~T\o nnesen-Friedman}
\address{Christina W. T\o nnesen-Friedman\\
Department of Mathematics\\ Union College\\
Schenectady\\  New York  12308\\ USA }
\email{tonnesec@union.edu}
\thanks{We would like to thank C.~Boyer, S.~Boyer, R.~Bryant, O.~Collin,
K.~Galicki, D.~Guan and S.~Maillot for stimulating discussions. The first
author was supported in part by FCAR grant NC-7264, and by NSERC grant
OGP0023879, the second author by the Leverhulme Trust, the William Gordon
Seggie Brown Trust and an EPSRC Advanced Research Fellowship. The first three
authors are members of EDGE, Research Training Network HPRN-CT-2000-00101,
supported by the European Human Potential Programme.}
\date{\today}
\begin{abstract} We present a classification of compact K\"ahler manifolds
admitting a hamiltonian $2$-form (which were classified locally in part I of
this work). This involves two components of independent interest.

The first is the notion of a rigid hamiltonian torus action.  This natural
condition, for torus actions on a K\"ahler manifold, was introduced locally in
part I, but such actions turn out to be remarkably well behaved globally,
leading to a fairly explicit classification: up to a blow-up, compact K\"ahler
manifolds with a rigid hamiltonian torus action are bundles of toric K\"ahler
manifolds.

The second idea is a special case of toric geometry, which we call orthotoric.
We prove that orthotoric K\"ahler manifolds are diffeomorphic to complex
projective space, but we extend our analysis to orthotoric orbifolds, where
the geometry is much richer.  We thus obtain new examples of
K\"ahler--Einstein $4$-orbifolds.

Combining these two themes, we prove that compact K\"ahler manifolds with
hamiltonian $2$-forms are covered by blow-downs of projective bundles over
K\"ahler products, and we describe explicitly how the K\"ahler metrics with a
hamiltonian $2$-form are parameterized. We explain how this provides a context
for constructing new examples of extremal K\"ahler metrics---in particular a
subclass of such metrics which we call weakly Bochner-flat.

We also provide a self-contained treatment of the theory of compact toric
K\"ahler manifolds, since we need it and find the existing literature
incomplete.
\end{abstract}
\maketitle

This paper is concerned with the construction of explicit K\"ahler metrics on
compact manifolds, and has several interrelated motivations. The first is the
notion of a \emph{hamiltonian $2$-form}, introduced in part I of this
series~\cite{ACG1}.

\begin{defn} Let $\phi$ be any (real) $J$-invariant $2$-form on the K{\"a}hler
manifold $(M, g, J,\omega)$ of dimension $2m$.  We say $\phi$ is {\it
hamiltonian} if
\begin{equation} \label{ham}
\nabla _X \phi = \tfrac12 (d \trace\phi \wedge JX - d^c\trace\phi  \wedge X)
\end{equation}
for any vector field $X$, where $\trace\phi=\langle \phi,\omega\rangle$ is the
trace with respect to $\omega$. When $M$ is a Riemann surface $(m=1)$, this
equation is vacuous and we require instead that $\trace \phi$ is a Killing
potential, i.e., a hamiltonian for a Killing vector field $J\grad_g
\trace\phi$.
\end{defn}

A second motivation is the notion of a \emph{weakly Bochner-flat} (WBF)
K\"ahler metric, by which we mean a K\"ahler metric whose Bochner tensor
(which is part of the curvature tensor) is co-closed. By the differential
Bianchi identity, this is equivalent (for $m\geq 2$) to the condition that
$\rho + \frac{\Scal}{2 (m + 1)} \, \omega$ is a hamiltonian $2$-form, where
$\rho$ is the Ricci form.  WBF K\"ahler metrics are extremal in the sense of
Calabi, i.e., the symplectic gradient of the scalar curvature is a Killing
vector field, and provide a class of extremal K\"ahler metrics which include
the Bochner-flat K\"ahler metrics studied by Bryant~\cite{bryant} and products
of K\"ahler--Einstein metrics. The geometry of WBF K\"ahler metrics is tightly
constrained, because the more specific the normalized Ricci form is, the
closer the metric is to being K\"ahler--Einstein, while the more generic it
is, the stronger the consequences of the hamiltonian property.

A hamiltonian $2$-form $\phi$ induces an isometric hamiltonian $\ell$-torus
action on $M$ for some $0\leq\ell\leq m$, which we call the order of
$\phi$. This says nothing for $\ell=0$, but for $\ell=m$, it means that $M$ is
toric.  Toric K\"ahler manifolds are well understood, and a third motivation
for our work is to extend this understanding to certain torus actions with
$0<\ell<m$.  We introduce the notion of a \emph{rigid} hamiltonian
$\ell$-torus action and prove that a compact K\"ahler manifold with such an
action has a blow-up which is biholomorphic to a bundle of toric K\"ahler
$2\ell$-manifolds.

We shall be particularly interested in the projective bundles of the form
$M=P(\cL_0\otimes\C^{d_0+1}\oplus\cdots\oplus\cL_\ell\otimes\C^{d_\ell+1})\to
S$, where $\cL_0,\ldots \cL_\ell$ are line bundles over a compact K\"ahler
manifold $S$ and the $\ell$-torus action is induced by scalar multiplication
on the vector bundles $\cL_j\otimes\C^{d_j+1}$, with $d_j\geq 0$.  The blow-up
of $M$ along the submanifolds determined by setting the $j$th fibrewise
homogeneous coordinate (in $\cL_j\otimes\C^{d_j+1}$) to zero, for
$j=0,\ldots\ell$, is a bundle of toric K\"ahler $2\ell$-manifolds: the
projective bundle $P(\tilde\cL_0\oplus\cdots\oplus\tilde\cL_\ell)\to \C
P^{d_0}\times\cdots \times\C P^{d_\ell}\times S$, where $\tilde
\cL_j=\cO(0,\ldots 0,-1,0,\ldots 0)\otimes\cL_j$ (with $\cO(-1)$ over the
$j$th factor $\C P^{d_j}$).

When $\ell=1$, projective line bundles have been well-used, since the seminal
work of Calabi~\cite{calabi1}, in the construction of explicit examples of
extremal K\"ahler metrics. The idea to consider blow-downs was introduced by
Koiso and Sakane~\cite{koi-sak1,koi-sak2}, who constructed K\"ahler--Einstein
metrics in this way. Our fourth motivation is to provide a general framework
for constructing extremal K\"ahler metrics on projective bundles and their
blow-downs, and in doing so we obtain new examples.

The toric K\"ahler $2m$-manifolds arising from hamiltonian $2$-forms of order
$m$ are of a special class, which we call orthotoric. Compact orthotoric
K\"ahler manifolds are necessarily biholomorphic to complex projective space,
but there are many more examples on orbifolds. Our final motivation is to
study K\"ahler metrics on toric orbifolds, especially orthotoric orbifolds,
and to obtain new examples.

The main goal of this paper is to show that a compact K\"ahler manifold with a
hamiltonian $2$-form of order $\ell$ is necessarily biholomorphic to a
projective bundle $M$ of the form described above, and conversely to show
precisely how to construct K\"ahler metrics with hamiltonian $2$-forms of
order $\ell$ on such bundles.

We hope however, that with the various motivations discussed above, the Reader
who does not share our enthusiasm for hamiltonian $2$-forms will find
something of interest in this paper.  Hamiltonian $2$-forms rather provide a
device that unifies and underlies the above themes. The journey to our main
result, and its consequences, yield a number of results of independent
interest.
\begin{bulletlist}
\item We obtain necessary and sufficient first order boundary conditions for
the compactification of compatible K\"ahler metrics on toric symplectic
orbifolds, clarifying work of Abreu~\cite{Abreu,Abreu1}, whose proofs we do
not understand (see Remark~\ref{angular-coordinates}
and~\S\ref{s:kahlercomp}).
\item We introduce and study rigid hamiltonian torus actions, and orthotoric
K\"ahler manifolds and orbifolds.
\item We construct new explicit K\"ahler--Einstein metrics on $4$-orbifolds.
\item We unify and extend constructions of K\"ahler metrics on projective
bundles, obtaining new weakly Bochner-flat and extremal K\"ahler metrics on
projective line bundles and on the projective plane bundle
$P(\cO\oplus\cO(1)\otimes\C^2)\to\C P^1$.
\end{bulletlist}

We have attempted to make this paper as independent as possible from the first
part~\cite{ACG1}. However, we shall make essential use of the local
classification of K\"ahler manifolds with a hamiltonian $2$-form of order
$\ell$, so we recall the result here. The Reader who is not interested in
hamiltonian $2$-forms {\it per se}, could take this local classification
result as a (rather complicated) definition of the class of K\"ahler metrics
that we wish to classify globally.

We define the {\it momentum polynomial} of a hamiltonian $2$-form $\phi$ to be
\begin{equation}\label{mompoly}
\Mp(t):=(-1)^m\pfaff (\phi-t\omega)
= t^m - (\trace\phi) \,t^{m-1} + \cdots + (-1)^m \pfaff\phi
\end{equation}
where the {\it pfaffian} is defined by $\phi \wedge \cdots \wedge \phi =
(\pfaff\phi)\omega\wedge\cdots\wedge\omega$.

\begin{thm}\label{ACGthm}\textup{\cite{ACG1}}
Let $(M,g,J,\omega)$ be a connected K\"ahler $2m$-manifold with a hamiltonian
$2$-form $\phi$. Then\textup:
\begin{numlist}
\item the functions $\Mp(t)$ on $M$ \textup(for each $t\in\R$\textup) are
Poisson-commuting hamiltonians for Killing vector fields $K(t) := J \grad_g
\Mp(t)$\textup;
\item there is a monic polynomial $\Mpc(t)$ with constant coefficients such
that $\Mp(t)=\Mpc(t)\Mpn(t)$ and, if $\Mpn(t)=\sum_{r=0}^\ell (-1)^r\sigma_r
t^{\ell-r}$ \textup(with $0\leq\ell\leq m$\textup), then the Killing vector
fields $K_r:= J\grad_g \sigma_r$ $(r=1,\ldots \ell)$ are linearly independent
on a connected dense open subset $M^0$ of $M$.  The integer $\ell$ is called
the \emph{order} of $\phi$.
\end{numlist}
On the open subset $M^0$, the roots $\xi_1,\ldots\xi_\ell$ of $\Mpn(t)$ are
smooth, functionally independent and everywhere pairwise distinct, and they
extend continuously to $M$. Denote by $\eta_a$, $a=1,\ldots N$ $(N\leq
m-\ell)$ the \emph{different} constant roots of $\Mpc(t)$ and by $d_a$ their
multiplicities. Then there are \textup(positive or negative definite\textup)
K\"ahler metrics $(g_a,\omega_a)$ of real dimension $2d_a$, functions
$\pF_1,\ldots \pF_\ell$ of one variable, and $1$-forms
$\theta_1,\ldots\theta_\ell$ with $\theta_r(K_s)=\delta_{rs}$ such that the
K\"ahler structure on $M^0$ is of the form
\begin{equation}\label{metric}\begin{split}
g&=\sum_{a=1}^N \Mpn(\eta_a) g_a
+\sum_{j=1}^\ell \frac{\Mp'(\xi_j)}{\pF_j(\xi_j)} d\xi_j^2
+\sum_{j=1}^\ell \frac{\pF_j(\xi_j)}{\Mp'(\xi_j)}\Bigl(\sum_{r=1}^\ell
\sigma_{r-1}(\hat\xi_j)\theta_r\Bigr)^2,\\
\omega&=\sum_{a=1}^N \Mpn(\eta_a)\omega_a
+\sum_{r=1}^\ell d\sigma_r\wedge \theta_r,\qquad\qquad
d\theta_r=\sum_{a=1}^N (-1)^r\eta_a^{\ell-r}\omega_a,
\end{split}\end{equation}
and the hamiltonian $2$-form $\phi$ is given by
\begin{equation}\label{phi}
\phi = \sum_{a=1}^{N} \eta_{a}\, \Mpn(\eta_a) \omega_{a}
+\sum_{r=1}^\ell (\sigma_r d\sigma_1 -d\sigma_{r+1})\wedge\theta_r
\end{equation}
with $\sigma_{\ell+1}=0$. \textup(Here $\sigma_{r-1}(\hat \xi_j)$ denote the
elementary symmetric functions of the roots with $\xi_j$ omitted. We remark
also that $\Mp'(\xi_j)=\Mpc(\xi_j)\prod_{k \neq j} (\xi_j-\xi_k)$.\textup)
\end{thm}

We shall obtain our global description of compact K\"ahler manifolds admitting
a hamiltonian $2$-form of order $\ell$ by exploiting three aspects of the
local geometry revealed by Theorem~\ref{ACGthm}.
\begin{numlist}
\item The components $\g{K_r,K_s}$ of the metric are constant on fibres of
the momentum map $(\sigma_1,\ldots\sigma_\ell)\colon M\to \R^\ell$. (This
holds on all of $M$ by continuity.)
\item The K\"ahler quotient metrics $\sum_{a=1}^N \Mpn(\eta_a) g_a$ are
simultaneously diagonalizable (with respect to $\sum_{a=1}^N g_a$) with
constant eigenvalues for each fixed $(\sigma_1,\ldots\sigma_\ell)$.
\item The roots $\xi_1,\ldots\xi_\ell$ of $\Mpn$ have orthogonal gradients.
\end{numlist}
In~\cite{ACG1}, these properties were interpreted by saying that
$(M,g,J,\omega)$ is given locally by a {\it rigid hamiltonian $\ell$-torus
action} with {\it semisimple K\"ahler quotient} and {\it orthotoric} fibres.
We shall see that this is not far from being true globally.

If $M$ is compact, the closure of the group of hamiltonian isometries of $M$
generated by $K_1,\ldots K_\ell$ is a torus $\T$ (with $\ell\leq \dim \T\leq
m$). When $\ell=m$, $K_1,\ldots K_m$ generate a torus action, and $M$ is a
toric K\"ahler manifold. In the first section we review the necessary
background of toric K\"ahler geometry and introduce a suitable invariant
language.  Then, in section~\ref{s:rigid}, we pursue a similar theory for
$\ell<m$ when property (i) holds. In particular we prove that $\dim\T=\ell$ so
there is a global rigid $\ell$-torus action. We provide a generalized Calabi
construction for such actions which classifies them up to covering when the
K\"ahler quotient is semisimple, i.e., when property (ii) is also
satisfied. In section~\ref{s:orthotoric}, we study toric K\"ahler manifolds
(and orbifolds) satisfying property (iii) in general, and here we exhibit new
explicit K\"ahler--Einstein metrics on compact $4$-orbifolds. In
section~\ref{s:global}, we obtain a complete description of compact K\"ahler
manifolds with hamiltonian $2$-forms, which we use to construct new examples
of compact weakly Bochner-flat and extremal K\"ahler manifolds. In subsequent
work we shall construct many more examples and classify weakly Bochner-flat
K\"ahler metrics in dimension $6$.

\tableofcontents

\section{Hamiltonian actions and toric geometry}\label{s:hamtoric}

We begin by reviewing hamiltonian torus actions, paying particular attention
to the theory of toric K\"ahler manifolds. Toric K\"ahler geometry can be
studied either from the complex or symplectic viewpoint, and we adopt,
primarily, the latter.  Furthermore, with a view to applications, we do not
restrict attention to manifolds, but also consider orbifolds: this is a
natural context in toric symplectic geometry~\cite{Abreu1,LT}. We refer
to~\cite{audin,GS,LT} for general information about torus actions on
symplectic manifolds, and to~\cite{Abreu,CDG,delzant,guillemin,guillemin0} for
further information about toric K\"ahler manifolds and orbifolds.

Our treatment has some novel features: in particular we obtain first order
boundary conditions for the compactification of compatible K\"ahler metrics on
toric symplectic manifolds.  Also, we present the theory in invariant
language, because for the torus actions generated by hamiltonian $2$-forms,
the natural basis of the Lie algebra $\mathfrak t$ is not (in general)
compatible with the lattice in $\mathfrak t$ defining the torus $\T$.

\subsection{Hamiltonian torus actions}

Let $\T$ be an $\ell$-dimensional torus, with Lie algebra $\mathfrak t$,
acting effectively on a symplectic $2m$-manifolds $(M,\omega)$, and for
$\xi\in\mathfrak t$ denote by $X_\xi$ the corresponding vector field on $M$.
Then we say that the action is \emph{hamiltonian} if there is a $\T$-invariant
smooth map $\mu\colon M\to {\mathfrak t}^*$, called a \emph{momentum map} for
the action, such that $\iota_{X_\xi}\omega=-\langle d\mu,\xi\rangle$ for any
$\xi\in\mathfrak t$.

\begin{rem} Note that our actions are hamiltonian in the strong sense
that $\mu$ is $\T$-invariant (if $\mu$ has a critical point---as it does in
the compact case---this is automatic). Since $\T$ is abelian, this implies
that $\omega(X_\xi,X_\eta)=0$ for any $\xi,\eta\in {\mathfrak t}$.  We also
remark that the action determines and is determined by $\mu$ up to a constant.
\end{rem}

We shall normally be interested in the case that $(M,\omega,\mu)$ has a
compatible almost K\"ahler structure, i.e., a $\T$-invariant metric $g$ and
almost complex structure $J$ with $\omega(X,Y)= g(JX,Y)$. Such compatible
metrics always exist.

We shall make significant use of the symplectic slice theorem for $\T$-orbits
in $M$, which we now recall. Let $\T\cdot x$ be such an orbit for $x\in M$.
Since $\T\cdot x$ is isotropic with respect to $\omega$, the isotropy
representation of $\T_x$ on $T_x M$ induces a $2(m-k)$-dimensional symplectic
representation on $V_x:=T_x(\T\cdot x)^0/T_x(\T\cdot x)$, where $T_x(\T\cdot
x)^0$ denotes the annihilator with respect to $\omega_x$ of $T_x(\T\cdot x)$
in $T_x M$. This is called the \emph{symplectic isotropy representation}.

Using the metric $g_x$, $T_x M$ is an orthogonal direct sum of the subspaces
\begin{equation*}
T_x(\T\cdot x)\cong {\mathfrak t}/{\mathfrak t}_x, \qquad
JT_x(\T\cdot x) \cong ({\mathfrak t}/{\mathfrak t}_x)^*
\cong \mathfrak t_x^0, \qquad V_x
\end{equation*}
where ${\mathfrak t}_x^0$ the annihilator of ${\mathfrak t}_x$ in ${\mathfrak
t}^*$ (identified with $JT_x(\T\cdot x)$ using $\omega_x$).

\begin{lemma} \label{sslice} Let $(M,g,J,\omega)$ be an almost K\"ahler
manifold with an isometric hamiltonian $\T$-action. Fix $x\in M$ and a
splitting $\chi\colon{\mathfrak t}\to {\mathfrak t}_x$ of the inclusion.

Then the action of $\T_x$ on the symplectic isotropy representation $V_x$ is
effective, and there is a symplectic form $\omega_0$ on the normal bundle
$N=\T\times_{\T_x} ({\mathfrak t}_x^0\oplus V_x)\to\T\cdot x$ and a
symplectomorphism $f$ from a neighbourhood of the zero section $0_N$ in $N$ to
a neighbourhood of $\T\cdot x$ in $M$ such that\textup:
\begin{bulletlist}
\item the obvious $\T$-action on $N$ by left multiplication is hamiltonian
with momentum map $\mu_0([\alpha,v])= \alpha+\mu_V(v)\circ\chi$, where
$\alpha$ is an element of the fibre belonging to ${\mathfrak t}_x^0 \subset
{\mathfrak t}^*$, and $\mu_V\colon V_x\to{\mathfrak t}_x^*$ is the momentum
map of the symplectic isotropy representation\textup;
\item $f$ is $\T$-equivariant, is equal to the bundle projection along $0_N$,
and its fibre derivative along $0_N$ is the natural identification of the
vertical bundle with $N$.
\end{bulletlist}
\end{lemma}
\begin{proof} The normal exponential map provides a $\T$-equivariant
diffeomorphism from a neighbourhood of $0_N\cong \T\cdot x$ of the normal
bundle $N=\T\times_{\T_x} ({\mathfrak t}_x^0\oplus V_x)\to\T\cdot x$ (with the
natural $\T$ action induced by left multiplication on $\T$) to a neighbourhood
of $\T\cdot x$ in $M$; then $\T$ acts effectively on $N$ while $\T_x$ acts
trivially on ${\mathfrak t}_x^0$, so $\T_x$ acts effectively on $V_x$.

The chosen projection $\chi\colon{\mathfrak t}\to {\mathfrak t}_x$ identifies
the normal bundle $N$ with the symplectic quotient of $T^*\T\times V_x$, by
the diagonal action of $\T_x$ (since $T^*\T\cong \T\times {\mathfrak t}^*$).
The induced symplectic form is $\T$-invariant with the given momentum map.

The pullback of $\omega$ by the normal exponential map gives another
symplectic form $\omega_1$ on a neighbourhood of $0_N$ in $N$, agreeing with
$\omega_0$ along $0_N$ ($\omega_1$ and $\omega_0$ both equal $\omega_x$ at
$T_{(x,0)} N\cong T_xM$).  By the equivariant relative Darboux theorem, there
is a $\T$-equivariant diffeomorphism $h$ of $N$ fixing $0_N$, with $dh=\Id$
there, and such that $h^*\omega_1=\omega_0$ on a neighbourhood $U$ of $0_N$ in
$N$. Then $f=\exp\circ h$ is the equivariant symplectomorphism we seek.
\end{proof}

This result easily generalizes to orbifolds---see~\cite[Lemma~3.5
and Remark~3.7]{LT}.

\subsection{Toric manifolds and orbifolds}\label{s:toric}

A connected $2m$-dimensional symplectic manifold or orbifold $(M,\omega)$ is
said to be {\it toric} if it is equipped with an effective hamiltonian
action of an $m$-torus $\T$ with momentum map $\mu\colon M\to{\mathfrak
t}^*$.  Compact toric symplectic manifolds were classified by
Delzant~\cite{delzant}, and this classification was extended to orbifolds by
Lerman--Tolman~\cite{LT}. Essentially, they are classified by the image
of the momentum map $\mu$, which is a compact convex polytope in $\mathfrak
t^*$, but this statement requires some interpretation, particularly in
the orbifold case.

\begin{defn} Let $\mathfrak t$ be an $m$-dimensional  real vector space.
Then a {\it rational Delzant polytope} $(\Delta,\Lambda,u_1,\ldots u_n)$ in
$\mathfrak t^*$ is a compact convex polytope $\Delta\subset \mathfrak t^*$
equipped with {\it normals} belonging to a lattice $\Lambda$ in $\mathfrak t$
\begin{equation}\label{lattice}
u_j\in\Lambda\subset\mathfrak t
\end{equation}
$(j=1,\ldots n,\ n>m)$ such that
\begin{gather}\label{rdp}
\Delta=\{x\in\mathfrak t^*: L_j(x)\geq 0,\ j=1,\ldots n\}\\
\tag*{with}
L_j(x)=\langle u_j,x\rangle+\lambda_j
\end{gather}
for some $\lambda_1,\ldots \lambda_n\in\R$, and such that for any vertex
$x\in\Delta$, the $u_j$ with $L_j(x)=0$ form a basis for $\mathfrak t$. If
the normals form a basis for $\Lambda$ at each vertex, then $\Delta$ is said
to be {\it integral}, or simply a {\it Delzant polytope}.
\end{defn}
The term {\it rational} refers to the fact that the normals span an
$m$-dimensional vector space over $\Q$. A rational Delzant polytope is
obviously $m$-valent, i.e., $m$ codimension one faces and $m$ edges meet at
each vertex: by~\eqref{rdp} the codimension one faces $\Fa_1,\ldots \Fa_n$ are
given by $\Fa_j=\Delta\cap\{x\in\mathfrak t^*: L_j(x)=0\}$, so that $u_j$ is
an inward normal vector to $\Fa_j$. In the integral case, the $u_j$ are
necessarily primitive, and so are uniquely determined by $(\Delta,\Lambda)$.
In general, the primitive inward normals are $u_j/m_j$ for some positive
integer labelling $m_j$ of the codimension one faces $\Fa_j$, so rational
Delzant polytopes are also called {\it labelled polytopes}~\cite{LT}. However,
it turns out to be more convenient to encode the labelling in the normals.
Note that $\lambda_1,\ldots \lambda_n$ are uniquely determined by
$(\Delta,\Lambda,u_1,\ldots u_n)$.

The rational Delzant theorem~\cite{delzant,LT} states that compact toric
symplectic orbifolds are classified (up to equivariant symplectomorphism) by
rational Delzant polytopes (with manifolds corresponding to integral Delzant
polytopes). Given such a polytope, $(M,\omega)$ is obtained as a symplectic
quotient of $\C^n$ by an $(n-m)$-dimensional subgroup $G$ of the standard
$n$-torus $(S^1)^n=\R^n/2\pi\Z^n$: precisely, $G$ is the kernel of the map
$(S^1)^n\to \T=\mathfrak t/2\pi\Lambda$ induced by the map $(x_1,\ldots
x_n)\mapsto \sum_{j=1}^n x_j u_j$ from $\R^n$ to $\mathfrak t$, and the
momentum level for the symplectic quotient is the image in $\mathfrak g^*$ of
$(\lambda_1,\ldots \lambda_n)\in \R^{n*}$ under the transpose of the natural
inclusion of the Lie algebra $\mathfrak g$ in $\R^n$.

Conversely, a toric symplectic orbifold gives rise to a rational Delzant
polytope as the image $\Delta$ of its momentum map $\mu$, where $\Lambda$ is
the lattice of circle subgroups, and the positive integer labelling $m_j$ of
the codimension one faces $\Fa_j$ is determined by the fact that the local
uniformizing group of every point in $\mu^{-1}(\Fa_j^0)$ is $\Z/m_j\Z$.
(Here, and elsewhere, for any face $\Fa$, we denote by $\Fa^{\smash 0}$ its
interior.)

\begin{rem}\label{cover} Toric symplectic manifolds and orbifolds are
simply connected (as topological spaces---the inverse image of the union of
the faces meeting a given vertex is contractible, and the complement has
codimension two).  However one can consider orbifold coverings and quotients:
a compact convex polytope with chosen normals (giving a basis for $\mathfrak
t$ at each vertex) is a rational Delzant polytope with respect to {\it any}
lattice satisfying~\eqref{lattice}. In particular, if $\Lambda$ is a (finite
index) sublattice of $\Lambda'$, then the torus $\T'=\mathfrak t/2\pi\Lambda'$
is the quotient of $\T=\mathfrak t/2\pi\Lambda$ by a finite abelian group
$\Gamma\cong \Lambda'/\Lambda$.  The corresponding toric symplectic orbifolds
$M$ and $M'$ (under the tori $\T$ and $\T'$) are related by a regular orbifold
covering: $M'=M/\Gamma$.

Clearly there is a `smallest' lattice $\Lambda$ satisfying~\eqref{lattice},
namely the lattice generated by the normals $u_1,\ldots u_n$.  This is a
sublattice of any other lattice $\Lambda'$ with $u_j\in\Lambda'$, so any toric
symplectic orbifold $M'$, corresponding to such a $\Lambda'$, is a quotient of
the toric symplectic orbifold $M$ (corresponding to $\Lambda$) by a finite
abelian group $\Gamma$.

In fact $M$ is the universal orbifold cover of $M'$ in the sense
of~\cite{thurston}. One may also characterize $\Lambda$ as the unique lattice
containing $u_1,\ldots u_n$ for which $G$ is connected, i.e., $M$ is a
symplectic quotient of $\C^n$ by a $(n-m)$-sub{\it torus} of $(S^1)^n$.
\end{rem}

\subsection{Compatible K\"ahler metrics: local theory}\label{s:compkahler}

We turn now to the study of compatible K\"ahler metrics on toric symplectic
orbifolds.  On the union $M^0:=\mu^{-1}(\Delta^0)$ of the generic orbits, such
metrics have an explicit description due to
Guillemin~\cite{guillemin,guillemin0}.  Orthogonal to the orbits is a rank $m$
distribution spanned by commuting holomorphic vector fields $JX_\xi$ for
$\xi\in\mathfrak t$. Hence there is a function $\ang\colon M^0\to \mathfrak
t/2\pi \Lambda$, defined up to an additive constant, such that
$d\ang(JX_\xi)=0$ and $d\ang(X_\xi)=\xi$ for $\xi\in\mathfrak t$. The
components of $\ang$ are `angular variables', complementary to the components
of the momentum map $\mu\colon M^0\to \mathfrak t^*$, and the symplectic form
in these coordinates is simply
\begin{equation}\label{toricomega}
\omega=\langle d\mu\wedge d\ang\rangle,
\end{equation}
where the angle brackets denote contraction of $\mathfrak t$ and $\mathfrak
t^*$.

These coordinates identify each tangent space with $\mathfrak t\oplus
\mathfrak t^*$, so any $\T$-invariant $\omega$-compatible almost K\"ahler
metric is given by
\begin{equation}\label{toricmetric}
g=\langle d\mu, {\bf G} , d\mu \rangle+ \langle d\ang,
{\bf H}, d\ang\rangle,
\end{equation}
where ${\bf G}$ is a positive definite $S^2\mathfrak t$-valued function on
$\Delta^0$, ${\bf H}$ is its inverse in $S^2\mathfrak t^*$---observe that
${\bf G}$ and ${\bf H}$ define mutually inverse linear maps $\mathfrak t^*
\to\mathfrak t$ and $\mathfrak t\to\mathfrak t^*$ at each point---and
$\langle\cdot,\cdot,\cdot\rangle$ denotes the pointwise contraction
${\mathfrak t}^* \times S^2\mathfrak t \times {\mathfrak t}^* \to \R$ or the
dual contraction.  The corresponding almost complex structure is defined by
\begin{equation}\label{toricJ}
Jd\ang = -\langle {\bf G}, d\mu\rangle
\end{equation}
from which it follows that $J$ is integrable if and only if ${\bf G}$ is the
hessian of a function on $\Delta^0$~\cite{guillemin}. 

\begin{rem}\label{angular-coordinates} The description of $\T$-invariant
$\omega$-compatible K\"ahler metrics on $M^0$ shows that they are
parameterized by functions on $\Delta^0$ with positive definite hessian.
There is a subtle point here, however, which is often overlooked in the
literature, namely that the angular coordinates $\ang$ depend on the
(lagrangian) orthogonal distribution to the $\T$-orbits in $M^0$, and there is
no reason for two metrics to have the same orthogonal distribution. This is
not a problem on $M^0$, since the obvious map sending one set of angular
coordinates to another is an equivariant symplectomorphism, but this
symplectomorphism may not extend to $M$.
\end{rem}

The Delzant construction realizes $(M,\omega)$ as a symplectic quotient of
$\C^n$, so there is an obvious choice of a `canonical' compatible K\"ahler
metric $g_0$, namely the one induced by the flat metric on $\C^n$.  An
explicit formula for this K\"ahler metric in symplectic coordinates was
obtained by Guillemin~\cite{guillemin}, and extended to the orbifold case by
Abreu~\cite{Abreu1}: on $M^0$, the canonical metric is given
by~\eqref{toricmetric} with ${\bf G}$ equal to
\begin{equation}\label{canonicalG}
\frac 12 \Hess\biggl(\sum_{j=1}^n L_j(\mu)\log|L_j(\mu)|\biggr)
=\frac 12 \sum_{j=1}^n\frac{u_j\otimes u_j}{L_j(\mu)}.
\end{equation}
Hence the induced metric on $\Delta^0$ is $\frac 12\sum_{j=1}^n
d(L_j)^2/L_j$. (See also~\cite{CDG}.)

\subsection{Compatible K\"ahler metrics: compactification}\label{s:kahlercomp}

On any compact toric symplectic manifold or orbifold, the canonical metric
$g_0$ is globally defined on $M$---by construction. The study of other
globally defined K\"ahler metrics is greatly facilitated by the following
elementary lemma (see also \cite{Abreu} and Remark~\ref{r:abreu}(ii) below).

\begin{lemma}\label{lem-abreu} Let $(M,\omega)$ be a toric symplectic
${2m}$-manifold or orbifold with momentum map $\mu\colon
M\to\Delta\subset\mathfrak t^*$, and suppose that $(g_0,J_0)$, $(g,J)$ are
compatible almost K\"ahler metrics on $M^0=\mu^{-1}(\Delta^0)$ of the
form~\eqref{toricmetric}--\eqref{toricJ}, given by ${\bf G_0},{\bf G}$ and the
same angular coordinates, and such that $(g_0,J_0)$ extends to an almost
K\"ahler metric on $M$. Then $(g,J)$ extends to an almost K\"ahler metric on
$M$ provided that
\begin{align}\label{abreu-bound-a}
{\bf G}-{\bf G_0} \quad&\text{is smooth on $\Delta,$}\\
{\bf G_0}{\bf H}{\bf G_0}-{\bf G_0} \quad&\text{is smooth on $\Delta$.}
\label{abreu-bound-b}
\end{align}
\end{lemma}

\begin{rem} \label{r:abreu} \begin{numlist}
\item We use here the fact that any $\T$-invariant smooth function on
$M$ is the pullback by $\mu$ of a smooth function on $\Delta$ (this follows
from the symplectic slice theorem and~\cite{schwarz}: see~\cite{LT}).
\item \label{p:abreu} For generators $X_\xi,X_\eta$ of the $\T$-action,
$g_0(X_\xi,X_\eta)$ is a $\T$-invariant smooth function on $M$, hence the
pullback of a smooth function on $\Delta$. Thus ${\bf H_0}$ is a smooth
$S^2\mathfrak t^*$-valued function on $\Delta$ (degenerating on
$\partial\Delta$).  Condition~\eqref{abreu-bound-a} thus implies that ${\bf
H_0}{\bf G}$ is smooth on $\Delta$. We claim that in the presence
of~\eqref{abreu-bound-a}, \eqref{abreu-bound-b} is equivalent to ${\bf
H_0}{\bf G}$ being nondegenerate on $\Delta$. Indeed, if ${\bf H_0}{\bf G}$ is
nondegenerate, its inverse ${\bf H}{\bf G_0}$ is smooth on $\Delta$; now
composing ${\bf G}-{\bf G_0}$ on the right by this we
obtain~\eqref{abreu-bound-b}.  Conversely, multiplying by ${\bf H_0}$ we
deduce from \eqref{abreu-bound-b} that ${\bf H}{\bf G_0}$ is smooth on
$\Delta$, so ${\bf H_0}{\bf G}$ is nondegenerate.
\end{numlist}
\end{rem}

\begin{proof}[Proof of Lemma \textup{\ref{lem-abreu}}]
The key point is that it suffices to show $g$ is smooth on $M$: it will then
be nondegenerate because it is compatible with $\omega$ (equivalently if $J$
extends smoothly to $M$, it is an almost complex structure on $M$ by
continuity).  For the smoothness of $g$, we simply compute the difference
\begin{align*}
g-g_0 &= \langle d\mu, {\bf G}-{\bf G_0},d\mu\rangle
+ \langle d\ang, {\bf H}-{\bf H_0},
d\ang\rangle\\
&=\langle d\mu, {\bf G}-{\bf G_0},d\mu\rangle
+ \langle J_0d\mu , {\bf G_0}{\bf H}{\bf G_0}-{\bf G_0},
J_0d\mu \rangle.
\end{align*}
Now $\mu$, $g_0$ and $J_0$ are smooth on $M$, hence so is $g$
by~\eqref{abreu-bound-a}--\eqref{abreu-bound-b}.
\end{proof}

According to Abreu~\cite{Abreu,Abreu1}, when $g_0$ is the canonical metric on
$(M,\omega)$, these conditions are not only sufficient but necessary for the
compactification of $g$. However, in our view there are some shortcomings in
his (rather sketchy) proof. In particular he does not address the issue of the
dependence of the angular coordinates on the metric. The following observation
only partially resolves this difficulty.

\begin{lemma}\label{lem:sameG} Let $(M,\omega)$ be a compact toric symplectic
manifold with two compatible almost K\"ahler metrics which induce the same
$S^2\mathfrak{t}$-valued function ${\bf G}$ on the interior of the Delzant
polytope. Then there is an equivariant symplectomorphism of $M$ sending one
metric to the other.
\end{lemma}
\begin{proof}
By Remark~\ref{angular-coordinates}, such a symplectomorphism exists on
$M^0$. It extends uniquely to $M$, since $M^0$ is dense and $(M,g)$ is a
complete. The extension is a distance isometry by continuity, and is therefore
smooth by a standard argument.
\end{proof}
Note that this lemma makes essential use of the completeness of $(M,g)$.  It
can, however, be extended to compact orbifolds, for instance by lifting the
distance isometry to compatible uniformizing charts.

On the other hand we learn nothing about the dependence of the angular
coordinates on metrics which induce different $S^2\mathfrak{t}$-valued
functions on the interior of the Delzant polytope.  Thus the above lemma does
not suffice to clarify Abreu's proof.

The issue of the compactification of toric K\"ahler metrics is an important
one.  We shall therefore establish precise necessary and sufficient
compactification conditions by a self-contained argument. Our proof also has
the merit of being elementary and, modulo the above lemma, local, in contrast
to~\cite{Abreu,Abreu1}, where the existence of a global biholomorphism is
used. Indeed, compactification is about boundary conditions, so it is a local
question. We shall present these boundary conditions in a form more closely
analogous to the well-known conditions in complex dimension one.  As a warm-up
for the rest of the subsection we first recall this case.

Let $(M,\omega)$ be a compact toric symplectic $2$-orbifold. This must be an
orbifold $2$-sphere (i.e., equivariantly homeomorphic to $\C P^1$ with the
standard circle action, but the two fixed points may be orbifold
singularities), equipped with a rotation invariant area form. On $M^0$, which
is diffeomorphic to $\C^\times$, a compatible K\"ahler metric takes the form
\begin{equation}\label{CP1}
g = \frac{d\mu^2}{\Theta(\mu)} + \Theta(\mu) dt^2,
\end{equation}
where $\omega=d\mu\wedge dt$. The rational Delzant polytope is an interval
$[\alpha,\beta]\in\mathfrak t^*$ with normals $u_\alpha, u_\beta\in\mathfrak
t$.  If we identify a generator of the lattice $\Lambda$ in $\mathfrak t$
with $1\in\R$ (chosen so that $u_\alpha$ is positive), then $t\colon M^0\to
\mathfrak t/2\pi \Lambda$ becomes a coordinate of period $2\pi$, and the
orbifold singularities have cone angles $2\pi/m_\alpha$, $2\pi/m_\beta$ where
$m_\alpha=u_\alpha, m_\beta=-u_\beta\in\Z^+$.

Since $\Theta(\mu)$ is the norm squared of the Killing vector field, $\Theta$
is smooth on $[\alpha,\beta]$, positive on the interior, and zero at the
endpoints. On the other hand, $\mu$ is a Morse function (i.e., the two
critical points are nondegenerate---this follows easily using a symplectic
slice) and $dd^c\mu=\Theta'(\mu)\omega$, so that $\Theta'(\alpha)$ and
$\Theta'(\beta)$ are nonzero.

Now let $\hat U\subset\R^2$ be an orbifold chart covering an $S^1$-invariant
neighbourhood $U = {\hat U}/{\Z}_{m_\alpha}$ of $\mu^{-1}(\alpha)$, where
$\Z_{m_\alpha}$ acts in the standard way on $\R^2$ and the covering map $\pi$
sends $0$ to $\mu^{-1}(\alpha)$. The $S^1$-action on $U$ lifts to one on
${\hat U}$, fixing $0$ and commuting with ${\mathbb Z}_{m_\alpha}$.  Now $\hat
t = t\circ \pi/m_{\alpha}$ is a coordinate of period $2\pi$ on $\hat
U\setminus \{0\}$ while $\hat \mu = m_\alpha (\mu \circ \pi)$ is the momentum
map of the $S^1$ action on $\hat U$, with respect to $\hat \omega = d\hat \mu
\wedge d\hat t = \pi^*\omega$. The pull back of $g$ to $\hat U \setminus
\{0\}$ is
\begin{equation*}
{\hat g} = \frac{d{\hat \mu}^2}{m_{\alpha}^2\Theta(\hat \mu/m_{\alpha})}
+ m_{\alpha}^2\Theta(\hat \mu/m_{\alpha})d{\hat t}^2.
\end{equation*}
If this metric compactifies smoothly at $0$ we must have
$m_\alpha\Theta'(\alpha)=2$ (see \cite{hwang-singer}). With an analogous
argument at $\mu^{-1}(\beta)$, we deduce that $u_\alpha\Theta'(\alpha)=2
=u_\beta\Theta'(\beta)$.

To show that these conditions are sufficient for the smooth extension of $g$
(in the orbifold sense) to $M$, we put $r^2/2=\mu-\alpha$ and let $t$ have
period $2\pi/m_\alpha$. Since $\Theta(\alpha)=0$, $g$ differs from a multiple
of $g_0=dr^2+\frac14 \Theta'(\alpha)^2 r^2 dt^2$ by a smooth bilinear form on
$M$, vanishing at $\mu=\alpha$. Clearly the condition
$\Theta'(\alpha)=2/m_\alpha$ provides a smooth (orbifold) extension of $g_0$
to $\mu^{-1}(\alpha)$ by considering $(r,t/m_\alpha)$ to be the polar
coordinates in a uniformising chart. The other endpoint is analogous.

To summarize, $g$ given by~\eqref{CP1} is globally defined on a toric orbifold
whose rational Delzant polytope is $[\alpha,\beta]\subset\mathfrak t^*$, with
normals $u_\alpha, u_\beta\in\mathfrak t$, if and only if $\Theta$ smooth on
$[\alpha,\beta]$, with
\begin{equation}\label{ell=1}\begin{split}
\Theta(\alpha) &= 0 = \Theta(\beta), \\
\Theta'(\alpha)u_\alpha &= 2 = \Theta'(\beta)u_\beta
\end{split}\end{equation}
and $\Theta$ positive on $(\alpha,\beta)$.  The derivative conditions make
invariant sense, since $\Theta$ takes values in $(\mathfrak t^*)^2$, so its
derivative takes values in $\mathfrak t^*$. Also note that the conditions are
manifestly independent of the choice of lattice (as they should be).

\smallbreak

In order to generalize this criterion to the case $m>1$, we introduce some
notation. For any face $\Fa\subset \Delta$, we denote by ${\mathfrak t}_\Fa
\subset {\mathfrak t}$ the vector subspace spanned by the inward normals $u_j
\in {\mathfrak t}$ to all codimension one faces of $\Delta$, containing $\Fa$;
thus the codimension of ${\mathfrak t}_\Fa$ equals the dimension of
$\Fa$. Furthermore, the annihilator ${\mathfrak t}^0_\Fa$ of ${\mathfrak
t}_\Fa$ in ${\mathfrak t}^*$ is naturally identified with $({\mathfrak
t}/{\mathfrak t}_\Fa)^*$.

\begin{prop}\label{first-order-conditions} Let $(M,\omega)$ be a compact
toric symplectic $2m$-manifold or orbifold with momentum map $\mu \colon M\to
\Delta \subset {\mathfrak t}^*$ and ${\bf H}$ be a positive definite
$S^2{\mathfrak t}^*$-valued function on $\Delta^0$.  Then ${\bf H}$ comes from
a $\T$-invariant, $\omega$-compatible almost K\"ahler metric $g$ via
\eqref{toricmetric} if and only if it satisfies the following
conditions\textup:
\begin{bulletlist}
\item \textup{[smoothness]} ${\bf H}$ is the restriction to $\Delta^0$ of a
smooth $S^2{\mathfrak t}^*$-valued function on $\Delta$\textup;
\item \textup{[boundary values]} for any point $y$ on the codimension one face
$\Fa_j \subset \Delta$ with inward normal $u_j$, we have
\begin{equation}\label{eq:toricboundary}
{\bf H}_y(u_j, \cdot) =0\qquad {and}\qquad (d{\bf H})_y(u_j,u_j) = 2 u_j,
\end{equation}
where the differential $d{\bf H}$ is viewed as a smooth $S^2{\mathfrak
t}^*\otimes {\mathfrak t}$-valued function on $\Delta$\textup;
\item \textup{[positivity]} for any point $y$ in interior of a face
$\Fa \subseteq \Delta$, ${\bf H}_y(\cdot, \cdot)$ is positive definite
when viewed as a smooth function with values in $S^2({\mathfrak
t}/{\mathfrak t}_\Fa)^*$.
\end{bulletlist}
\end{prop}
\begin{proof} We first prove the necessity of these conditions.
Let $(M,\omega,\mu)$ be a compact toric symplectic orbifold with polytope
$\Delta$, and $(g,J)$ a compatible K\"ahler metric. For any $x \in M$ and
$\xi,\eta \in \mathfrak t$, we put ${\bf H}_{\mu(x)} (\xi,\eta) = g_x (X_\xi,
X_\eta)$. Clearly ${\bf H}$ is an $S^2{\mathfrak t}^*$-valued function on
$\Delta$ and the smoothness and positivity properties follow immediately from
the definition.

It remains to establish the boundary values~\eqref{eq:toricboundary} for
$y=\mu(x)$ in a codimension one face $F_j$. The vanishing of ${\bf H}_y(u_j,
\cdot) =0$ is immediate from the definition (the Killing vector field
corresponding to $u_j$ vanishes on $\mu^{-1}(F_j)$). This implies in
particular that $d{\bf H}_y(u_j, u_j)$ is proportional to $u_j$. To obtain the
correct constant, we use a symplectic slice, as in Lemma~\ref{sslice}, to
pullback the metric $g$ to the normal bundle $N$ of the orbit $\T\cdot x$ for
a point $x\in M$ with $\mu(x)=y$, and restrict to the symplectic isotropy
representation $V_x$.  By construction, the Killing vector field corresponding
to $u_j$ induces the generator $X$ of the standard circle action on $V_x$, and
the metric induced by $g$ agrees to first order at $0$ with the constant
metric $g_0$ given by $g_x$. It is now straightforward to check that the
constant is $2$ (indeed, $(V_x,g_0,\omega_0)$ is a toric K\"ahler
$2$-orbifold, so we have already computed this above).

\smallbreak

Now we explain why the given conditions are sufficient to conclude ${\bf H}$
that comes from a smooth compatible metric on $(M,\omega)$.

We know that the function ${\bf H_0}= {\bf G_0}^{-1}$, with ${\bf G_0}$
defined by~\eqref{canonicalG}, does correspond to a globally defined invariant
K\"ahler metric on $(M,\omega)$ (and so it satisfies the given conditions, as
one can easily check directly). By virtue of Lemma~\ref{lem-abreu}, it is
enough to show that for any ${\bf H}={\bf G}^{-1}$ satisfying the given
conditions, the sufficient conditions
\eqref{abreu-bound-a}--\eqref{abreu-bound-b} are satisfied. As explained in
Remark~\ref{r:abreu}, we have to check that both ${\bf H} {\bf G_0}$ and ${\bf
G} - {\bf G_0}$ are smoothly extendable about each point $y_0 \in \partial
\Delta$. We shall establish this by a straightforward argument using Taylor's
Theorem.

Suppose that $y_0$ belongs to the interior of a $k$-dimensional face $\Fa$ of
$\Delta$. Let us choose a vertex of $\Fa$. Since $\Delta$ is a rational
Delzant polytope, the affine functions $L_i(y)=\langle u_i,y\rangle+\lambda_i$
which vanish at this vertex form a coordinate system on $\Delta$.  By
reordering the inward normals $u_1,\ldots u_n$, we can suppose that these
coordinate functions are $L_1(y), \ldots L_m(y)$ (so $u_1,\ldots u_m$ form a
basis for $\mathfrak t$) and that $L_1(y),\ldots L_{m-k}(y)$ vanish on $\Fa$
(so $u_1,\ldots u_{m-k}$ span ${\mathfrak t}_\Fa$). We set
$y_i=L_i(y)-L_i(y_0)$ for $i=1,\ldots m$.  These functions also form a
coordinate system on $\Delta$, with $y_0$ corresponding to the origin, and
$y_1,\ldots y_{m-k}$ vanish on $\Fa$.

We now let $H_{ij}(y) = {\bf H}_y(u_i,u_j)$ and let $(G_{ij}(y))$ be the
inverse matrix to $(H_{ij}(y))$ (which is the matrix of ${\bf G}$ with respect
to the dual basis). Similarly we define inverse matrices $(H^0_{ij}(y))$ and
$(G^0_{ij}(y))$. The conditions (i)--(iii) imply:
\begin{bulletlist}
\item $H_{ij}(y)$ are smooth functions on $\Delta$;
\item on any codimension one face $\Fa_i$ containing $\Fa$ (with inward normal
$u_i$, $i=1,\ldots m-k$), we have
\begin{equation}\label{eq:explicittoricbdy}
H_{ij}(y) = H_{ji}(y) = 0 \qquad \text{for all $j=1,\ldots m$}
\qquad \text{and} \qquad \partial H_{ii}/\partial y_i = 2.
\end{equation}
\item the matrix $(H_{ij}(y))_{i,j=m-k+1}^m$ is positive definite
on the interior of $\Fa$;
\end{bulletlist}
We conclude from~\eqref{eq:explicittoricbdy} that for $i=1,\ldots m-k$,
$H_{ij}(y) =H_{ji}(y) = O(y_i)$ (for all $j=1,\ldots m$) and
$H_{ii}(y)=2y_i(1+O(y_i))$, where $O(y_i)$ denotes the product of $y_i$ with a
smooth function of $y$.

Putting these conditions together, we then have:
\begin{align*}
H_{ij}(y) &= 2y_i\delta_{ij}+y_i y_j F_{ij}(y)
&&\text{for }i,j=1,\ldots m-k \\
H_{ij}(y)&= y_i F_{ij}(y) && \text{for } i =1,\ldots m-k
\text{ and } j=m-k+1,\ldots m,
\end{align*}
where $F_{ij}$ are smooth functions. (Recall also that $H_{ij}=H_{ji}$.)

It follows that $\det (H_{ij}(y)) = 2^{m-k} y_1y_2\cdots y_{m-k} P(y)$ where
the function $P(y)=\det (H_{ij}(y))_{i,j=m-k+1}^m + O(y_1) + O(y_2)+\cdots
O(y_{m-k})$ is positive at the origin. Since the same holds for $H^0_{ij}(y)$
it follows that $\det (H_{ij}(y))/\det (H^0_{ij}(y))$ can be extended to the
origin as a smooth and positive function.

On the other hand $G_{pq}(y)$ is the determinant of a cofactor matrix of
$(H_{ij}(y))$ divided by $\det(H_{ij}(y))$. This will be smooth if the
determinant of the cofactor is $O(y_i)$ for each $i=1,\ldots m-k$.  We see
that this is true unless $1\leq p=q\leq m-k$, in which case we obtain
$G_{pp}(y) = (1+O(y_p))/2y_p$. The same holds for $G^0_{pq}(y)$.

We deduce that ${\bf G} - {\bf G_0}$ is smooth at $y_0$, and hence
${\bf H_0}{\bf G}$ is smooth at $y_0$. Since it is nondegenerate there,
its inverse ${\bf H}{\bf G_0}$ is also smooth.
\end{proof}

\begin{rem} \label{abreuok} By continuity, it suffices that the boundary
conditions~\eqref{eq:toricboundary} hold on the interior of the codimension
one faces. However they and their tangential derivatives imply that for a
point $y$ on \emph{any} face $\Fa \subset \Delta$, we have
\begin{equation}\label{eq:strongtoricboundary}
{\bf H}_y(u_j, \cdot) =0\qquad {and}\qquad (d{\bf H})_y(u_j,u_k) = 2\delta_{jk}
u_j
\end{equation}
for any inward normals $u_j,u_k$ in ${\mathfrak t}_F$.

The proof also shows that our first order conditions are equivalent to
\eqref{abreu-bound-a}--\eqref{abreu-bound-b} with ${\bf G}_0$ given by
\eqref{canonicalG}, thus establishing the validity of \cite{Abreu1}---see
Remark~\ref{r:abreu}(ii).
\end{rem}

\subsection{Toric complex manifolds and orbifolds}\label{s:complextoric}

We now turn briefly to the complex point of view on toric K\"ahler manifolds
and orbifolds.  Given a rational Delzant polytope $(\Delta,\Lambda,u_1,\ldots
u_n)$, we obtain a complex subgroup $G^c$ of $(\C^\times)^n$ as the
complexification of $G$. The relation between complex quotients and symplectic
quotients then shows~\cite{audin,guillemin,kirwan} that the canonical complex
structure on the toric symplectic orbifold $(M,\omega)$ constructed from
$\Delta$ is equivariantly biholomorphic to the quotient by $G^c$ of a dense
open subset $\C^n_s$ of $\C^n$ given by
\begin{equation}\label{complexquot}
\C^n_s={\textstyle\bigcup_\Fa \C^n_\Fa},\quad
\C^n_\Fa = \{(z_1,\ldots z_n)\in\C^n: z_j=0 \text{ iff } L_j(x)=0
\text{ for } x\in \Fa^0\}.
\end{equation}
Thus $\C^n_s$ is $\C^n$ with the coordinate subspaces removed that do not
correspond to faces of $\Delta$. Observe that the complex quotient only
depends on the inward normals (which determine $G^c$) and the combinatorics of
the faces (which determine $\C^n_s$), i.e., by specifying which sets of
codimension one faces have nonempty intersection.  These data can be encoded
in a family of convex simplicial cones called a \emph{fan}.

Furthermore, {\it any} $\T$-invariant $\omega$-compatible complex structure on
$M$ is equivariantly biholomorphic to the standard one (see~\cite{LT} for the
result in the general orbifold case). Of course this biholomorphism does not
preserve $\omega$ in general. Thus two toric K\"ahler manifolds (or orbifolds)
are equivariantly biholomorphic if and only if they have the same fan.

\subsection{Restricted toric manifolds}\label{s:restricttoric}

Toric K\"ahler manifolds can be used to provide examples of K\"ahler manifolds
with non-toric isometric hamiltonian torus actions simply by restricting the
action to a subtorus. These torus actions can be surprisingly complicated in
general. However, the subtori generated by a subset of the normals to the
Delzant polytope have much simpler actions.

\begin{ex}\label{ex:cp2} We can illustrate this in the simplest
nontrivial case of $S^1$ actions on $\C P^2$, which is toric under the action
of $\T\cong S^1\times S^1$ given by $(\lambda_1,\lambda_2)\colon
[z_0,z_1,z_2]\mapsto [z_0,\lambda_1 z_1,\lambda_2 z_2]$. The `tame' $S^1$
subgroups generated by the normals are given by $\lambda_1=0$, $\lambda_2=0$
or $\lambda_1=\lambda_2$. The momentum map of the $S^1$ action is then the
projection of the momentum map of $\T$ along the corresponding face of the
Delzant polytope $\Delta$ (which is a simplex). The momentum map of `wild'
$S^1$ subgroups, such as $\lambda_1=\lambda^2$, $\lambda_2=\lambda^3$, are
given by more general projections. We wish to draw attention to two
distinctions between these two types of $S^1$ action.
\begin{numlist} 
\item For tame actions, the momentum map of the $S^1$ action has no critical
values on the interior of the momentum interval, whereas for wild actions it
does.
\item For tame actions, the orbits of the complexified action (of $\C^\times$)
have smooth closures, whereas for wild actions, they do not---for instance
they are singular cubics for the case $\lambda_1=\lambda^2$,
$\lambda_2=\lambda^3$
\end{numlist}

\begin{figure}[ht]
\begin{center}
\includegraphics[width=.5\textwidth]{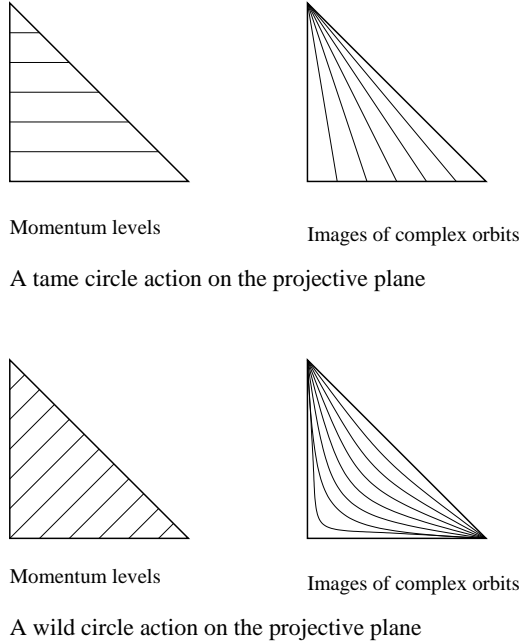}
\caption{Circle actions on $\C P^2$}
\end{center}
\end{figure}

The blow up of $\C P^2$ at a point is the first Hirzebruch surface
$F_1=P(\cO\oplus \cO(1))\to \C P^1$. If this point is one of the three fixed
points of the $\T$-action corresponding to a vertex of $\Delta$, then the
standard fibrewise $S^1$ action on $F_1$ descends to the tame $S^1$ action on
$\C P^2$ corresponding to the opposite edge. Thus a tame $S^1$ action realises
$\C P^2$ as the blowdown of a toric bundle (of projective lines) over $\C
P^1$.
\end{ex}

We generalize this by considering torus actions on blowdowns of toric bundles
(with fibre any toric K\"ahler manifold) over a product of complex projective
spaces.

Let $\cV$ be a toric K\"ahler $2\ell$-manifold, under a torus $\T$, with
Delzant polytope $\Delta$.  By the Delzant construction $\cV$ is
($\T$-equivariantly symplectomorphic to) a symplectic quotient of $\C^n$ by an
$n-\ell$ dimensional subgroup $G$ of the standard $n$-torus $\T^n$ (with
$\T=\T^n/G$). From \S\ref{s:complextoric}, ${\cV}$ is also
($\T^c$-equivariantly biholomorphic to) the holomorphic quotient
${\C}^n_s/G^c$, where ${\C}^n_s$ is the set of stable points in $\C^n$.

Given integers $d_1,\ldots d_n\geq 0$, there are now two constructions we can
make.
\begin{numlist}
\item Let $\C^D=\C^{d_1+1}\times \C^{d_2+1}\times\cdots \C^{d_n+1}$.  Then we
have a block diagonal action of $\T^n$ on $\C^D$ as a subtorus of the standard
torus $\T^N$: the $i$th circle in $\T^n$ acts by scalar multiplication on
$\C^{d_i+1}$ and trivially on the other factors. Since $G$ is a subtorus of
$\T^n$, we can form the symplectic quotient of $\C^D$ by $G$ and this will be
diffeomorphic to the stable quotient by $G^c$. Let us denote the corresponding
manifold by $M$.

The standard K\"ahler structure on $\C^D$ can be written in block diagonal
momentum coordinates $(x_1,\ldots x_n)$ of $\T^n$ as
\begin{align*}
\tilde g_0
&=\sum_{j=1}^n\Bigl(x_j g_j+ \frac{dx_j^2}{2x_j}+2x_j\theta_j^2\Bigr)\\
\tilde \omega
&= \sum_{j=1}^n \bigl(x_j\omega_j+dx_j\wedge\theta_j\bigr),
\qquad d\theta_j=\omega_j
\end{align*}
where $x_j=r_j^2/2$, for the radial coordinate $r_j$ on $\C^{d_j+1}$, and
$g_j$ is the Fubini--Study metric on $\C P^{d_j}$, normalized so that
$\theta_j$ is the connection $1$-form of the Hopf fibration and $\frac12
g_j+\theta_j^2$ is the round metric on the unit sphere $S^{2d_j+1}$: we obtain
the flat metric in spherical polar coordinates on each $\C^{d_j+1}$ factor by
substituting $x_j=r_j^2/2$.

This induces a K\"ahler structure on $M$ by writing the momentum coordinates
$x_j=L_j(\mu)$ of $\T^n$ in terms of the momentum map $\mu$ of $\T$, where
$L_1,\ldots L_n$ are the affine functions defining $\Delta$.  The resulting
K\"ahler metric, in the notation of \S\ref{s:compkahler}, is
\begin{equation}\label{eq:standard}\begin{split}
g_0' &=
\sum_{j=1}^n L_j(\mu) g_j + \langle d\mu, {\bf G_0}, d\mu \rangle
+ \langle\boldsymbol\theta,{\bf H_0}, \boldsymbol \theta\rangle,\\
\omega' &= \sum_{j=1}^n L_j(\mu)\omega_j
+ \langle d\mu\wedge\boldsymbol\theta\rangle, \qquad 
d\boldsymbol\theta=\sum_{j=1}^n \omega_j\otimes u_j
\end{split}\end{equation}
with ${\bf G_0}$ given by~\eqref{canonicalG} and ${\bf H_0}$ is inverse to
${\bf G_0}$. This reduces to the \emph{canonical} toric K\"ahler structure on
$\cV$ when $d_j=0$ for all $j$.

Our aim is to show that there is a compatible K\"ahler structure on $M$
generalizing the \emph{given} toric K\"ahler structure of $\cV$, which is
determined by an arbitrary matrix ${\bf H}$ satisfying the necessary and
sufficient compactification conditions of
Proposition~\ref{first-order-conditions}. To do this, and to understand better
the holomorphic geometry of $M$, we consider another construction.

\item Let $\tilde M=\bigoplus_{j=1}^n \cO(-1)_j \to \prod_{j=1}^n\C P^{d_j}$,
where $\cO(-1)_j=\cO(0,\ldots 0,-1,0,\ldots 0)$ is the line bundle which is
$\cO(-1)$ over $\C P^{d_j}$ and trivial over the other factors. Let $\tilde
M^0=\prod_{j=1}^n \cO(-1)^\times_j$ be the associated holomorphic principal
$(\C^\times)^n$-bundle (given by removing the zero section from each line
bundle). Now $\tilde M=\tilde M_0\times_{(\C^\times)^n} \C^n$ admits a
fibre-preserving holomorphic action of $(\C^\times)^n$.

Since $G$ is a subtorus of ${\T}^n$ we can form the holomorphic stable
quotient of $\tilde M$ by $G^c$ to obtain a complex manifold $\hat M$. We see
immediately that $\hat M = M^0\times_{\T^c}\cV$ where $M^0=\tilde
M^0/G^c$. Thus $\hat M$ is a bundle of toric complex manifolds.
\end{numlist}

It is easy to see how constructions (i) and (ii) are related, since
$\cO(-1)\to\C P^{d_j}$ is the blow-up of $\C^{d_j+1}$ at the origin, so that
$\tilde M$ is (equivariantly biholomorphic to) the blow-up of $\C^D=\prod_j
\C^{d_j+1}$ along the union over $j$ of the coordinate subspaces with zero in
the $j$th factor. The stable quotients of $\tilde M$ and $\C^N$ that we
consider are related by this blow-up (by construction), and so $\hat M$ is
(equivariantly biholomorphic to) a blow-up of $M$.

The K\"ahler structure~\eqref{eq:standard} on $M$ therefore pulls back to give
a K\"ahler structure on $\smash{\hat M}$, except that the metric and
symplectic form degenerate on the exceptional divisor. Again by construction,
this induces the canonical toric K\"ahler structure of $\cV$ on each fibre of
$\hat M$.

Let ${\bf G}={\bf H}^{-1}$ be the matrices inducing the given toric K\"ahler
structure on $\cV$. Then we obtain a new K\"ahler structure on $\hat M$,
degenerating on the exceptional divisor and inducing the given toric K\"ahler
structure on each fibre:
\begin{equation}\label{eq:given}\begin{split}
g' &=
\sum_{j=1}^n L_j(\mu) g_j + \langle d\mu, {\bf G}, d\mu \rangle
+ \langle\boldsymbol\theta,{\bf H}, \boldsymbol \theta\rangle,\\
\omega' &= \sum_{j=1}^n L_j(\mu)\omega_j
+ \langle d\mu\wedge\boldsymbol\theta\rangle, \qquad 
d\boldsymbol\theta=\sum_{j=1}^n \omega_j\otimes u_j.
\end{split}\end{equation}

There is no reason \textit{a priori} why this should descend to $M$ (in
particular, the complex structure is different). Nevertheless, it does,
because of the strong control over the boundary behaviour of ${\bf H}$ given
by Proposition~\ref{first-order-conditions}.

\begin{prop} The degenerate K\"ahler structure~\eqref{eq:given} on $\hat M$
descends to give a \textup(nondegenerate\textup) K\"ahler structure on $M$.
\end{prop}
\begin{proof}
We know that~\eqref{eq:standard} is globally defined smooth K\"ahler structure
on $M$. We shall show that~\eqref{eq:given} defines a compatible K\"ahler
metric on the same symplectic manifold (with the same angular coordinates).
For this, it suffices to show that the difference $g'-g_0'$ is smooth on
$M$. However, since the compatible K\"ahler metrics defined on $\cV$ by ${\bf
H}$ and ${\bf H_0}$ are smooth, Proposition~\ref{first-order-conditions} and
Remark~\ref{abreuok} show that ${\bf G}-{\bf G_0}$ and ${\bf G_0}{\bf H}{\bf
G_0}-{\bf G_0}$ are smooth functions on the Delzant polytope $\Delta$ of
$\cV$. Now the momentum map $\mu$ on $(M,\omega')$ is smooth, with image
$\Delta$. It therefore follows, as in the proof of Lemma~\ref{lem-abreu}, that
$g'-g_0'$ is smooth.
\end{proof}

\section{Rigid hamiltonian torus actions}\label{s:rigid}

In this section, we introduce the notion of a {\it rigid hamiltonian torus
action}. Toric K\"ahler manifolds automatically carry such an action: our goal
is to extend some of the rigid properties of toric K\"ahler manifolds to rigid
torus actions in general, and to classify them.  In the first three
subsections we study respectively the differential topology, symplectic
geometry and biholomorphism type of compact (smooth) K\"ahler manifolds with
such an action, then we combine these threads to describe the K\"ahler
geometry. In the final subsection, we specialize to the case that the torus
action is `semisimple' and give a generalized Calabi construction of all
compact K\"ahler manifolds with a semisimple rigid torus action.

\subsection{Stratification of the momentum polytope}\label{s:stratification}

Before defining the torus actions we will consider, we establish a couple of
basic facts.  We shall make essential use the convexity theorem of Atiyah
and Guillemin--Sternberg~\cite{atiyah,GS}.

\begin{lemma}\label{tricky}
Let $\T$ be a torus in the group of hamiltonian isometries of a compact
connected K\"ahler manifold $(M,g,J,\omega)$, which is the closure of the
group generated by $\ell$ hamiltonian Killing vector fields
$K_r=J\grad_g\sigma_r$ $(r=1,\ldots\ell)$ that are independent on a dense
open set.  Suppose that $\g{K_r,K_s}$ depends only on
$(\sigma_1,\ldots\sigma_\ell)$ for $r,s=1,\ldots\ell$. Then\textup:
\begin{numlist}
\item the torus $\T$ has dimension $\ell$\textup;
\item the image of the momentum map $\mu\colon M\to{\mathfrak t}^*$ of $\T$
is a compact convex polytope such that $\mu$ is regular
\textup(i.e., submersive\textup) as a map to the interior of any of its
faces.
\end{numlist}
\end{lemma}
\begin{proof}
By the Atiyah--Guillemin--Sternberg convexity theorem~\cite{atiyah,GS}, the
image of $\mu$ is a compact convex polytope $\Delta$ in ${\mathfrak t}^*$, the
convex hull of the finite image $I$ of the fixed point set of $\T$.  The
momentum coordinates $\boldsymbol\sigma = (\sigma_1,\ldots \sigma_{\ell})$ are
related to $\mu$ by the natural inclusion
\begin{equation*}
\R^{\ell} \cong {\rm span}(K_1,\ldots K_{\ell}) \subseteq {\mathfrak t},
\end{equation*}
which in turn gives rise to a linear projection $\pi\colon {\mathfrak t}^*
\to \R^{\ell*}$ such that $\boldsymbol\sigma = \pi \circ \mu$.

Let us first consider the image of $\Delta$ by $\pi$. We claim that $\pi$ is
injective on $I$. Indeed, since $K_1,\ldots K_\ell$ generate $\T$, the fixed
point set is precisely the set of common zeros of $K_1,\ldots K_\ell$, and
since $\g{K_r,K_s}$ depends only on $\boldsymbol\sigma$, the preimage of an
element of $\pi(\Delta)$, containing an element of $I$, consists entirely of
elements of $I$. Now $I$ is finite and the preimages of $\pi$ are convex, so
each such preimage has just one point.

Second, we note that the set of regular values of $\boldsymbol\sigma$ is
connected.  Indeed, the critical point set of $\boldsymbol\sigma$ in $M$ has
codimension at least two---it is the set where the holomorphic $\ell$-vector
$K_1^{1,0}\wedge\cdots \wedge K_\ell^{1,0}$ vanishes---so the set of regular
points $U$ is connected. Now as $\g{K_r,K_s}$ depends only on
$\boldsymbol\sigma$, the inverse image of a critical value consists entirely
of critical points, so the set of regular values is $\boldsymbol\sigma(U)$.

Third, consider the orbits of the commuting vector fields $JK_1,\ldots
JK_\ell$---this is the gradient flow of $\boldsymbol\sigma$, and so the orbit
of any regular point consists entirely of regular points and its boundary
points are all critical. Now regular points map to regular values and critical
points to critical values, so by the connectivity of the regular values, all
regular orbits have the same image---and the closure is the image of
$\boldsymbol\sigma$ since regular values are dense.

These facts implies the conclusions of the lemma as follows.

\begin{numlist}
\item Suppose $x$ is a regular point of $\boldsymbol\sigma$ and $\mu(x)$
belongs to a closed face $\Fa$ of $\Delta$. Then the $\T^c$ orbit of $x$ also
maps to $\Fa$, where $\T^c$ is the complexification of $\T$. Since the orbit
under $JK_1,\ldots JK_\ell$ is contained in the $\T^c$ orbit, $\pi$ maps $\Fa$
onto ${\rm Im}\,\boldsymbol\sigma$. Now $\pi$ is bijective on vertices, so
$\Fa=\Delta$. In other words the inverse image (under $\pi$) of a regular
value of $\boldsymbol\sigma$ meets no proper face of $\Delta$: this clearly
implies $\pi$ is bijective, hence $\mu=\boldsymbol\sigma$ and $\dim \T=\ell$.

\item We have seen that the image of the closure $C$ of any regular $\T^c$
orbit is the whole of $\Delta$. Atiyah~\cite{atiyah} shows that the inverse
image in $C$ of any open face $\Fa^0$ is a single $\T^c$-orbit and $\mu$ is a
submersion from this orbit to $\Fa^0$.  Since this is true for all regular
orbits, and the union of the regular orbits is dense, the claim follows.
\endnumlproof

\begin{defn}\label{rigid} Let $(M,g,J,\omega)$ be a connected K\"ahler
$2m$-manifold with an effective isometric hamiltonian action of an
$\ell$-torus $\T$ with momentum map $\mu\colon M\to{\mathfrak t}^*$.  We say
the action is {\it rigid} iff for all $x\in M$, $R_x^*g$ depends only on
$\mu(x)$, where $R_x\colon \T\to \T\cdot x\subset M$ is the orbit map.
\end{defn}

In other words, for any two generators $X_\xi,X_\eta$ of the action
($\xi,\eta\in\mathfrak t$), $\g{X_\xi,X_\eta}$ is constant on the levels of
the momentum map $\mu$.  We remark that the inverse image of a critical value
of $\mu$ can be approximated (to first order on a dense open subset) by
inverse images of nearby regular values. Hence it suffices to know that the
generators have constant inner products on the generic level sets of
$\mu$. Thus part (i) of Lemma~\ref{tricky} implies that on a compact manifold,
a {\it local} rigid torus action (as in~\cite{ACG1}) is necessarily a global
one. In particular, on a compact K\"ahler manifold with a hamiltonian $2$-form
of order $\ell$, the associated Killing vector fields $K_1,\ldots K_\ell$
generate a rigid $\ell$-torus action. Another example is any toric K\"ahler
manifold.

Part (ii) of Lemma~\ref{tricky} has further consequences for compact
K\"ahler manifolds with a rigid torus action.

\begin{prop}\label{strata} Suppose $(M,g,J,\omega)$ is a
compact connected K\"ahler manifold of dimension $2m$, with a rigid
hamiltonian $\ell$-torus action with momentum map $\mu$ whose image is a
compact convex polytope $\Delta$.
\begin{numlist}
\item If $\Fa$ is a $k$-dimensional closed face $(0\leq k\leq \ell)$ of
$\Delta$, then $M_\Fa:=\mu^{-1}(\Fa)$ is a compact totally geodesic K\"ahler
submanifold of $M$ of dimension $2(m_\Fa+k)$ $(0\leq m_\Fa\leq m-\ell)$ with
a rigid hamiltonian action of a $k$-torus $\T/\T_\Fa$, where $\T_{\Fa}$ is
the intersection of the isotropy subgroups of points in $M_\Fa$.
\item If $\Fa^0$ is the interior of $\Fa$, then
$M_{\Fa}^0:=\mu^{-1}(\Fa^0)\cong\Fa^0\times P_\Fa$ where $P_\Fa$ is a compact
manifold of dimension $2m_{\Fa}+k$ with a locally free action of $\T/\T_\Fa$.
Moreover, the levels of $\mu$ are compact connected submanifolds of $M$.
\end{numlist}
\end{prop}
\begin{numlproof}

\item Let $\T_\Fa$ be the intersection of the isotropy subgroups of points
in $M_\Fa$. Then the connected component of the identity in $\T_\Fa$ is an
$(\ell-k)$-dimensional subtorus of $\T$, and $M_\Fa$ is a connected
component of its fixed point set. Since $\T_\Fa$ acts on $M$ effectively
by hamiltonian isometries, $M_\Fa$ is a compact totally geodesic K\"ahler
submanifold of $M$, of dimension at most $2m-2(\ell-k)$. By definition,
$M_\Fa$ carries an effective hamiltonian action of $\T/\T_\Fa$ (which is
connected, hence a $k$-torus), so it has dimension at least $2k$. The
momentum map is essentially $\mu$, viewed as a map from $M_\Fa$ to the
affine span of $\Fa$, so the action is rigid.

\item By Lemma~\ref{tricky}, the critical values of $\mu$, regarded in the
above way, are precisely the boundary points of $\Fa$, and $\mu$ is regular as
a map from $M_{\Fa}^0$ to $\Fa^0$. The gradient flow of $\mu$ commutes with
$\T$ and hence provides an equivariant trivialization of $M_{\Fa}^0$. Thus
$M_\Fa^0$ is diffeomorphic to $\Fa^0\times P_\Fa$ and the action of
$\T/\T_\Fa$ is given by an effective locally free action on $P_\Fa$, with
trivial action on $\Fa^0$. The levels of $\mu$ are smooth since any point in
the image of $\mu$ is in some open face; they are connected by~\cite{atiyah}.
\end{numlproof}

This shows that `wild' $S^1$ actions on $\C P^2$ (as a symplectic manifold) of
Example~\ref{ex:cp2} cannot be rigid with respect to any compatible K\"ahler
metric. One can easily check that the `tame' actions are rigid with respect to
the Fubini--Study metric.

\subsection{The symplectic isotropy representations}\label{s:isorep}

We now wish to obtain precise information about the symplectic isotropy
representations of the torus action.  If $\mu(x)$ belongs to an open
$k$-dimensional face $\Fa^0$, then the Lie algebra ${\mathfrak t}_x$ of the
isotropy group $\T_x\geq \T_\Fa$ of $x$ is the vector subspace of elements of
$\mathfrak t$, annihilated by the elements of the vector subspace of
$\mathfrak t^*$ parallel to $\Fa$: indeed this is clearly the image of
$d\mu_x$, and $\mathfrak t_x$ is the kernel of the transpose of $d\mu_x$.

Since the orbit $\T\cdot x$ is $k$-dimensional, the symplectic isotropy
representation $V_x=T_x(\T\cdot x)^0/T_x(\T\cdot x)$ of $\T_x$ (and its Lie
algebra ${\mathfrak t}_x$) has dimension $m-k$. Hence it is an orthogonal
direct sum of $m-k$ complex $1$-dimensional representations with (not
necessarily distinct) characters $\T_x\to S^1$. Differentiating this action
gives the weights $\alpha_1,\ldots \alpha_{m-k}$ of the action of ${\mathfrak
t}_x$, which are integral elements of ${\mathfrak t}_x^*$.

Since the ${\mathfrak t}_x$ action is effective, the weights $\alpha_1,\ldots
\alpha_{m-k}$ span ${\mathfrak t}_x^*$, and we order them so that
$\alpha_1,\ldots \alpha_{\ell-k}$ form a basis for ${\mathfrak t}_x^*$.

\begin{lemma}\label{weights}
Suppose $\mu(x)$ belongs to an open $k$-dimensional face $\Fa^0$ of $\Delta$
and let $V_x$ be the symplectic isotropy representation of $\T_x$ at $x\in M$.
\begin{numlist}
\item The induced ${\mathfrak t}_x$ action has exactly $\ell-k$ distinct
nonzero weights.
\item $\T_x$ is connected.
\end{numlist}
\end{lemma}
\begin{numlproof}
\item We choose a projection $\chi\colon{\mathfrak t}\to {\mathfrak t}_x$ and
introduce a symplectic slice as in Lemma~\ref{sslice}. Thus there is a
$\T$-equivariant symplectomorphism from a neighbourhood $U$ of the zero
section $0_N$ in the normal bundle $N\to\T\cdot x$ to a neighbourhood of
$\T\cdot x$ in $(M,\omega)$, where the normal bundle $N=\T
\times_{\T_x}({\mathfrak t}_x^0\oplus V_x) \to \T\cdot x$ is realised as a
symplectic quotient of $T^*\T\times V_x$ by the diagonal action of $\T_x$.
The symplectomorphism identifies $0_N$ with $\T\cdot x$ and its differential
along the zero section is essentially the identity map. Let us denote the
pullback of $(g,J,\omega)$ by $(g_0,J_0,\omega_0)$.  We then have that $g_0$
agrees with $g_x$ at $x$.

We now bring in the rigidity condition that the induced metric on $\T$ depends
only on $\mu$. This implies that for any vector fields $X_\xi,X_\eta$
($\xi,\eta\in\mathfrak t$) induced by the action of $\T_x$ on
$(U,g_0,J_0,\omega_0)$, $g_0(X_\xi,X_\eta)$, as a function on $U$, depends
only on the momentum map $\mu_0$ of $N$, $\mu_0([\alpha,v])=
\alpha+\mu_V(v)\circ\chi$ with $\mu_V=\frac12\sum_{i=1}^{m-k} |z_i|^2
\alpha_i$, where $z_1,\ldots z_{m-k}$ are the standard complex coordinates on
the weight spaces in $V_x$.  It follows from~\cite{schwarz} that (being smooth
on $U$) $g_0(X_\xi,X_\eta)$ is a smooth function of $\mu_0$. In particular,
for $\alpha=0 \in {\mathfrak t}_x^0$, $g_0(X_\xi,X_\eta)$ is a smooth function
of $\mu_V$. Thus, on $V_x\cap U$, $d(g_0(X_\xi,X_\eta))$ is a pointwise linear
combination of the components of
\begin{equation}\label{dmu0}
d\mu_V = \frac12\sum_{i=1}^{m-k} (z_i d\overline z_i+\overline z_i dz_i)
\alpha_i.
\end{equation}
In other words (since it vanishes at the origin of $V_x$) it equals $\langle
d\mu_V, B(\xi,\eta)\rangle$ for a smooth bilinear form $B\colon V_x\cap U\to
S^2{\mathfrak t}_x^*\otimes{\mathfrak t}_x$.  Now since $X_\xi$ and $X_\eta$
vanish at the origin, $g_0(X_\xi,X_\eta)$ differs from
$g_x(X_\xi,X_\eta)=\sum_{i=1}^{m-k} \alpha_i(\xi)\alpha_i(\eta) |z_i|^2$ by a
smooth function vanishing to second order at the origin, so its exterior
derivative on $V_x\cap U$ is, to first order, equal to
\begin{equation}\label{dg0}
d(g_x(X_\xi,X_\eta))= \sum_{i=1}^{m-k}
(z_i d\overline z_i+\overline z_i dz_i)\alpha_i(\xi)\alpha_i(\eta).
\end{equation}
If we differentiate $d(g_0(X_\xi,X_\eta))=\langle d\mu_V, B(\xi,\eta)\rangle$
with respect to $\overline z_i$, using~\eqref{dmu0} and~\eqref{dg0}, and
evaluate at the origin of $V_x$, the error terms and derivative of $B$ go
away. Equating coefficients of $dz_1,\ldots dz_{m-k}$ therefore gives
\begin{equation*}
2\alpha_i(\xi)\alpha_i(\eta) = \alpha_i(B_0(\xi,\eta))
\end{equation*}
for all $i$, i.e., $B_0^*\alpha_i = 2\alpha_i\otimes \alpha_i$. (We remark
that this generalizes the conditions~\eqref{eq:strongtoricboundary} in the
toric case.)  Now $\alpha_1,\ldots \alpha_{\ell-k}$ is a basis for ${\mathfrak
t}_x^*$, so we may write $\alpha_{\ell-k+1},\ldots \alpha_{m-k}$ as
$\alpha_i=\sum_{j=1}^{\smash{\ell-k}} \lambda_{ij} \alpha_j$.  We then deduce from
$B_0^*\alpha_i = 2\alpha_i\otimes \alpha_i$ that
\begin{equation*}
\lambda_{ij} \lambda_{ik} = \delta_{jk} \lambda_{ij}.
\end{equation*}
Thus for each $i$, $\lambda_{ij}$ is nonzero for at most one $j$, and then
equal to one, i.e., for any $i=\ell-k+1,\ldots m-k$, the weight $\alpha_i$ is
either zero, or it is one of $\alpha_1,\ldots \alpha_{\ell-k}$.

\item We prove that all isotropy groups of the $\T$-action are connected.
Since the gradient flow of $\mu$ commutes with $\T$, it suffices to prove this
near a fixed point $y$ of the $\T$-action, where the symplectic slice gives a
$\T$-equivariant symplectomorphism with a neighbourhood of the origin in a
symplectic vector space $V_y$. Now since the $\T$-action on $V_y$ is
effective, with $\ell$-distinct nonzero weights, these form a basis for the
dual lattice. This ensures the isotropy groups of points in $V_y$ are
connected.
\end{numlproof}

Part (i) of Lemma~\ref{weights} is the key to the theory of rigid hamiltonian
torus actions. In particular it allows us to refine Proposition~\ref{strata}.

\begin{prop} \label{delzant}
Suppose $(M,g,J,\omega)$ is a compact connected K\"ahler manifold with a
rigid hamiltonian $\ell$-torus action, as in
Proposition~\textup{\ref{strata}}.
\begin{numlist}
\item If $\Fa^0$ is an open $k$-dimensional face, then the isotropy group of
all points in $M_\Fa^0$ is an $(\ell-k)$-torus $\T_\Fa$, and the isotropy
representations are all equivalent, with the distinct nonzero weights in
${\mathfrak t}_\Fa^*$ forming a basis for the lattice dual to the lattice of
circle subgroups of ${\mathfrak t}_\Fa$.
\item The image $\Delta$ of $\mu$ is a Delzant polytope.
\item $P_\Fa$ is a principal $k$-torus bundle \textup(under
$\T/\T_\Fa$\textup) over a compact manifold $S_\Fa$ of dimension $2m_\Fa$,
with a family of K\"ahler structures parameterized by $\Fa^0$.
\end{numlist}
\end{prop}
\begin{numlproof}
\item This is immediate from Lemma~\ref{weights}: the distinct nonzero weights
form a basis for ${\mathfrak t}_\Fa$, the Lie algebra of the (connected)
isotropy group of any point in $M_\Fa^0$.

\item Applying this to a fixed point, observe that the directions of the
distinct nonzero weights are the edges meeting the corresponding vertex of
$\Delta$. There are $\ell$ of these and the dual basis gives a basis for the
lattice of circle subgroups of $\T$ consisting of normals to the faces
meeting the vertex.

\item By Proposition~\ref{strata}, $P_\Fa$ has a locally free action of
$\T/\T_{\Fa}$, and by Lemma~\ref{weights}, the isotropy groups are
connected, so the action is free. Hence $P_\Fa$ is a principal $\T/\T_{\Fa}$
bundle over a compact manifold $S_\Fa$.  Choosing a point $v$ in $\Fa^0$
identifies $S_\Fa$ with the K\"ahler quotient of $M_\Fa$ at momentum level
$v$.
\end{numlproof}

\subsection{The complexified torus action}\label{s:complextorus}

We now turn to the structure of the orbits of the complexified torus
action. If the $\T$ action is generated by vector fields $K_1,\ldots K_\ell$,
then the complexified action of $\T^c$ is generated by the (real) holomorphic
vector fields $K_1,\ldots K_\ell, JK_1,\ldots JK_\ell$.  These are linearly
independent on a dense open set (since the $\T$ action is hamiltonian) and
generate a foliation of $M$ by complex orbits, whose generic leaf is
$2\ell$-dimensional. As we have already remarked in \S \ref{s:stratification},
$JK_1,\ldots JK_\ell$ generate the gradient flow of $\mu$, and therefore the
momentum image of a $2k$-dimensional leaf is a $k$-dimensional open face
$\Fa^0$ of $\Delta$; the isotropy group of any point in this leaf is the
complexification $\T^c_\Fa$ of $\T_\Fa$ and the closure (in $M$) of the leaf
maps onto the closed face $\Fa$.

To understand the complex orbits further, we reinterpret $V_x$ as the fibre
of the normal bundle to $\T^c\cdot x$ at $x$, carrying the complex isotropy
representation, and we linearize the $\T^c$ action using a holomorphic slice
rather than a symplectic one.

In general, let $G$ be a compact Lie group of hamiltonian isometries of a
K\"ahler manifold $M$, and let $G^c$ be the complexification, which acts
holomorphically on $M$.  Then the {\it holomorphic slice
theorem}~\cite{heinloos,sjamaar} states that if $G^c\cdot x$ is the orbit
through $x\in M$ with isotropy representation $(G^c_x,V_x)$, then there is a
$G^c$-equivariant biholomorphism from a neighbourhood of $G^c\cdot x$ in $M$
to a neighbourhood of the zero section in $G^c\times_{G^c_x} V_x\to G^c\cdot
x$.

\begin{rem}
For many purposes, it suffices to know that a neighbourhood of $x$ is {\it
locally} $G^c$-equivariantly biholomorphic to a neighbourhood of the zero
section in $G^c\times_{G^c_x} V_x$. This is quite easy to establish.
Indeed, let $\psi\colon U\to M$ be a holomorphic chart with $\psi(0)=x$ and
$d\psi_0=\Id$, where $U$ is an open neighbourhood of the origin in $T_x M$.
We can assume $U$ and $\psi$ are $G_x$-equivariant by averaging, since $G_x$
is compact. Now by acting with $G$, we obtain a $G$-equivariant
biholomorphism $\tilde\psi$ from a neighbourhood $\tilde U$ of $G\cdot x$ in
$M$ to a neighbourhood of the zero section in $G\times_{G_x}\tilde V_x\to
G\cdot x$.  Here $\tilde V_x$ is the orthogonal complement of $T_x(G\cdot
x)$: note $\tilde V_x = V_x\oplus W_x$ where $V_x$ is the orthogonal
complement of $T_x(G^c\cdot x)$, and $W_x = JT_x(G\cdot x)$.

Now since $\tilde\psi$ is holomorphic and $G$-equivariant, it is (locally)
$G^c$-equivariant. This is only a local result, because the domain $\tilde U$
is \textit{a priori} only $G$-invariant, not $G^c$-invariant. The hard part of
the holomorphic slice theorem is to show such a `local' slice can be
analytically continued to a $G^c$-invariant neighbourhood of $G^c\cdot x$.
\end{rem}

\begin{lemma} \label{slice}
Suppose $\mu(x)$ belongs to an open $k$-dimensional face $\Fa^0$ of $\Delta$
and let $\T^c_1,\T^c_2,\ldots \T^c_{\ell-k}$ be the complexifications of the
circle subgroups of the isotropy subgroup $\T_\Fa$ dual to the basis of
distinct nonzero weights in the symplectic isotropy representation of
$\T_\Fa$.

Then $\T^c_\Fa=\T^c_1\times\cdots\times \T^c_{\ell-k}$ and there is a
$\T^c$-equivariant biholomorphism from a neighbourhood $U$ of $\T^c\cdot x$
in $M$ to a neighbourhood $W$ of the zero section in
\begin{equation*}
\T^c\times_{\T^c_{\Fa}}\bigl(
V_0\oplus V_1\oplus\cdots\oplus V_{\ell-k}\bigr)\to \T^c/\T^c_{\Fa}
\end{equation*}
where $V_0$ is the trivial representation \textup(possibly zero\textup),
while for $i=1,\ldots\ell-k$, $V_i$ is a nonzero vector space carrying the
standard action of $\T^c_i\cong\C^\times$ by scalar multiplication, with
$\T^c_j$ acting trivially for $j\neq i$.  Under this biholomorphism\textup:
\begin{numlist}
\item the $p$-dimensional faces $\tFa$ meeting $\Fa$ correspond
bijectively to $(p-k)$-element subsets $\J_{\tFa}\subseteq\{1,\ldots
\ell-k\}$ in such a way that
\begin{bulletlist} \item
$M_{\tFa}\cap U$ is the intersection of $W$ with those elements whose
$V_j$ component vanishes for $j\in \{1,\ldots \ell-k\} \smallsetminus
\J_{\tFa}$\textup;
\end{bulletlist}

\item if $Y$ is a $p$-dimensional complex orbit with $x\in \overline
Y\subseteq M_{\tFa}$, $\dim \tFa=p$ then there are one dimensional subspaces
of $V_j$ for $j\in \J_{\tFa}$ such that
\begin{equation}
\overline Y\cap U\cong\T^c\times_{\T^c_{\Fa}}
\textstyle\bigoplus_{j\in \J_{\tFa}} L_j
\end{equation}
under the obvious inclusion into $\T^c\times_{\T^c_{\Fa}}\bigl( V_0\oplus
\cdots\oplus V_{\ell-k}\bigr)$.
\end{numlist}
\end{lemma}
\begin{proof}
By the holomorphic slice theorem there is a $\T^c$-equivariant
biholomorphism from a neighbourhood of $\T^c\cdot x$ to neighbourhood of the
zero section in $\T^c\times_{\T^c_{\Fa}} V_x$ where $V_x$ is normal to
$\T^c\cdot x$ at $x$. Equivalently, $V_x$ is the symplectic isotropy
representation of $\T_\Fa$, now equipped with the natural complexified
action of $\T^c_\Fa$.  By Lemma~\ref{weights}, the distinct nonzero weights
of the ${\mathfrak t}_\Fa$ action on $V_x$ are dual to a basis for the
lattice of circle subgroups of $\T_\Fa$, and we take the $V_i$'s to be the
weight spaces (with $V_0$ the zero weight space). This gives what we want.

\begin{numlist}
\item It is clear that the faces $\tFa$ containing $\Fa$ correspond to
subsets $\J_{\tFa}$ of $\{1,\ldots \ell-k\}$ with $\T^c_j$ acting
nontrivially on $M_{\tFa}$ for $j\in \J_{\tFa}$.  The biholomorphism
identifies $M_{\tFa}\cap U$ with those elements of $W$ whose isotropy group
is contained in $\T^c_{\tFa}$. Since the latter is the product of the
$\T^c_j$ for $j\in \{1,\ldots \ell-k\}\smallsetminus \J_{\tFa}$, the
result follows.

\item Under the biholomorpism, the complex orbits $Y$ near $\T^c\cdot x$ are
all of the form $\T^c\times_{\T^c_{\Fa}} (v_0+U_1\times\cdots\times
U_{\ell-k})$, where $v_0\in V_0$ and either
$U_j=L_j^\times:=L_j\smallsetminus\{0\}$, where $L_j$ is a one-dimensional
subspace of $V_j$, or $U_j=\{0\}\subset V_j$.

If $Y$ is a $p$-dimensional orbit in $M_{\smash{\tFa}}^{0}$, then these two
cases occur accordingly as $j\in \J_{\tFa}$ or not.  Clearly $x\in
\overline Y$ if and only if $v_0=0$, and then the biholomorpism
identifies $\overline Y\cap U$ with $\bigoplus_{j\in \J_{\tFa}} L_j$ as
stated.  \endnumlproof

Lemma~\ref{slice} gives a lot of information about the equivariant holomorphic
geometry of $M$. For instance, applying it at a fixed point gives a
$\T^c$-equivariant chart from a neighbourhood of the fixed point to
$U_0+V_1\oplus\cdots \oplus V_\ell$, where $V_1,\ldots V_\ell$ are the
nontrivial weight spaces associated to the corresponding vertex $v$ of
$\Delta$, and $U_0$ is a neighbourhood of the origin in the trivial weight
space $V_0$. In the toric case, $V_0=0$ and $\dim V_j=1$ for all $j$, and we
obtain the linear charts underlying the toric complex manifold. In the general
case, such charts provide a finite atlas, since there are finitely many
vertices $v$ and they have compact preimages $S_v=\mu^{-1}(v)$.

\begin{prop} \label{complex}
Suppose $(M,g,J,\omega)$ is a compact connected K\"ahler manifold with a rigid
hamiltonian $\ell$-torus action, as in Proposition~\textup{\ref{strata}}.
\begin{numlist}
\item The closure of a $2k$-dimensional complex orbit in $M$ is a toric
K\"ahler submanifold of $M$ whose Delzant polytope is a $k$-dimensional face
$\Fa$ of $\Delta$.

\item For any $k$-dimensional face $\Fa$ of $\Delta$, $M_\Fa^0=\Fa^0\times
P_\Fa$ is a holomorphic principal $\T^c/\T^c_\Fa$-bundle over a complex
manifold $S_\Fa$.

\item The blow-up of $M_\Fa$ along the inverse images of the codimension one
faces of $\Fa$ is equivariantly biholomorphic to the total space of
$M_\Fa^0\times_{\T^c/\T^c_\Fa} \cV_\Fa\to S_\Fa$ for some smooth toric
complex manifold $\cV_{\Fa}$.

\item If $\Fa$ is a $k$-dimensional face, with the $(k-1)$-dimensional face
$\tFa$ in its boundary, then $S_{\Fa}$ is a holomorphic $\C P^d$-bundle over
$S_{\tFa}$ with $d=m_{\Fa}-m_{\tFa}\geq 0$.

Furthermore if $Q_\Fa$ denotes the fibrewise Hopf fibration over the $\C
P^d$-bundle $S_{\Fa}\to S_{\tFa}$, then $P_{\Fa}\to S_{\Fa}$ is the pullback
of $P_{\tFa}\to S_{\tFa}$ along the $S^{2d+1}$-bundle map $Q_\Fa\to
S_{\tFa}$ composed with the $S^1$-bundle map $Q_\Fa\to S_\Fa$.
\end{numlist}
\end{prop}
\begin{numlproof}
\item For all $x\in M$, any complex orbit has a smooth closure along
$\T^c\cdot x$ by Lemma~\ref{slice}. Hence the closures of the complex orbits
are smoothly embedded, and become toric K\"ahler manifolds under the induced
metric.  We have already remarked that $\mu$ maps any such orbit closure to
a face $\Fa$ of $\Delta$, and clearly $\mu$, viewed as a map to the affine
span of $\Fa$ (with a choice of origin), is a momentum map for the induced
toric action.

\item For convenience, we prove this result for $\Fa=\Delta$: the general
result follows by replacing $M$ with $M_\Fa$ and $\T^c$ by $\T^c/\T^c_\Fa$.

Since $\T^c$ acts freely on $M^0$ it defines a holomorphic fibration over
$S_\Delta$. To verify that the fibration is locally trivial, observe that a
neighbourhood of a $\T^c$ orbit in $M^0$ is equivariantly biholomorphic to a
neighbourhood of the zero section in $\T^c\times V_0\to \T^c$.  The latter,
being $\T^c$-invariant, is of the form $\T^c\times U_0$, and the projection
to $U_0$ gives the required local trivialization. Since $\T^c$ acts simply
transitively on the fibers, $M^0$ is a principal $\T^c$-bundle over
$S_{\Delta}$.

\item We again prove the result when the face is the whole polytope
$\Delta$.

We first consider the blow-up $\hat M$ of $M$ along all $M_\Fa$ with $\Fa$
codimension one in $\Delta$. (Of course the blow-up is
trivial if $M_\Fa$ already has complex codimension one in $M$). Thus, 
$\hat M$ is the complex manifold  obtained from $M$
by replacing each $M_\Fa$ by its projectivized normal bundle $\hat 
M_\Fa$; these become divisors (i.e.
 of complex codimension one) in 
$\hat M$, and the $\T^c$ action
lifts naturally to $\hat M$.  Lemma~\ref{slice} shows that the generic
$\T^c$ orbits for the lifted action have disjoint smooth closures in 
$\hat M$, and this gives a
holomorphic fibration of $\hat M$ whose fibres are all toric K\"ahler
manifolds with Delzant polytope $\Delta$. In particular (forgetting the
symplectic structure) they are all isomorphic toric
complex manifolds~\cite{guillemin0,LT}.

Let $\cV_\Delta$ be a toric complex manifold in this isomorphism class, and
choose a basepoint on the generic orbit $\cV_\Delta^0$ to identify it with
$\T^c$.  Then there is an equivariant biholomorphism $M^0\times_{\T^c}
\cV_\Delta^0 \to M^0= \smash{\hat M}^0$ (here ${\hat M}^0$ stands for the
subset of points of $\hat M$ with generic $\T^c$ orbits; it is the same as
$M^0$ because the blow-up is the identity on the complement of the exceptional
divisor). Since $\cV_\Delta$ has the same isotropy representations as the
fibres of $\hat M$, this extends to an equivariant biholomorphism
$M^0\times_{\T^c} \cV_\Delta \to \hat M$ (indeed the holomorphic slices of
Lemma~\ref{slice} provide the extension).

\item Consider, as in (iii), the blow-up $\hat M$ of $M$ along its
codimension one faces.  This is equivariantly biholomorphic to
$M^0\times_{\T^c} \cV_\Delta$ and for \emph{any} face $\Fa$, the inverse
image $\hat M_\Fa$ of $M_\Fa$ in $\hat M$ is $M^0\times_{\T^c} \cV_\Fa$
(where only $\T^c/\T^c_\Fa$ acts effectively on $\cV_{\Fa}$, which is the
inverse image of $\Fa$ in $\cV_\Delta$).

Now $\cV_{\Fa}^0$ is equivariantly biholomorphic to a $\T^c_{\tFa}/\T^c_{\Fa}$
bundle over $\cV_{\tFa}^0$, namely the punctured normal bundle of
$\cV_{\tFa}^0$ in $\cV_{\Fa}$, so it follows that the same is true for
$\smash{\hat M}_{\Fa}^0$: it is equivariantly biholomorphic to the punctured
normal bundle of $\smash{\hat M}_{\tFa}^0$ in $\smash{\hat M}_{\Fa}$.  Passing
to the blow-down, we deduce that $M_{\Fa}^0$ is equivariantly biholomorphic to
the punctured normal bundle of $M_{\tFa}^0$ in $M_{\Fa}$, which is a
$\T^c/\T^c_{\Fa}$-equivariant bundle with $\T^c_{\tFa}/\T^c_{\Fa}$ acting by
scalar multiplication on the fibres.

The quotient by $\T^c/\T^c_{\Fa}$ identifies $S_{\Fa}$ biholomorphically with
a bundle over $S_\Fa$. To describe this bundle, we first divide the punctured
normal bundle of $M_{\tFa}^0$ by $\T^c_{\tFa}/\T^c_{\Fa}$ to obtain the
projectivized normal bundle as a $\T^c/\T^c_{\tFa}$-equivariant $\C P^d$
bundle over $M_{\tFa}^0$ with trivial action on the fibres. Now the quotient
by $\T^c/\T^c_{\tFa}$ shows that $S_{\Fa}\to S_{\tFa}$ is a holomorphic $\C
P^d$-bundle.

The unit normal bundle of $M_{\tFa}^0$ is the sphere bundle induced by the
Hopf fibtration over the projectivized normal bundle and the result follows.
\end{numlproof}

This shows that `wild' $S^1$ actions on $\C P^2$ (as a complex manifold)
discussed in Example~\ref{ex:cp2} cannot be rigid with respect to any
compatible K\"ahler metric.  On the other hand, we noted there that the
complex orbits of `tame' $S^1$ actions do indeed have smooth closures.

\subsection{K\"ahler geometry of rigid hamiltonian torus actions}
\label{s:kahlerrigid}

Given a K\"ahler $2m$-manifold $M$ with a rigid hamiltonian action of an
$\ell$-torus $\T$, we have obtained a description of the equivariant
biholomorphism type of $M$, stratified by the inverse images of the faces of
the momentum polytope $\Delta$: $M^0$ is a principal $\T^c$-bundle over a
complex manifold $S_\Delta$ of dimension $2m_\Delta$, with $m_\Delta=m-\ell$,
and there is a toric complex manifold $\cV_\Delta$ such that the blow up of
$M$ along the codimension one faces of $\Delta$ is biholomorphic to
$M^0\times_{\T^c} \cV_\Delta\to S_\Delta$; {\it mutatis mutandis}, the inverse
image $M_\Fa=\mu^{-1}(\Fa)$ of a face of $\Delta$ has the same structure;
further if $\Fa_1,\ldots \Fa_n$ denote the codimension one faces of $\Delta$,
then $S_{\Fa_j}$ has dimension $2m_{\Fa_j} \leq 2m_\Delta$ and $S_\Delta$ is a
$\C P^{\smash{d_j}}$-bundle over $S_{\Fa_j}$ with $d_j=m_\Delta-m_{\Fa_j}$,
and we say a \emph{blow-down occurs} over $\Fa_j$ if $d_j>0$. (We remark that
if $\tFa$ is a codimension one face of $\Fa$, it must be $\Fa\cap\Fa_j$ for
some codimension one face $\Fa_j$ of $\Delta$. We then have
$m_{\Fa}-m_{\tFa}=d_j=m_\Delta-m_{\Fa_j}$.)

It remains to descibe the K\"ahler structure of $M$ in terms of this
equivariant biholomorphism type. To do that we first recall some equivalent
formulations of the rigidity condition established (locally) in \cite{ACG1}.

Suppose, generally, that $M$ is a K\"ahler manifold endowed with an isometric
hamiltonian action of an $\ell$-torus $\T$ with momentum map $\mu$.  For a
contractible open subset $U$ of the regular values of $\mu$, the gradient flow
of $\mu$ identifies $\mu^{-1}(U)$ with $\mu^{-1}(v)\times U$ for any $v$ in
$U$, and hence $\mu^{-1}(U)/\T \cong S\times U$ for a complex manifold $S$,
with a family $\omega_h$ of compatible symplectic forms on the fibres of
$S\times U\to U$.  We can therefore define the derivative $d_\mu\omega_h$ with
respect to $\mu$, and this will be a $2$-form on $S$ with values in $\mathfrak
t$. Now $\mu^{-1}(U)$ is a principal $\T$-bundle with connection over $S\times
U$, so it has a curvature form $\boldsymbol\Omega$, which is also a closed
$2$-form with values in $\mathfrak t$. If $d_\mu\omega_h=\boldsymbol\Omega$ on
$S\times U$ we say that the {\it rigid Duistermaat--Heckman property} holds
(so-called because it holds in cohomology by work of Duistermaat and
Heckman). We then have the following global version of \cite[Proposition
8]{ACG1}.

\begin{lemma}\label{l:rigid} For an isometric hamiltonian $\T$-action the
following are equivalent.
\begin{numlist}
\item The action is rigid.
\item The $\T^c$-orbits are totally geodesic.
\item The orthogonal distribution to the $\T^c$-orbits is $\T^c$-invariant.
\item The rigid Duistermaat--Heckman property holds.
\end{numlist}
\end{lemma}
\begin{proof} This is essentially the same as \cite[Proposition 8]{ACG1}.
Let $X$ denote a vector field which is orthogonal to a $\T^c$-orbit. The
rigidity condition is equivalent to the statement that $\partial_X
(g(K_r,K_s))=-2g(\nabla_{K_r} K_s,X)$ vanishes along the given orbit for all
such vector fields $X$. Since $J$ is parallel and $K_s$ is holomorphic this is
equivalent to the fact that the $\T^c$ orbit is parallel. It is easy to
compute that this condition is equivalent to the fact that ${\cL}_{K_r} X$ and
${\cL}_{JK_r} X$ are orthogonal to the given $\T^c$ orbit for all $X$,
$r=1,\ldots\ell$, i.e., the orthogonal distribution is $\T^c$-invariant.

(iv) is equivalent to the local rigidity of the action on $M^0$ by the
Pedersen--Poon construction (see~\cite{ped-poon,ACG1}); this implies rigidity
on $M$ by continuity.
\end{proof}

We next show that a compact K\"ahler manifold with a rigid hamiltonian action
of a torus gives rise in a natural way to the following data.

\begin{defn} Let $\cV$ be a compact toric K\"ahler manifold under an
$\ell$-torus $\T$ with Delzant polytope $\Delta$.  Then {\it rigid
hamiltonian data} for $\cV$ consists of a quadruple
$(\cV_\Fa,S_\Fa,P_\Fa,\omega_\Fa)$ for each face $\Fa$ of $\Delta$, where:
\begin{numlist}
\item $\cV_\Fa$ is the inverse image of $\Fa$ in $\cV$, which is a
compact toric K\"ahler manifold under $\T/\T_\Fa$, where $\T_\Fa$ is the
isotropy subgroup of $\T$ associated to $\Fa$;

\item $S_\Fa$ is a compact complex manifold which is a holomorphic
projective space bundle over $S_{\tFa}$ for any codimension one face
${\tFa}$ of $\Fa$;

\item $\pi\colon P_\Fa\to S_\Fa$ is a principal $\T/\T_\Fa$-bundle
with connection $\boldsymbol \theta_\Fa\colon TP_\Fa\to\mathfrak t/\mathfrak
t_\Fa$, whose curvature $\boldsymbol\Omega_\Fa\in
C^\infty(S_\Fa,\Lambda^{1,1}S_\Fa\otimes\mathfrak t/\mathfrak t_\Fa)$ pulls
back to the fibres of $S_\Fa\to S_{\tFa}$ to give the Fubini--Study metric
in $2\pi c_1(\cO(1))$ tensored with the \textup(primitive inward\textup)
normal to the codimension one face $\tFa$;

\item $\omega_\Fa$ is a section of (the pullback of) $\Lambda^{1,1}S_\Fa$
over $S_\Fa\times \Fa$, which
\begin{itemize}
\item is positive on $S_\Fa\times \Fa^0$,
\item satisfies $d_\mu\omega_\Fa = \boldsymbol\Omega_\Fa$ on $S_\Fa\times\{v\}$
for all $v\in\Fa^0$,
\item and whose restriction to $S_\Fa\times \tFa$, for any codimension one
face $\tFa$ of $\Fa$, is the pullback of $\omega_{\tFa}$ along the map
$S_\Fa\times\tFa\to S_\tFa\times \tFa$.
\end{itemize}
\end{numlist}
\end{defn}

\begin{prop}\label{data}
Let $M$ be a compact connected K\"ahler $2m$-manifold with a rigid
hamiltonian action of an $\ell$-torus $\T$ and momentum map $\mu\colon M\to
\Delta$. Then there are rigid hamiltonian data
$(\cV_\Fa,S_\Fa,P_\Fa,\omega_\Fa)$ \textup(for the faces $\Fa$ of
$\Delta$\textup) associated to a toric K\"ahler manifold $\cV$ with Delzant
polytope $\Delta$ such that\textup:
\begin{bulletlist}
\item the pullback of the K\"ahler metric on $M_\Fa=\mu^{-1}(\Fa)$ to the
fibres of the blow-up $\smash{\hat M_\Fa}\cong P_\Fa\times_\T \cV_\Fa$
\textup(see Proposition~\textup{\ref{complex}}\textup) is
induced by the K\"ahler metric on $\cV_\Fa$\textup;
\item $S_\Fa$ is the K\"ahler quotient of $M_\Fa$ by $\T/\T_\Fa$
and the K\"ahler quotient metric at momentum level $v\in\Fa^0$ is induced by
$\omega_\Fa$ on $S_\Fa\times\{v\}$\textup;
\item the orthogonal distribution to the generic $\T^c/\T^c_\Fa$ orbits
in $M_\Fa$ is the joint kernel of $\boldsymbol\theta_\Fa$ and $d\mu$.
\end{bulletlist}
In particular, on $M^0\cong P_\Delta\times_\T \cV_\Delta^0$, the K\"ahler
structure is given by
\begin{equation}\label{hkf}\begin{split}
g &= h_0 + \langle \mu, \boldsymbol h\rangle+ \langle d\mu, {\bf G} , d\mu
\rangle+ \langle\boldsymbol \theta,{\bf G}^{-1},
\boldsymbol\theta\rangle,\\
\omega &= \Omega_0 + d \langle \mu, \boldsymbol\theta \rangle = \Omega_0 +
\langle \mu,\boldsymbol \Omega\rangle +\langle d\mu\wedge \boldsymbol\theta
\rangle,
\end{split}\end{equation}
where $\omega_\Delta= \Omega_0 + \langle \mu,\boldsymbol \Omega\rangle$, $h_0
+ \langle \mu, \boldsymbol h\rangle$ is the corresponding family of hermitian
metrics, $\boldsymbol\theta=\boldsymbol\theta_\Delta$,
$\boldsymbol\Omega=\boldsymbol\Omega_\Delta$, the toric K\"ahler metric on
$\cV_\Delta^0$ is given by~\eqref{toricmetric}, for ${\bf G} \colon
\Delta^0\to S^2\mathfrak t$, and \textup(as before\textup) angled brackets
denote pointwise contractions.
\end{prop}
\begin{proof} It suffices to prove the result for the whole polytope
$\Delta$. We know by Propositions~\ref{delzant} and \ref{complex} that the
blow-up $\hat M$ is equivariantly biholomorphic to $M^0\times_{\T^c}
\cV_\Delta\cong P_\Delta\times_{\T}\cV_\Delta$ for a toric complex manifold
$\cV_\Delta$, and the fibres of $P_\Delta\times_{\T}\cV_\Delta\to S_\Delta$
map biholomorphically onto the complex orbit closures in $M$.  The K\"ahler
metric of $M$ induces a K\"ahler structure on each complex orbit closure,
which depends only on the momentum map $\mu$. Since $\mu$ is $\T$-invariant,
there is a toric K\"ahler structure on $\cV_\Delta$, with Delzant polytope
$\Delta$, such that the fibres of $P_\Delta\times_{\T}\cV_\Delta$, with the
metric induced from $\cV_\Delta$, map \emph{isometrically} onto the complex
orbit closures in $M$.

The K\"ahler metric on $M^0$ induces a principal $\T$-connection on $M^0\to
B_\Delta=S_\Delta\times\Delta^0$ (the orthogonal distribution to the fibres),
and by Lemma~\ref{l:rigid}, this is the pullback of a principal
$\T$-connection $\boldsymbol\theta$ on $\pi\colon P_\Delta\to S_\Delta$.  The
lemma also shows that the family $\omega_\Delta$ of K\"ahler forms induced on
$S_\Delta$ depends affinely on $\mu\in\Delta^0$ and $\pi^*
d_\mu\omega_\Delta=d\boldsymbol\theta$, so the linear part is the curvature
$\boldsymbol\Omega$ of the connection $\boldsymbol\theta$ for all
$\mu\in\Delta^0$; $\omega_\Delta$ is therefore smoothly defined for all $\mu$.

The K\"ahler form on $M$ pulls back to the blow-up $\hat M$ to give a $2$-form
which degenerates on the exceptional divisor. Using the description of this
divisor given in Proposition~\ref{complex} and the smooth dependence of the
K\"ahler form on $\mu$, it follows that the K\"ahler form $\omega_\Delta$
approaches to the pullback of $\omega_\Fa$ along $S_\Delta\to S_\Fa$ as
$\mu\rightsquigarrow\Fa^0\subset\Fa$, for a codimension one face $\Fa$ of
$\Delta$. We then deduce that the pullback of $d_\mu\omega_\Delta$ to a fibre
of $S_\Delta\to S_{\Fa}$ takes values, for $\mu\in\Fa^0$, in the annihilator
$\mathfrak t_\Fa$ of $T_\mu \Fa$, i.e., is of the form $\Omega\otimes u_\Fa$,
where $u_\Fa$ is the primitive inward normal to $\Fa$, and $\Omega$ is a
$(1,1)$-form on $S_\Fa$.  Since the normal bundle to the divisor $\hat M_\Fa$
in $\hat M$ must have degree $-1$ on each fibre of $S_\Delta\to S_{\Fa}$ and
$\Omega$ is the curvature of a connection on this degree $-1$ line bundle, we
must have $[-\Omega/2\pi]\in c_1(\cO(-1))$.

To show that $\Omega$ is the Fubini--Study metric in its K\"ahler class, we
take $v\in\Fa^0$, the interior of a codimension one face of $\Delta$, and
construct a symplectic slice, as in Lemma~\ref{sslice}, to a point $x$ in
$\mu^{-1}(v)$ projecting to the given fibre of $S_\Delta\to S_{\Fa}$. Thus a
neighbourhood of $\T\cdot x$ in $M$ is equivariantly symplectomorphic to a
neighbourhood $U$ of the zero section $0_N\cong \T\cdot x$ of the normal
bundle $N=\T\times_{\T_\Fa} ({\mathfrak t}_x^0\oplus V_x)\to \T\cdot x$, with
the obvious $\T$-action, and canonical symplectic form $\omega_0$. Pulling
back the K\"ahler structure of $M$, and restricting to the fibre $V_x$ at $x$,
gives a K\"ahler metric on a neighbourhood of the origin in $V_x$ with a rigid
hamiltonian circle action of $\T_\Fa$ and constant symplectic form. Observe
that the K\"ahler quotient $P(V_x)$ of $V_x\smallsetminus\{0\}$ by $\T_\Fa$ is
a fibre of $S_\Delta\to S_\Fa$.

Let $z=r^2/2$ be half the distance squared to the origin in $V_x$---which is
the momentum map of the $\T_\Fa$ action contracted with $u_\Fa\in \mathfrak
t_\Fa$. Then the K\"ahler structure on $V_x$ may be written
\begin{equation*}
g = z h + \frac{dz^2}{H(z)} + H(z)\theta^2, \qquad\qquad\omega = z \Omega +
dz\wedge\theta,
\end{equation*}
for some function $H(z)$, where $d\theta=\Omega$ and $\Omega$ is as before,
and $(h,\Omega)$ is independent of $z$ (the K\"ahler quotient depends
affinely on $z$ and degenerates at $z=0$). The vector field dual to $\theta$
generates the $S^1$ action, and this preserves $z$, so it is tangent to the
level surfaces of $z$ (which are spheres), and generates a (topological)
Hopf fibration of them. Now $z$ is a function of the geodesic distance to
$z=0$ (the geodesic distance is obtained by integrating $1/\sqrt{H(z)}$).
For smooth compactification at $z=0$, the metric on geodesic spheres must
have constant curvature when $z\to 0$.  Hence $(h,\Omega)$ must tend to the
Fubini--Study metric, so that $\theta$ tends to the standard connection as
$z\to 0$. Since $(h,\Omega)$ is independent of $z$, it is the Fubini--Study
metric.

The explicit form of the metric on $M^0$ easily follows from Lemma
\ref{l:rigid} and Proposition \ref{data}. Note that a similar formula can be
established on $M^0_\Fa=\mu^{-1}(\Fa^0)$ for any face $\Fa$, but an origin
needs to be chosen in $\Fa$ so that $\mu\restr{M_\Fa}$ can be considered to
take values in $(\mathfrak t/\mathfrak t_\Fa)^*=\mathfrak t_\Fa^0$.
\end{proof}

\begin{rem} \label{construct}
In the absense of blow-downs, $\Omega_0+\langle \mu, \boldsymbol\Omega\rangle$
is positive for all $\mu$ in $\Delta$, and for all $\Fa$, $S_\Fa=S_\Delta$,
$P_\Fa=P_\Delta/\T_\Fa$, with the induced K\"ahler metrics and connections;
then the data of this proposition clearly \emph{do} define (uniquely) a
K\"ahler metric on $M$ with a rigid hamiltonian action of $\T$. However, the
existence of the connection $\boldsymbol\theta$ implies {\it integrality
conditions} on the curvature form $\boldsymbol\Omega$, and the
compactification of the toric K\"ahler metric on $\cV_\Delta$ implies {\it
boundary conditions} on ${\bf G}$.
\end{rem}

When there are blow-downs, it is difficult to describe the data needed to
construct the K\"ahler metric on $M$, because of the family of fibrations
$S_{\tFa}\to S_\Fa$: the K\"ahler quotient metrics are related by pullback,
and the fibrations and pullbacks must commute with. Rather than attempt this
in full generality, we restrict attention to a special case, which is all we
shall need for the application to hamiltonian $2$-forms.

\subsection{Semisimple actions and the generalized Calabi construction}
\label{s:generalcalabi}

\begin{defn} A hamiltonian torus action is {\it semisimple} if for any
regular value $v$ of the momentum map $\mu$, the derivative with respect to
$\mu$ of the family $\omega_h$ of K\"ahler forms on the complex quotient $S$
is parallel and diagonalizable with respect to $\omega_h$ at $\mu=v$.
(Observe that $S$ is well defined, as a complex orbifold at least, for $\mu$
in the connected component $U_v$ of $v$ in the regular values, since the
gradient flow of $\mu$ is transitive on $U_v$.)
\end{defn}

Integrating this condition, we deduce that on any connected component of the
regular values of $\mu$, the corresponding K\"ahler quotient metrics
$\omega_h$ are simultaneously diagonal with the same Levi-Civita connections.
Thus for a semisimple rigid hamiltonian torus action, there is a symplectic
$(1,1)$-form $\Omega_S$ on $S_\Delta$ such that the family of K\"ahler forms
induced by $\mu\in\Delta^0$ are parallel and simultaneously diagonalizable
with respect to $\Omega_S$.

\begin{defn}\label{gcc}
By {\it generalized Calabi data} of dimension $m$, rank $\ell$, we mean:
\begin{numlist}
\item a $2(m-\ell)$-dimensional product $S$ of $N\geq 0$ K\"ahler manifolds
$(S_a,\pm g_a,\pm\omega_a)$ of dimension $2m_a>0$ (if $\ell=m$, $N=0$);
\item a compact toric $2\ell$-dimensional K\"ahler manifold $\cV$ with Delzant
polytope $\Delta\subset\mathfrak t^*$ and momentum map $\mu_\cV\colon
\cV\to\Delta$;
\item a principal $\T$-bundle $P\to S$, with a principal connection of
curvature $\boldsymbol\Omega\in C^\infty(S,\Lambda^{1,1}S\otimes\mathfrak t)$,
where $\T$ is the $\ell$-torus acting on $\cV$;
\item a $(1,1)$-form $\Omega_0$ on $S$ such that $\Omega_0 + \langle v,
\boldsymbol\Omega\rangle$ is positive for $v\in\Delta^0$;
\item constants $c_{a0}\in\R$ and $\boldsymbol c_a\in\mathfrak t$ such that
$\Omega_0=\sum_{a=1}^N c_{a0}\omega_a$ and
$\boldsymbol\Omega=\sum_{a=1}^N\boldsymbol c_a\omega_a$;
\item a subset $\cC\subset\{1,\ldots N\}$ such that for $a\notin\cC$,
$\{v\in\Delta:c_{a0}+\langle v,\boldsymbol c_a\rangle=0\}$ is empty, while
for $a\in \cC$ they are distinct codimension one faces of $\Delta$ with
(primitive) inward normals $u_a\in\mathfrak t$, and $S_a = \C P^{d_a}$ with
$d_a>0$, $\pm g_a$ is a Fubini--Study metric and $\boldsymbol
c_a\otimes\omega_a/2\pi\in u_a\otimes c_1(\cO(-1))$.
\end{numlist}
Given these data we define the manifold $\hat M=P\times_{\T}\cV=
M^0\times_{\T^c}\cV\to S$, where $M^0=P\times_{\T}\mu^{-1}_\cV(\Delta^0)$.
Since the curvature $2$-form of $P$ has type $(1,1)$, $M^0$ becomes a
holomorphic principal $\T^c$-bundle with connection and $\smash{\hat M}$ is a
complex manifold. The toric K\"ahler structure on $\cV$ endows $\smash{\hat
M}$ with a fibrewise metric and `momentum map' $\hat\mu\colon\smash{\hat
M}\to\Delta$: indeed, being $\T$ invariant, the momentum map $\mu_\cV$ of
${\cV}$ can be defined on ${\hat M}=P\times_{\T}\cV$.

According to (vi), the set $\cC$ corresponds bijectively to a subset $\cB$ of
the codimension one faces of $\Delta$, and for $\Fa \in {\mathcal B}$
corresponding to $a\in\cC$, the connection on $\hat M_\Fa:=\hat\mu ^{-1}(\Fa)$
is flat over each fibre of $S\to\prod_{\smash{b\neq a}} S_b$.  This gives a
$\C P^{d_a}$ fibration of $\hat M_\Fa$ such that the normal bundle to $\hat
M_\Fa$ in $\hat M$ is a line bundle which has degree $-1$ on each $\C P^{d_a}$
fibre. Since a tubular neighbourhood of $\hat M_\Fa$ in $\hat M$ is
diffeomorphic to a neighbourhood of the zero section in the normal bundle, it
follows that the topological space $M$, obtained by contracting $\hat M$ along
the $\C P^{d_a}$ fibration of each such $\hat M_\Fa$, is a smooth manifold and
$M^0$ is an open dense submanifold.

If the K\"ahler structure given by~\eqref{hkf} (which pulls back to the
fibrewise metric on the fibres of $\smash{\hat M}\to S$) extends smoothly to
$M$, then we say that this K\"ahler manifold $(M,g,J,\omega)$ is given by {\it
the generalized Calabi construction} (with blow-downs).
\end{defn}

We shall see that the contraction $\hat M\to M$ realises $\hat M$ (with the
complex structure described above) as a blow-up of $M$. We therefore refer to
this contraction as a blow-down. Our main result shows that all generalized
Calabi data give rise to a generalized Calabi construction, and that this
classifies compact K\"ahler manifolds with a semisimple rigid hamiltonian
torus actions up to a covering.

\begin{thm}\label{generalizedcalabi} Let $M$ be a compact connected K\"ahler
$2m$-manifold with a semisimple rigid hamiltonian action of an $\ell$-torus
$\T$ and momentum map $\mu\colon M\to \Delta \subset{\mathfrak t}^*$. Then
some cover of $M$ is given by the generalized Calabi construction.

Conversely, for any generalized Calabi data \textup{(i)--(vi)}, the
generalized Calabi construction produces a smooth K\"ahler manifold with a
semisimple rigid hamiltonian action of an $\ell$-torus.
\end{thm}
\begin{proof} We construct the generalized Calabi data from
Proposition~\ref{data}, imposing the condition that the action is semisimple.
As remarked in~\cite[\S 3.3]{ACG1}, the condition that $\Omega_0$ and the
components of $\boldsymbol\Omega$ are simultaneously diagonalizable and
parallel (with respect to some K\"ahler metric $\Omega_S$) implies that the
(distinct) eigendistributions $\cH_a$ $(a=1,\ldots N)$ are parallel. By the
deRham decomposition theorem, some cover of $S_\Delta$ (for instance the
universal cover), is a K\"ahler product $(S,\Omega_S)=\prod_{a=1}^N
(S_a,\omega_a)$ (note that $S$ may not be compact). The generalized Calabi
data (i)--(v) are then obtained from Proposition~\ref{data} by setting
$\cV=\cV_\Delta$, pulling back $P_\Delta$, $\theta_\Delta$, $\Omega_0$ and
$\boldsymbol\Omega$ to give a principal bundle $P$ with connection over $S$,
and defining the constants $c_{a0}$ and $\boldsymbol c_a$ by (v).

Let $\cB$ be the set of codimension one faces $\Fa$ of $\Delta$ such that a
blow-down occurs (i.e., $M_{\Fa}$ is not a divisor); then $\hat M$ is the
blow-up of $M$ along $M_{\Fa}$ with $\Fa\in\cB$. The pullback of the metric to
$\hat M$ degenerates on the fibres of a $\C P^{d}$-bundle $S_\Delta\to
S_{\Fa}$ for some $d>0$. Now $\C P^d$ is simply connected, so this is covered
by a $\C P^d$-bundle with total space $S$, whose base is a cover of $S_\Fa$.
Hence there must be at least one $a$ such that $c_{a0}+\langle v,\boldsymbol
c_a\rangle=0$ for $v\in\Fa^0$; since $\C P^d$ does not admit a K\"ahler
product metric, this $a$ is unique, and $S_\Fa$ is covered by $\prod_{b\neq a}
S_b$, while $S_a = \C P^{d_a}$ with $d_a=d$.  On the other hand
$c_{a0}+\langle v,\boldsymbol c_a\rangle$ is an affine function of $v$, so it
can vanish on at most one codimension one face of the Delzant polytope
$\Delta$.  Thus $\cB$ corresponds bijectively to a subset
$\cC\subset\{1,\ldots N\}$. Now note that for any face $\Fa$, with
$v\in\Fa^0$, the metric induced on $S_\Fa$ is nondegenerate, so
$c_{a0}+\langle v,\boldsymbol c_a\rangle$ does not vanish on $\Delta$ for
$a\notin\cC$. This establishes (vi).

The pullback of $\hat M$ to $S$ is a cover of $\hat M$, and by construction,
this descends to $M$. Hence, up to a cover, $M$ is obtained from the
generalized Calabi construction.

\smallbreak

Conversely, given the data of Definition~\ref{gcc}, we will prove that there
exists a smooth compact K\"ahler manifold $(M,g,J,\omega)$ with a semisimple
rigid hamiltonian action of the $\ell$-torus $\T$ given by the generalized
Calabi construction. The main difficulty is to deal with the blow-downs.

Let us suppose there are $k\geq 0$ blow-downs: then, after reordering, we may
assume ${\mathcal C} = \{1,\ldots k\}$ and that $S = {\C}P^{d_1} \times
\cdots \times {\C}P^{d_k} \times S''$ for some K\"ahler product $S''$.  The
conditions (iii) and (v) of Definition~\ref{gcc} imply that
$\boldsymbol\Omega'':=\sum_{a=k+1}^N\boldsymbol c_a\omega_a$ is the curvature
of a principal $\T$-bundle $P'' \to S''$. We are going to let $M$ be of the
form $P''\times_{\T} M'$, where $M'$ is a $2(\ell +d_1 + \cdots + d_k)$
dimensional K\"ahler manifold with a rigid semisimple isometric hamiltonian
action of $\T$, obtained from the generalized Calabi construction with respect
to the following data:
\begin{numlist}
\item $S'= {\C}P^{d_1} \times \cdots \times {\C}P^{d_k}$; 
\item $(\cV, \omega, \mu, \Delta)$, with the given compatible toric K\"ahler
metric;
\item a principal $\T$-bundle $P'$ with curvature form ${\boldsymbol \Omega}'
= \sum_{a=1}^{k}{\boldsymbol c}_{a} \omega_a$;
\item $\Omega_0' = \sum_{a=1}^k c_{a0}\omega_{a}$;
\item the given constants $c_{a0}$ and ${\boldsymbol c}_a$ for $a=1,\ldots
k$;
\item $\cC = \{1,\ldots k \}$.
\end{numlist}
Since the data for $M$ are generalized Calabi data, so are these data for
$M'$.  If $M'$ can be constructed with these data, it follows from
Proposition~\ref{data} and Remark~\ref{construct} that $M$ is equipped with a
K\"ahler metric and a semisimple rigid action of $\T$; using the first part of
the theorem we also see that $M$ is given by the generalized Calabi
construction associated to the initial data.

Thus it remains only to establish the generalized Calabi construction for
$M'$. However, such an $M'$ is obtained as a restricted toric K\"ahler
manifold, the construction of which we discussed already in
\S\ref{s:restricttoric}.
\end{proof}

Just as toric complex manifolds may be described in terms of linear charts,
i.e., in terms of a family of vector spaces, each with a decomposition into
one dimensional subspaces, glued together by Laurent monomials, so bundles of
toric complex manifolds (arising in the generalized Calabi construction
without blow-downs) may be decribed (by the holomorphic slice theorem) in
terms of families of vector bundles, each a direct sum of line bundles, glued
together in a similar way.  The simplest case is the case of projective
bundles $P(\cL_0\oplus\cL_1\oplus\cdots\oplus \cL_\ell)\to S$, which are
obtained by gluing together the vector bundles $\bigoplus_{k\neq j} \cL_k$ for
$j=0,\ldots \ell$. This is the only case we shall need in the sequel.

\section{Orthotoric geometry}\label{s:orthotoric}

We now return to our primary aim: the classification of compact K\"ahler
manifolds endowed with a hamiltonian $2$-form. In this section we treat the
case when the order of the hamiltonian $2$-form is maximal, and therefore the
corresponding K\"ahler manifolds are toric. Motivated by the orthogonality of
the gradients of the roots of the momentum polynomial, see Theorem
\ref{ACGthm}, we define orthotoric K\"ahler manifolds and orbifolds, and
classify the compact ones.

\subsection{The polytope of an orthotoric orbifold}
\label{s:orthopoly}

\begin{defn} An \emph{orthotoric} K\"ahler manifold (or orbifold)
$M$ is a toric K\"ahler $2m$-manifold (or orbifold) with a momentum map
${\boldsymbol\sigma}=(\sigma_1,\ldots \sigma_m)$ and (rational) Delzant
polytope $\Delta = {\boldsymbol \sigma}(M)$, such that on the dense open set
$M^0 = {\boldsymbol \sigma}^{-1}(\Delta^0)$ of regular points of ${\boldsymbol
\sigma}$, the roots $\xi_1,\ldots\xi_m$ of the momentum polynomial
$\sum_{r=0}^m (-1)^r\sigma_r t^{m-r}$ ($\sigma_0=1$) are smoothly defined,
pairwise distinct and functionally independent, and the K\"ahler metric has
the explicit form
\begin{equation}\label{orthotoric}\begin{split}
g &= \sum_{j=1}^m \frac{\Delta_j}{\oF_j (\xi_j)}\, d\xi_j^2
+ \sum_{j=1}^m \frac{\oF_j (\xi_j)}{\Delta_j}
\biggl(\sum_{r=1}^m \sigma_{r-1}(\hat{\xi}_j)\, dt_r \biggr)^2\\
&= \sum_{r,s,j=1}^m\biggl(
\frac{(-1)^{r+s}\Delta_j\xi_j^{2m-r-s}}{\oF_j(\xi_j)}d\sigma_r\,d\sigma_s
+\frac{\oF_j(\xi_j)\sigma_{r-1}(\hat{\xi}_j)\sigma_{s-1}(\hat{\xi}_j)}
{\Delta_j}dt_r\,dt_s\biggr)\\
\omega &= \sum_{j=1}^m d \xi_j \wedge
\biggl(\sum_{r=1}^{m} \sigma_{r-1}(\hat{\xi}_j)\, dt_r \biggr )
= \sum_{r=1}^{m} d\sigma_r \wedge d t_r,
\end{split}\end{equation}
for functions $\oF_1,\ldots\oF_m$ of one variable.  Here
$\Delta_j=\prod_{k\neq j} (\xi_j-\xi_k)$.
\end{defn}
Clearly the gradients of $\xi_1,\ldots\xi_m$ are orthogonal with respect to
$g$. Conversely, it was shown in~\cite[\S 3.4]{ACG1} that this property
characterizes orthotoric K\"ahler manifolds (and the result applies equally to
orbifolds).

Note that the basis $K_1,\ldots K_m$ of the Lie algebra of the torus
identifies it with $\R^m$, and we view the invariant $1$-forms $dt_1,\ldots
dt_m$ as the dual basis of $\R^{m*}$.

\begin{prop}\label{orthopoly} Let $M$ be a compact orthotoric K\"ahler
$2m$-manifold or orbifold with momentum map $\boldsymbol\sigma=
(\sigma_1,\ldots \sigma_m)$ and rational Delzant polytope $\Delta$.
\begin{numlist}
\item $\Delta$ is the \textup(one to one\textup) image under the elementary
symmetric functions of a domain of the form
\begin{align}
D&=\{(\xi_1,\ldots \xi_m)\in\R^m:\alpha_{j}\leq\xi_j\leq\beta_j\}\\
\tag*{\it where} &\alpha_1<\beta_1\leq
\alpha_2<\beta_2\leq\cdots<\beta_{m-1}\leq \alpha_m<\beta_m.
\end{align}
Thus, setting $\sigma_0=1$, $\Delta=\{(\sigma_1,\ldots\sigma_m):
(-1)^{m-j}\sum_{r=0}^m (-1)^r \sigma_r\alpha_j^{m-r}\leq 0 \ {\rm and}$
$(-1)^{m-j}\sum_{r=0}^m (-1)^r \sigma_r\beta_j^{m-r}\geq 0 \ {\rm for}\
j=1,\ldots m \}$. This is a simplex if and only if $\alpha_{j+1}=\beta_{j}$
for $j=1,\ldots m-1$.

\item If $M$ is nonsingular \textup(i.e., a manifold\textup), then
$\Delta$ is a simplex.
\end{numlist}
\end{prop}
\begin{numlproof}
\item $\sigma_1,\ldots\sigma_m$ are the elementary symmetric functions of the
roots $\xi_1,\ldots\xi_m$ of the momentum polynomial, and we want to find the
domain $D$ in the $\xi_j$ coordinates corresponding to $\Delta$.  We first
remark that this domain must be bounded.  Also, the functions $\oF_j(\xi_j)$
must be nonzero on the interior $D^0$ of $D$ in order that the metric be
finite and nondegenerate.

Now consider in particular the metric on the torus given by
\begin{equation*}
\g{K_r,K_s} = \sum_{j=1}^m \frac{\oF_j(\xi_j)
\,\sigma_{r-1}(\hat \xi_j)\sigma_{s-1}(\hat\xi_j)}{\Delta_j}.
\end{equation*}
The determinant of this matrix is (up to a sign) $\prod_{j=1}^m
\oF_j(\xi_j)$. As we approach a special orbit of the $m$-torus action, i.e.,
as $\boldsymbol\sigma$ approaches the boundary of $\Delta$, this must tend to
zero, i.e., at least one of the functions $\oF_j$ of one variable must tend to
zero. Since these functions are nonvanishing on $D^0$, it follows that $D^0$
is a domain of the form $\prod_{j=1}^m (\alpha_j,\beta_j)$, where $\oF_j$ is
nonvanishing on the interval $(\alpha_j,\beta_j)$ and tends to zero at the
endpoints. Now $\xi_1,\ldots\xi_m$ must be pairwise distinct on $D^0$, so we
may assume (after reordering) that $\xi_1<\cdots<\xi_m$ on $D^0$.  Hence
\begin{equation*}
\alpha_1<\beta_1\leq \alpha_2<\beta_2\leq\cdots<\beta_{m-1}\leq
\alpha_m<\beta_m.
\end{equation*}
Noting that the elementary symmetric functions are affine in each variable, we
readily check that this domain does indeed map bijectively to a convex
polytope. Indeed any $(\xi_1,\ldots\xi_m)$ in $D$ satisfy
\begin{equation}\label{deltabound}
(-1)^{m-j} \prod_{k=1}^m (\alpha_j-\xi_k) \leq 0,
\qquad (-1)^{m-j} \prod_{k=1}^m (\beta_j-\xi_k)\geq 0
\end{equation}
for all $j=1,\ldots m$; equality is attained in one of these expressions on
any face, and in any of these expressions on some face. Expanding in terms of
the elementary symmetric functions of $(\xi_1,\ldots\xi_m)$ gives the explicit
description of $\Delta$.

A compact convex polytope in $\R^{m*}$ is a simplex if and only if it has
$m+1$ vertices. The vertices of $D$ are the points where
$\xi_j\in\{\alpha_j,\beta_j\}$ for all $j=1,\ldots m$. Now observe that a
vertex of $D$ maps to a vertex of $\Delta$ if and only if it does not lie on
one of the diagonals $\xi_j=\xi_k$ for $j\neq k$.

\item We shall show that $\alpha_{j+1}=\beta_j$ for $j=1,\ldots m-1$.  Suppose
for contradiction that this does not hold for some $j\in\{1,\ldots m-1\}$ and
consider the four vertices
\begin{gather*}\notag
(\alpha_1,\ldots \alpha_{j-1},
\alpha_j,\alpha_{j+1},\alpha_{j+2},\ldots \alpha_m),\qquad
(\alpha_1,\ldots \alpha_{j-1},
\alpha_j,\beta_{j+1},\alpha_{j+2},\ldots \alpha_m),\\
(\alpha_1,\ldots \alpha_{j-1},
\beta_j,\alpha_{j+1},\alpha_{j+2},\ldots \alpha_m),\qquad
(\alpha_1,\ldots \alpha_{j-1},
\beta_j,\beta_{j+1},\alpha_{j+2},\ldots \alpha_m).
\end{gather*}
Since $\alpha_j<\beta_j<\alpha_{j+1}<\beta_{j+1}$ these four points map to
four distinct vertices spanning a two dimensional face of $\Delta$
($\xi_k=\alpha_k$ defines a hyperplane).  Now any face of a Delzant polytope
is Delzant (as one easily checks) and the Delzant property is invariant under
affine transformation. Hence we may as well map this $2$-dimensional face into
$\R^2$ by sending $(\sigma_1,\ldots \sigma_m)$ to
$(\sigma_1-a_1,\sigma_2-a_1\sigma_1+a_1^2-a_2)$, where $a_1$ and $a_2$ are the
first two elementary symmetric functions of $\{\alpha_k:k\neq j,j+1\}$: in
terms of $\xi_j,\xi_{j+1}$ (fixing $\xi_k=\alpha_k$ for $k\neq j,j+1$), this
formula gives $(\xi_j+\xi_{j+1},\xi_j\xi_{j+1})$, and so our face gets mapped
to the quadrilateral with vertices
\begin{equation*}\notag
(\alpha_j + \alpha_{j+1}, \alpha_j\alpha_{j+1}),\;
(\alpha_j + \beta_{j+1}, \alpha_j\beta_{j+1}),\;
(\beta_j + \alpha_{j+1}, \beta_j\alpha_{j+1}),\;
(\beta_j + \beta_{j+1}, \beta_j\beta_{j+1})
\end{equation*}
and normals (up to scale)
\begin{equation*}
(\alpha_j,-1), \quad (\beta_j,-1),\quad (\alpha_{j+1},-1), \quad
(\beta_{j+1},-1).
\end{equation*}
Again, $\alpha_j<\beta_j<\alpha_{j+1}<\beta_{j+1}$, so these four normals
point in distinct directions, and so cannot be scaled to form a basis for the
same lattice at each vertex. Our quadrilateral is therefore not Delzant, hence
neither is $\Delta$, a contradiction.
\end{numlproof}

For the rest of this subsection we suppose $\Delta$ is a simplex: the above
proposition shows that this is is necessarily true if $M$ is nonsingular.

By the Delzant construction, any symplectic orbifold $M$ whose rational
Delzant polytope is a simplex is a symplectic quotient of $\C^{m+1}$ by a one
dimensional subgroup $G$ of $(S^1)^{m+1}$. From the relation between complex
and symplectic quotients, cf.~\eqref{complexquot}, it follows that $M$ is a
quotient of a {\it weighted projective space} $\C P^{m}_{a_0, \ldots
a_{m}}$---here $a_0,\ldots a_{m}\in\Z^+$ have highest common factor $1$ and
$\C P^{m}_{a_0, \ldots a_{m}}$ is the quotient of $\C^{m+1} \smallsetminus \{
0 \}$ by the holomorphic action
\begin{equation*}
(z_0, \ldots z_{m}) \to (\zeta^{a_0}z_0, \ldots \zeta^{a_{m}}z_{m})
\quad {\rm for} \quad\zeta \in \C^{\times};
\end{equation*}
note that $\C P^{m}_{1, \ldots 1}$ is the usual (nonsingular) $\C P^{m}$.

We want to describe $\Delta$ more explicitly as a rational Delzant simplex.
We put $\beta_0=\alpha_1$, so $\Delta$ is the image under the elementary
symmetric functions of the domain
\begin{align}\label{simplexD}
D&=\{(\xi_1,\ldots \xi_m)\in\R^m: \beta_{j-1}\leq\xi_j\leq\beta_j\}\\
\tag*{where} &\beta_0<\beta_1<\cdots<\beta_{m-1}<\beta_m.
\end{align}

\begin{prop}\label{orbi-simplex}
Let $M$ be a compact orthotoric K\"ahler $2m$-orbifold whose Delzant polytope
$\Delta$ is the image of~\eqref{simplexD} under the elementary symmetric
functions.
\begin{numlist}
\item $\Delta = \{\boldsymbol\sigma:\langle v_j,\boldsymbol\sigma\rangle
+\kappa_j\geq 0\}$, where $\kappa_j={\beta_j^{m}}/{\prod_{k\neq j}
(\beta_j-\beta_k)}$ and
\begin{equation}\label{vj}
v_j=\biggl(\frac{-\beta_j^{m-1}}{\prod_{k\neq j}(\beta_j-\beta_k)},
\ldots \frac{(-1)^r\beta_j^{m-r}}{\prod_{k\neq j}(\beta_j-\beta_k)},
\ldots \frac{(-1)^m}{\prod_{k\neq j}(\beta_j-\beta_k)}\biggr).
\end{equation}
The codimension one faces of $\Delta$ are $\Fa_0,\ldots \Fa_m$, where
\begin{bulletlist}
\item $\Fa_0$ is the image of the boundary component $\xi_1=\beta_0$ of $D$,
\item $\Fa_m$ is the image of the boundary component $\xi_m=\beta_m$ of $D$,
and
\item $\Fa_j$, for $j=1,\ldots m-1$, is the union of the images of the
boundary components $\xi_j=\beta_{j}$ and $\xi_{j+1}=\beta_j$ of $D$.
\end{bulletlist}

\item The normals are of the form $u_j=2 n_j v_j/c$, where $c>0$ and
$n_j\in\Z^+$ $(j=0,\ldots m)$ have highest common factor $1$\textup; then $M$
is equivariantly biholomorphic to an orbifold quotient of $\C
P^{m}_{a_0,\ldots a_m}$, where $n_j=\prod_{k\neq j} a_k$.

\item $M$ is nonsingular if and only if it is biholomorphic to $\C P^m$ if and
only if $n_j=1$ \textup(for all $j$\textup) and the lattice of circle
subgroups is generated by $u_0,\ldots u_m$.  The dual lattice in $\R^{m*}$ is
then generated by
\begin{equation}\label{duallat}
\theta_{p,q} = \sum_{r=0}^m \tfrac 1 2 \cc\bigl(\sigma_{r}^{\smash\beta}(\hat
\beta_q)-\sigma_{r}^{\smash\beta}(\hat \beta_p)\bigr) dt_r
\end{equation}
where $\sigma_{r}^{\smash\beta}(\hat \beta_p)$ denotes the $r$th elementary
symmetric function of the $m$ variables $\{\beta_j:j=0,\ldots m,\, j\neq p\}$.
\end{numlist}
\end{prop}
\begin{numlproof}
\item When $\Delta$ is a simplex, the inequalities in~\eqref{deltabound} may
be written
\begin{equation*}
\frac{\prod_{k=1}^m (\beta_j-\xi_k)}
{\prod_{k\neq j} (\beta_j-\beta_k)}\geq 0
\end{equation*}
for all $j=0,\ldots m$, which immediately gives the stated form of $\Delta$.
(Note that the apparent codimension two face $\xi_j=\beta_{j}=\xi_{j+1}$ is
`straightened out' by the elementary symmetric functions; this is why $\Delta$
has only $m+1$ faces, not $2m$.)

\item From the form of the simplex $\Delta$ it is immediate that the normals
$u_0,\ldots u_m$ are positive multiples of $v_0,\ldots v_m$.  They belong to a
common lattice if and only if the linear dependence relation among them can be
written $\sum_{j=0}^m u_j/n_j =0$, where $n_0,\ldots n_m$ are nonzero rational
numbers. We now observe that the $v_j$'s already satisfy $\sum_{j=0}^{m}
v_j=0$ by the Vandermonde identity (cf.~\cite[Appendix B]{ACG1}). Hence we
must have $u_j=C n_j v_j$ for some nonzero constant $C$ and without loss of
generality we can take $C$ and $n_j$'s to be positive and suppose $n_0,\ldots
n_m$ are integers with highest common factor $1$. We then put $C=2/c$.

We have already seen that any toric K\"ahler orbifold with polytope a simplex
is equivariantly biholomorphic to an orbifold quotient of a weighted
projective space. It remains to show that the integers $n_j$ are related to
the weights $a_k$ by $n_j=\prod_{k\neq j} a_k$.  For this, we
note~\cite{Abreu1} that any weighted projective space has an orbifold quotient
whose simplex is standard with respect to the lattice $\Lambda$, i.e., the
primitive normals sum to zero. The primitive normals are $u_j/m_j$ and Abreu
shows that the labels (in this case) are given by $m_j=\prod_{k\neq j}
a_k$. Since $\sum_{j=0}^m v_j=0$, and the $m_j$ have highest common factor
$1$, we have $m_j=n_j$.

\item The only (orbifold quotient of a) weighted projective space which is
nonsingular is $\C P^m$.  Clearly $M$ is equivariantly biholomorphic to $\C
P^m$ if and only if the $n_j$ all equal $1$ and the lattice $\Lambda$ of
circle subgroups is the minimal one.  In terms of the vector fields
$K_1,\ldots K_m$, it follows that vector fields generating $\Lambda$ are
\begin{equation}\label{euler}
X_j = \frac{2}{\cc} \sum_{r=1}^m
\frac{(-1)^r\beta_j^{m-r} K_r}{\prod_{k\neq j}(\beta_j-\beta_k)},
\end{equation}
with $\sum_{j=0}^m X_j=0$. To see that~\eqref{duallat} generate the dual
lattice, we note that
\begin{equation*}
\theta_{p,q} =  \sum_{r=1}^m \tfrac 1 2 \cc\,
\sigma_{r-1}^{\smash\beta}(\hat \beta_p, \hat\beta_q) (\beta_p-\beta_q) dt_r
\end{equation*}
for $0\leq p<q\leq m$, where $\sigma_{r-1}^{\smash\beta}(\hat \beta_p,
\hat\beta_q)$ is the $(r-1)$st elementary symmetric function of the
$m-1$ variables $\{\beta_j:j=0,\ldots m, j\neq p,q\}$.  We compute
\begin{align*}
\theta_{p,q}(X_j)&=\sum_{r=1}^m \frac{(-1)^r\sigma_{r-1}^{\smash\beta}(\hat
\beta_p, \hat\beta_q)\, (\beta_p-\beta_q)\beta_j^{m-r}}{\prod_{k\neq
j}(\beta_j-\beta_k)}\\
&=\frac{\prod_{k\neq p,q}(\beta_j-\beta_k)}{\prod_{k\neq
j}(\beta_j-\beta_k)} (\beta_q-\beta_p) = \delta_{jq}-\delta_{jp}
\end{align*}
and the result follows.
\end{numlproof}

The constant $c$ determines the scale of $M$: the symplectic volume is
proportional to $1/c$. The other constants $\beta_0,\ldots \beta_m$ are
related to the fact that the Killing vector fields $K_1,\ldots K_m$ do not
necessarily form an integral basis.

We remark that all simplices are equivalent under affine transformation, and
so for any $\beta_0<\cdots< \beta_m$, any rational Delzant simplex is
equivalent to the simplex of this proposition for some lattice $\Lambda$ in
$\R^{m}$ and some normals $u_j=2n_j v_j/c \in\Lambda$.

\subsection{Compactification of orthotoric K\"ahler metrics}

We next establish necessary and sufficient conditions for the compactification
of the orthotoric K\"ahler metric \eqref{orthotoric} on a compact
$2m$-orbifold $M$. We obtain these conditions by specializing those of
Proposition~\ref{first-order-conditions} to the orthotoric case.

\begin{prop}\label{orthocompact}
Let $M$ be a compact symplectic $2m$-orbifold such that the rational Delzant
polytope $\Delta\subset \R^{m*}$ is the image of $\prod_{j=1}^m
[\alpha_j,\beta_j]$ under the elementary symmetric functions, where
\begin{equation*}
\alpha_1<\beta_1\leq\alpha_2<\beta_2\leq\cdots
<\beta_{m-1}\leq\alpha_m<\beta_m.
\end{equation*}
Let $L_j^\alpha(\boldsymbol\sigma)=\langle u^\alpha_j,\boldsymbol\sigma\rangle
+ \lambda^{\alpha}_j$ and $L_j^{\smash\beta}(\boldsymbol\sigma)=\langle
u^{\smash\beta}_j,\boldsymbol\sigma \rangle + \lambda^{\smash\beta}_j$ where
\begin{align*}
\lambda_j^{\alpha}&= -c_j^{\alpha}\alpha_j^{m}, &
u_j^{\alpha} &= c_j^{\alpha} \bigl(\alpha_j^{m-1}, \ldots
(-1)^{r-1}\alpha_j^{m-r}, \ldots (-1)^{m-1}\bigr),\\
\lambda_j^{\smash\beta}&= -c_j^{\smash\beta}\beta_j^{m}, &
u_j^{\smash\beta} &= c_j^{\smash\beta} \bigl(\beta_j^{m-1}, \ldots
(-1)^{r-1}\beta_j^{m-r}, \ldots (-1)^{m-1}\bigr),
\end{align*}
and the constants $c_j^{\alpha},c_j^{\smash\beta}\in\R$ are such that the
normals of $\Delta$ are the distinct elements among
$u_j^{\alpha},u_j^{\smash\beta}$, i.e., $\Delta = \{\boldsymbol\sigma \in
\R^{m*}: L_j^\alpha(\boldsymbol\sigma)\geq 0 \ {\rm and}\
L_j^{\smash\beta}(\boldsymbol\sigma)\geq 0$ {\rm for} $j=1,\ldots m\}$, but if
$\alpha_{j+1}=\beta_{j}$, we have $c_{j+1}^{\alpha}=c_j^{\smash\beta}$ as the
normals $u_{j+1}^{\alpha},u_j^{\smash\beta}$ are then not distinct.

Then the K\"ahler metric~\eqref{orthotoric}, defined for $\xi_j\in
(\alpha_j,\beta_j)$, extends to an orthotoric K\"ahler metric on $M$ if and
only if for $j=1,\ldots m$, $\Theta_j$ is the restriction to
$(\alpha_j,\beta_j)$ of a smooth function $\Theta$ on $\bigcup_{j=1}^m
[\alpha_j,\beta_j]$ satisfying \textup(for $j=1,\ldots m$\textup)\textup:
\begin{gather}
\label{bound}\begin{split}
\Theta(\alpha_j) &= 0 = \Theta(\beta_j),\\
\Theta'(\alpha_j) c_j^{\alpha} &= 2 = \Theta'(\beta_j) c_j^{\smash\beta};
\end{split}\\
(-1)^{m-j}\Theta>0\quad\text{on}\quad (\alpha_j,\beta_j).\label{positivity}
\end{gather}
\end{prop}

\begin{proof} By~\eqref{orthotoric}, ${\bf H}$ is given by
\begin{equation*}
H_{rs} = \sum_{j=1}^m \frac{\oF_j(\xi_j)
\,\sigma_{r-1}(\hat \xi_j)\sigma_{s-1}(\hat\xi_j)}{\Delta_j}.
\end{equation*}
This is a smooth and symmetric function of $\xi_1,\ldots\xi_m$, so by
Glaeser~\cite{glaeser}, it is a smooth function of $\sigma_1,\ldots\sigma_m$.
The positivity condition is clear, so it remains to consider the boundary
conditions~\eqref{bound}. We must show these are equivalent
to~\eqref{eq:toricboundary}.

The form of the normals shows that ${\bf H}(u_i^\alpha,\cdot)$ is given by
\begin{align*}
\sum_{r=1}^m H_{rs} (u_i^\alpha)_r &= \sum_{j,r=1}^m
\frac{c_i^{\alpha} \oF_j(\xi_j)
\,(-1)^{r-1} \sigma_{r-1}(\hat \xi_j)\alpha_i^{m-r}
\sigma_{s-1}(\hat\xi_j)}{\Delta_j}\\
&= \sum_{j=1}^m \frac{c_i^{\alpha} \oF_j(\xi_j) \sigma_{s-1}(\hat\xi_j)
\prod_{k\neq j} (\alpha_i-\xi_k)}{\Delta_j}.
\end{align*}
On the codimension one face $\xi_i=\alpha_i$, this reduces to $c_i^{\alpha}
\oF_i(\alpha_i) \sigma_{s-1}(\hat\xi_i)$, which vanishes for all $s$ if and
only if $\oF_i(\alpha_i)=0$. For the derivative conditions we differentiate
\begin{equation*}
{\bf H}(u_i^\alpha,u_i^\alpha)=\sum_{j=1}^m \frac{(c_i^{\alpha})^2
\oF_j(\xi_j) \prod_{k\neq j} (\alpha_i-\xi_k)^2}{\Delta_j}
\end{equation*}
and evaluate along $\xi_i=\alpha_i$ to obtain
\begin{equation*}
(c_i^{\alpha})^2 \oF_i'(\alpha_i) \prod_{k\neq i}
(\alpha_i-\xi_k)\,d\xi_i
= c_i^{\alpha}\oF_i'(\alpha_i) d\prod_{k=1}^m
c_i^\alpha(\alpha_i-\xi_k)\bigg|_{\xi_i=\alpha_i}.
\end{equation*}
This equals $2u_i^\alpha$ if and only if $c_i^\alpha\oF_i'(\alpha_i)=2$.
The boundary conditions at the $\beta$ endpoints are analogous.
\end{proof}

Note that~\eqref{bound} could be taken as the definition of the constants
$c_j^{\alpha}$ and $c_j^{\smash\beta}$.  However, these are then required to
satisfy positivity and integrality conditions, since $(-1)^{m-j}c_j^{\alpha}$
and $(-1)^{m-j+1}c_j^{\smash\beta}$ must be positive for the normals to be
inward pointing, while $u_j^{\alpha}$ and $u_j^{\smash\beta}$ must belong to a
common lattice in $\R^m$.

\smallbreak

We summarize our results for the case that the rational Delzant polytope is a
simplex.  The following is immediate from Propositions~\ref{orthopoly},
\ref{orbi-simplex} and~\ref{orthocompact}.

\begin{thm}\label{CPorthotoric}
Let $M$ be a compact orthotoric $2m$-manifold or orbifold with momentum map
$\boldsymbol\sigma$, whose rational Delzant polytope is a simplex $\Delta$
with normals $u_0,u_1,\ldots u_m\in\R^m$, where the K\"ahler metric is given
by~\eqref{orthotoric} on $M^0=\boldsymbol\sigma ^{-1}(\Delta ^0)$.

\begin{numlist}
\item $M$ is equivariantly biholomorphic to a toric orbifold quotient of $\C
P^m_{a_0,\ldots a_m}$ and, with $n_j=\prod_{k\neq j} a_k$, there are
constants $\beta_0<\beta_1<\cdots <\beta_m$, $\cc>0$, and a smooth function
$\Theta$ on $[\beta_0,\beta_{m}]$, such that for $j=0,\ldots m$\textup:
\begin{gather}\label{CPnormals}
u_j=\frac{2n_j}{\cc}
\biggl(\frac{-\beta_j^{m-1}}{\prod_{k\neq j}(\beta_j-\beta_k)},\ldots
\frac{(-1)^r\beta_j^{m-r}}{\prod_{k\neq j}(\beta_j-\beta_k)},\ldots
\frac{(-1)^m}{\prod_{k\neq j}(\beta_j-\beta_k)}\biggr); \displaybreak[0]\\
\label{Theta}
\Theta_j=\Theta\quad\text{on}\quad [\beta_{j-1}, \beta_j]; \displaybreak[0]\\
\label{CPpositivity}
(-1)^{m-j}\Theta>0\quad\text{on}\quad (\beta_{j-1},\beta_{j});
\displaybreak[0]\\
\label{CPbound}
\Theta(\beta_j) = 0,\qquad
\Theta'(\beta_j) = - \frac{\cc}{n_j} \prod_{k\neq j}(\beta_{j}-\beta_k).
\end{gather}
\item Conversely, given constants $\beta_0<\beta_1<\cdots<\beta_m$, $\cc>0$
and a smooth function $\Theta$ on $[\beta_0,\beta_{m}]$ satisfying
~\eqref{CPpositivity}--\eqref{CPbound}, the K\"ahler metric given
by~\eqref{orthotoric} and~\eqref{Theta} defines an orthotoric structure on $\C
P^m_{a_0,\ldots a_m}$ and its toric orbifold quotients, such that the rational
Delzant polytope is the image of $[\beta_0,\beta_1]\times
[\beta_1,\beta_2]\times\cdots\times [\beta_{m-1},\beta_m]$ under the
elementary symmetric functions, with normals given by~\eqref{CPnormals}.

\item Any {\rm (nonsingular)} compact orthotoric K\"ahler $2m$-manifold $M$
arises in this way \textup(with $n_j=1$ for $j=0,\ldots m$\textup) and is
equivariantly biholomorphic to $\C P^m$.
\end{numlist}
\end{thm}

\subsection{Examples on weighted projective spaces}\label{s:orth-examples}

\subsubsection{The Fubini--Study metric}
We recall from~\cite[\S 5.4]{ACG1} that an orthotoric K\"ahler metric has
constant holomorphic sectional curvature $\cc$ if and only if $\oF_j=\oF_0$
for all $j=1,\ldots m$, where $\oF_0$ is a polynomial of degree $m+1$ with
distinct roots and leading coefficient $-\cc$. We then have
$\oF_0(t)=-\cc\prod_{j=0}^m (t-\beta_j)$ with $\beta_0<\cdots<\beta_m$, which
clearly satisfies~\eqref{CPpositivity}--\eqref{CPbound}.  Thus we see directly
that this orthotoric metric is defined on $\C P^m$, in accordance
with~\cite[\S 2.4]{ACG1}, where it was shown more generally that the
Fubini--Study metric on $\C P^m$ admits hamiltonian $2$-forms of arbitrary
order $\leq m$, in one to one correspondence with Killing potentials.

This form of the Fubini--Study metric is familiar for $m=1$,
when~\eqref{orthotoric} yields
\begin{equation*}
g = \frac{ d\xi^2 }{\cc(\xi-\beta_0)(\beta_1-\xi)} +
{\cc(\xi-\beta_0)(\beta_1-\xi)} dt^2.
\end{equation*}
Setting $2\xi=(\beta_1-\beta_0)z+\beta_0+\beta_1$,
$t=2\psi/\cc(\beta_1-\beta_0)$ and rescaling $g$ by $\cc$, we get
\begin{equation*}
g_{FS} = \frac{ dz^2 }{1-z^2} + (1-z^2) d\psi^2.
\end{equation*}
In arbitrary dimension $m$, the Fubini--Study metric is the `canonical' metric
associated to its simplex, hence is given here by~\eqref{toricmetric} with
\begin{align*}
{\bf G} &= \frac{1}{2} \Hess \biggl( \sum_{j=0}^{m} L_j(\boldsymbol\sigma)
\log |L_j(\boldsymbol\sigma)|\biggr)\\ \tag*{where} L_j(\boldsymbol\sigma) &=
\langle u_j,\boldsymbol\sigma \rangle + \lambda_j = \frac{2}{\cc}
\frac{\prod_{k=1}^m (\beta_j-\xi_k)}{\prod_{k\neq j} (\beta_j-\beta_k)}.
\end{align*}
It follows that
\begin{align*}\notag
\sum_{r,s=1}^{m} {\bf G}_{rs}d\sigma_r\,d\sigma_s &=
\frac 12 \sum_{j=0}^m L_j\,\Bigl( \frac{dL_j}{L_j}\Bigr)^2
=\frac 1\cc \sum_{j=0}^m
\frac{\prod_{k=1}^m (\beta_j-\xi_k)}{\prod_{k\neq j} (\beta_j-\beta_k)}
\biggl(\sum_{k=1}^m \frac{d\xi_k}{\xi_k-\beta_j}\biggr)^2\\
&= \sum_{p,q=1}^{m} \frac 1{\oF_0(\xi_p)} \sum_{j=0}^m
\biggl(\prod_{k\neq q} (\beta_j-\xi_k)\biggr)
\biggl(\prod_{k\neq j} \frac{\xi_q-\beta_k}{\beta_j-\beta_k}\biggr)
d\xi_p\,d\xi_q
\end{align*}
which immediately yields the orthotoric description~\eqref{orthotoric}, since
the inner sum over $j$ is $\Delta_q\,\delta_{pq}$ by the Lagrange
interpolation formula.

\subsubsection{Bochner-flat metrics}
More generally, any {\it extremal} orthotoric metric~\eqref{orthotoric} for
which $\Theta_j = \Theta$ is necessarily Bochner-flat~\cite[\S 5.4]{ACG1}; in
this case $\Theta$ must be a polynomial of degree $\leq m+2$, and the boundary
conditions~\eqref{bound} imply that $\Theta$ has $m+1$ or $m+2$ distinct
roots. The former case gives the Fubini--Study metric and its orbifold
quotients, while the latter recovers the Bochner-flat examples of
\cite{bryant}, which are defined on $\C P^{m}_{a_0, \ldots a_{m}}$ for
distinct weights $a_0,\ldots a_m$. Indeed, for any positive integers $a_0 >
\cdots >a_m$ we take the metric \eqref{orthotoric} with
\begin{equation*}
\Theta_j(t)=\Theta(t) = - (t-\beta) \Theta_0(t)= c (t - \beta)
\prod_{j=0}^m (t-\beta_j),
\end{equation*}
where $c>0$ is a homothety factor for the metric and we deduce
from~\eqref{CPbound} that the real numbers $\beta_0 < \cdots < \beta_m <
\beta$ and $c>0$ satisfy
\begin{equation} \label{betas}
\beta_j = \beta-\frac{a_j}{\prod_{k=0}^m a_k}.
\end{equation}
This metric is Bochner-flat (see \cite[Proposition 16]{ACG1}) and compactifies
on $\C P^{m}_{a_0, \ldots a_{m}}$ (see Theorem~\ref{CPorthotoric}). As shown
by Bryant~\cite{bryant}, there are actually Bochner-flat metrics on $\C
P^{m}_{a_0, \ldots a_{m}}$ for \emph{any} choice of weights. There is an
easy way to see this~\cite{LP} using the relation between
Bochner-flat metrics and flat CR structures found by Webster~\cite{webster}.
Indeed $\C P^m_{a_0,\ldots a_m}$ is a quotient of $S^{2m+1}$ by a weighted
$S^1$-action by CR automorphisms of the flat CR structure, and the Sasakian
structure induced by the associated Reeb field gives rise to a Bochner-flat
K\"ahler metric on the quotient.

The Bochner-flat metrics on weighted projective spaces are all toric
(see~\cite{Abreu1} for the general form in momentum coordinates).  However,
when the weights are not distinct, they are not orthotoric (apart from the
Fubini--Study metric): from our point of view, the Bochner-flat metric is
endowed with a natural hamiltonian $2$-form which is (an affine deformation
of) the normalized Ricci form~\cite{ACG1} and it has order $m$ if and only if
the weights $a_j$ are distinct.

\begin{rem} Note that the orthotoric Bochner-flat K\"ahler metric on
a weighted projective space is unique (up to isomorphism and scale):
$\beta_0,\ldots \beta_m$ are determined as above (the choice of $\beta$ can be
absorbed in the coordinate freedom). In fact a stronger uniqueness result is
true: the Bochner-flat metric is the unique \emph{extremal} K\"ahler metric
(up to isomorphism and scale) on \emph{any} weighted projective space. To see
this, recall that the second deRham cohomology group of $\C P^m_{a_0,\ldots
a_m}$ is one dimensional, so there is only one K\"ahler class up to scale
(this follows, for instance, by the Smith--Gysin sequence for the space of
orbits, $\C P^m_{a_0,\ldots a_m}$, of the weighted $S^1$-action on the
$(2m+1)$-sphere); therefore the uniqueness result of Guan~\cite{guan} (which
readily generalizes to orbifolds) applies to the K\"ahler class of $\C
P^m_{a_0,\ldots a_m}$.

The uniqueness implies that any toric $2m$-orbifold of constant scalar
curvature, whose rational Delzant polytope is a simplex, is an orbifold
quotent of $\C P^{m}$. Note that the Futaki invariant of $\C P^m_{a_0,\ldots a_m}$
vanishes if and only if $a_0=a_1=\cdots =a_m$.
\end{rem}

\subsection{K\"ahler--Einstein orthotoric surfaces}\label{s:KEorbifolds}

In this subsection we present new examples of K\"ahler--Einstein metrics on
compact orbifolds. As we have seen in the previous subsection, we have to work
beyond the context of weighted projective spaces, so we consider polytopes
with more than $m+1$ codimension one faces.  We restrict attention to complex
orbifold surfaces ($m=2$) in order to make the construction completely
explicit. In this case, a polytope with more than $m+1=3$ faces necessarily
has $2m=4$ faces and we are in the `generic' case where the roots
$\xi_1,\xi_2$ are {\it everywhere} distinct on $\Delta$.

According to \cite[\S5.3]{ACG1}, an orthotoric K\"ahler metric on a
$4$-orbifold is K\"ahler--Einstein if and only if $\Theta_j(t) = -P_j(t)/C$,
$j=1,2$ for some positive constant $C$, and some $\pm$-monic polynomials $P_j$
of degree $3$, such that $P_1(t) - P_2(t) = c$ where $c$ is a constant.  The
Bochner tensor vanishes precisely when $c=0$, and the metric is then the
Fubini--Study metric.  We therefore assume that $c\neq 0$ in order to obtain
new examples. Also, for compactness, the scalar curvature must be positive
(otherwise the Ricci tensor would be nonpositive, contradicting the existence
of Killing vector fields with zeros), which implies that the polynomials $P_j$
are monic.

It remains to solve the compactification conditions. For simplicity, we shall
take the lattice $\Lambda\subset\R^2$ to be $\Z^2$ or a sublattice.  The
conditions of Proposition~\ref{orthocompact} can then be satisfied by
supposing that $P_j$ has integer roots (including the endpoints $\alpha_j$ and
$\beta_j$) and $C$ is chosen so that $2/\Theta'(\alpha_j) = c_j^{\alpha}= -2C
/ P_j'(\alpha_j)$ and $2/\Theta'(\beta_j) = c_j^{\beta}=-2C / P_j'(\beta_j)$
are all integers for $j=1,2$.

The condition \eqref{bound} implies that $P_1$ and $P_2$ have three distinct
roots, $p_1<q_1<r_1$ and $p_2>q_2>r_2$, respectively.  The condition $P_1-P_2
= c$ reads
\begin{align*}
p_1 + q_1 + r_1 &= p_2 + q_2 + r_2,\\
p_1^2 + q_1^2 + r_1^2 &= p_2^2 + q_2^2 + r_2^2.
\end{align*}
Positivity and~\eqref{bound} give $\alpha_1=p_1$, $\beta_1=q_1$,
$\alpha_2=q_2$, $\beta_2=p_2$ and hence (without loss) $q_1<q_2$.  Taking the
roots to be all integral, we note that, up to an affine deformation of the
hamiltonian $2$-form and orbifold coverings/quotients, we can also assume that
$\gcd(p_1,q_1,r_1,p_2,q_2,r_2)=1$ and
\begin{equation*}
p_1 + q_1 + r_1 = 0 = p_2 + q_2 + r_2.
\end{equation*}
A class of solutions to this problem is obtained by taking any coprime
positive integers $(p,q)$ with $p>q$ and putting
\begin{equation*}
p_1=-p, \quad q_1 = -q, \quad r_1 =p+q, \quad
p_2= p, \quad q_2 = q, \quad r_2 = -p-q.
\end{equation*}
With these assumptions, we have
\begin{align}  \nonumber
\alpha_1 &= -p, \quad \beta_1 = -q, \quad \alpha_2  = q, \quad \beta_2 =
p;\\ \label{solution1}
\Theta_1(\xi) &= -\frac{(\xi+p)(\xi+q)(\xi-p-q)}{C};\\ \label{solution2}
\Theta_2(\xi) &= -\frac{(\xi-p)(\xi-q)(\xi+p+q)}{C}.
\end{align}
The corresponding Delzant polytope $\Delta$ is the quadrilateral with vertices
\begin{equation*}
(0,-p^2), \quad (0, - q^2), \quad (p-q, -pq), \quad (q-p, -pq)
\end{equation*}
and one-dimensional faces $\Fa_j^{\alpha}$, $\Fa_j^{\smash\beta}$, $j=1,2$
determined by the lines
\begin{equation*}
\{ \boldsymbol\sigma : \ell_j^{\alpha}(\boldsymbol\sigma)= 0 \}, \quad
\{ \boldsymbol\sigma : \ell_j^{\smash\beta}(\boldsymbol\sigma)= 0 \},
\end{equation*}
where
\begin{align*}
\ell_1^{\alpha} (\boldsymbol\sigma) &= p^2  + p\sigma_1 + \sigma_2,&
\ell_1^{\smash\beta}(\boldsymbol\sigma) &= q^2 + q \sigma_1 + \sigma_2,\\
\ell_2^{\alpha} (\boldsymbol\sigma) &= q^2 - q \sigma_1 + \sigma_2,&
\ell_2^{\smash\beta}(\boldsymbol\sigma) &= p^2 - p\sigma_1 + \sigma_2.
\end{align*}
Furthermore, letting $2C = (p-q)(2q+p)(2p+q)$ in \eqref{solution1} and
\eqref{solution2}, we get
\begin{equation*}\notag
\Theta_1'(\alpha_1) = \Theta_2'(\beta_2) = 2/(2p+q), \quad
\Theta_1'(\beta_1) = \Theta_2'(\alpha_2) = -2/(2q+p),
\end{equation*}
so the conditions of Proposition~\ref{orthocompact} are satisfied with
\begin{equation*}
c_1^{\alpha}= c_2^{\smash\beta}= 2q+p, \quad
c_2^{\alpha} = c_1^{\smash\beta}= 2p+q.
\end{equation*}
Thus, according to Proposition~\ref{orthocompact}, the corresponding
K\"ahler--Einstein orthotoric metric $g_{p,q}$ compactifies on the toric
orbifold \ka surface $M(p,q)$ classified by
\begin{equation*}
(\Delta, \Lambda, c_1^{\alpha},c_1^{\smash\beta},
c_2^{\alpha},c_2^{\smash\beta}),
\end{equation*}
where $\Lambda$ is the standard lattice $\Z^2 \subset \R^2$ (in which case
$c_j^{\alpha}, c_j^{\smash\beta}$ are nothing but the integer labels
corresponding to the 1-dimensional faces of $\Delta$, see \S\ref{s:toric}).

We claim that two orbifold surfaces $M(p,q)$ and $M(p',q')$ are
biholomorphically equivalent iff $p=p'$ and $q=q'$. Indeed, in order to be
biholomorphic as complex orbifolds, $M(p,q)$ and $M(p',q')$ must be isomorphic
as toric varieties. Therefore, the cooresponding polytopes $\Delta$
and $\Delta'$ must determine congruent fans~\cite[Thm.9.4]{LT}. One easily
checks that the latter happens iff $(p,q) = (p',q')$; alternatively, using the
uniqueness of the hamiltonian $2$-form established in Propostion \ref{prop1}
below, one can see that the K\"ahler--Einstein metrics $g_{p,q}$ and
$g_{p',q'}$ are locally isometric if and only if $(p,q) = (p',q')$.

We summarize our construction as follows.
\begin{thm} There is a family of nonequivalent compact K\"ahler--Einstein
orthotoric orbifold surfaces $(M(p,q),g_{p,q})$, depending on coprime
positive integers $q<p$.
\end{thm}

\begin{rem} (i) According to the results of \cite{ACG0}, the primitive part
of the hamiltonian $2$-form $\phi$ associated to $g_{p,q}$ defines an
integrable almost-complex structure $I$ on $M(p,q)$, which is compatible with
$g_{p,q}$ but induces the opposite orientation to the one of $M(p,q)$.  With
respect to this structure, $(M(p,q),g_{p,q},I)$ become a compact, Einstein,
non-\ka hermitian complex orbifold surface (see \cite{lebrun0} for a
classification in the smooth case).

\smallbreak \noindent (ii) A similar construction yields a countable family of
compact orbifold complex surfaces supporting orthotoric weakly Bochner-flat
metrics which are neither Bochner-flat nor K\"ahler--Einstein (see \cite{ACG0}
for a classification in the smooth case).

\smallbreak \noindent (iii) According to \cite{BG1}, any K\"ahler--Einstein
orbifold $(M,g,J,\omega)$ of complex dimension $m$ gives rise to a
Sasaki--Einstein structure on the total space $S$ of a principal $S^1$
$V$-bundle over $M$ (which is suitably associated to the canonical bundle of
$M$). In general, $S$ is an $(2m + 1)$-dimensional orbifold rather than a
manifold, but it may happen that $S$ is nonsingular even though $M$ is
singular~\cite{BG1}: in fact, $S$ is nonsingular if and only if all local
uniformizing groups of $M$ inject into the structure group $S^1$ (see
\cite[Theorem 2.3]{BG1}). In the case of {\it toric} K\"ahler orbifolds all
local uniformizing groups are abelian \cite{LT} so that if $S$ is nonsingular,
then all local uniformizing groups of $M$ must be cyclic. Using this
observation one can show that the universal orbifold covers $\widehat M(p,q)$
of $M(p,q)$ (the one which corresponds to the lattice generated by the normals
of $\Delta$, see Remark ~\ref{cover}) give rise only to {\it singular}
$5$-dimensional Sasaki--Einstein orbifolds.
\end{rem}

\section{Compact K\"ahler manifolds with hamiltonian $2$-forms}
\label{s:global}

We now combine the work of the previous three sections to classify, up to a
covering, compact K\"ahler manifolds with a hamiltonian $2$-form. In the case
of K\"ahler surfaces we can refine the classification. We end by giving some
examples of extremal and weakly Bochner-flat K\"ahler metrics with hamiltonian
$2$-forms.

\subsection{General classification}\label{s:genclass}

There are two parts to a general classification of compact K\"ahler manifolds
with a hamiltonian $2$-form. First, we must classify the possible equivariant
biholomorphism types of manifolds which can admit a hamiltonian $2$-form.
Second, we describe the compatible K\"ahler structures on such manifolds which
\emph{do} admit hamiltonian $2$-forms.

The equivariant biholomorphism type is described in parts (ii)--(iv) of the
following theorem: we show that, up to a blow-up and a covering, a compact
K\"ahler manifold with a hamiltonian $2$-form of order $\ell$ is biholomorphic
to a projective bundle of the form $P(\cL_0\oplus \cL_1\oplus\cdots \oplus
\cL_\ell)\to S$ where $\cL_j$ are holomorphic line bundles over a product $S$
of K\"ahler manifolds $S_a$. Such a bundle admits an action of a complex
$\ell$-torus $\T^c$, defined by scalar multiplication in each line bundle (an
$(\ell+1)$-torus action on $\cL_0\oplus\cdots \oplus \cL_\ell$) modulo overall
scalar multiplication \textup(which acts trivially on the
projectivization\textup). In part (v) of the theorem, we show that the
relevant K\"ahler structures are given by a special case of the generalized
Calabi construction with $\cV = \C P^{\ell}$. Conversely we show that this
construction produces compact K\"ahler manifolds with a hamiltonian $2$-form.

\begin{thm} \label{thm:main}
Let $(M,g,J,\omega)$ be a compact connected K\"ahler $2m$-manifold with a
hamiltonian $2$-form $\phi$ of order $\ell\geq 0$, with nonconstant roots
$\xi_1,\ldots \xi_\ell$ and \textup(distinct\textup) constant roots
$\eta_1,\ldots \eta_N$, $N\geq 0$.

\begin{numlist}
\item The elementary symmetric functions $(\sigma_1,\ldots\sigma_\ell)$ of
$(\xi_1,\ldots\xi_\ell)$ are the components of the momentum map
$\boldsymbol\sigma\colon M\to\R^{\ell*}$ of an $\ell$-torus $\T\leq
\mathrm{Isom}(M,g)$. The image $\Delta$ of $\boldsymbol\sigma$ is a Delzant
simplex in $\R^{\ell*}$, whose interior is the image under the elementary
symmetric functions of a domain $D=\prod_{j=1}^{\smash\ell}
(\beta_{j-1},\beta_j)$ with $\beta_0<\beta_1<\cdots<\beta_\ell$.

\item Let $S_\Delta$ be the stable quotient of $M$ by the complex torus $\T^c$
and let $\smash{\hat M}$ be the blow-up of $M$ along the inverse image of the
codimension one faces $\Fa_0,\Fa_1,\ldots \Fa_\ell$ of $\Delta$.  Then there
are holomorphic line bundles $\cL_0,\cL_1,\ldots \cL_\ell$ over $S_\Delta$
\textup(uniquely determined up to overall tensor product with a holomorphic
line bundle\textup) such that $\hat M$ is $\T^c$-equivariantly biholomorphic
to $P(\cL_0\oplus \cL_1\oplus\cdots \oplus \cL_\ell)\to S_\Delta$.

\item $S_\Delta$ is covered by a product $S$ of $N$ Hodge K\"ahler manifolds
$(S_a,\pm g_a,\pm \omega_a)$ of dimension $2d_a$, indexed by the constant
roots $\eta_a$ \textup($S$ is a point if $N=0$\textup). There are constants
$\cc,C_1,\ldots C_N$ such that for $j=0,\ldots\ell$, $a=1,\ldots N$
\begin{equation} \label{chernclasses}
\tfrac12 \cc \biggl(\prod_{k\neq j} (\eta_a-\beta_k) \biggr)
\bigl(C_a(\eta_a-\beta_j)+1\bigr)[\omega_a/2\pi]
\end{equation}
is an integral cohomology class on $S_a$, and the pullback of $\cL_j$ to $S$
is a tensor product $\bigotimes_{a=1}^N \pi_a^* \cL_{j,a}$, where $\pi_a$ is
the projection of $S$ to $S_a$ and $\cL_{j,a}\to S_a$ is a holomorphic line
bundle with first Chern class given by~\eqref{chernclasses}.

\item The subset $\cB$ of those $j\in \{0,\ldots\ell\}$ for which the blow up
over the face $\Fa_j$ is nontrivial corresponds bijectively to a subset $\cC$
of $\{1,\ldots N\}$ such that for $j\in\cB$ corresponding to $a\in\cC$,
$\eta_a=\beta_j$, $S_a = \C P^{d_a}$, $\pm g_a$ is the Fubini--Study metric on
$S_a$ of constant holomorphic sectional curvature $\pm \cc \prod_{k\neq
j}(\beta_{j}-\beta_k)$, and \textup(without loss\textup) $\cL_{j,a}=\cO(-1)$
and $\cL_{k,a}=\cO$ for $k\neq j$.

For $a\notin\cC$ either $\eta_a<\beta_0$ or $\eta_a>\beta_\ell$.

\item The K\"ahler metric on $M$ and its pullback to $\hat M$ are determined
by the explicit metric~\eqref{metric} on $M^0$, where\textup:
\begin{itemize}
\item the pullback to $S=\prod_{a=1}^N S_a$ of the K\"ahler quotient metric on
$S_\Delta$ induced by $\boldsymbol\sigma(\xi_1,\ldots\xi_\ell) \in\Delta^0$ is
the K\"ahler product metric
\begin{equation}\label{kpmetric}
\sum_{a=1}^N \biggl(\prod_{j=1}^\ell (\eta_a-\xi_j)\biggr)g_a;
\end{equation}
\item $\theta_1,\ldots\theta_\ell$ are the components of a connection on $\hat
M\to S_\Delta$ associated to a principal $\T$-connection\textup;
\item for $j=1,\ldots \ell$, $F_j(t)=\Mpc(t)\Theta(t)$, where
$\Mpc(t)=\prod_{a=1}^N (t-\eta_a)^{d_a}$,
\begin{gather}\label{Theta-pos}
(-1)^{m-j}\Theta>0\quad\text{on}\quad (\beta_{j-1},\beta_{j}),\\
\Theta(\beta_j) = 0,\qquad
\Theta'(\beta_j) = - \cc \prod_{k\neq j}(\beta_{j}-\beta_k),
\label{Theta-bound}
\end{gather}
and the metric on the $\C P^\ell$-fibres of $\hat M\to S_\Delta$ is the
orthotoric K\"ahler metric~\eqref{orthotoric} with
$\oF_j(t)=\Theta(t)$\textup;
\end{itemize}
\end{numlist}

Conversely, suppose $S$ is a product of Hodge K\"ahler manifolds $(S_a,\pm
g_a, \pm\omega_a)$ and constants $\beta_0,\ldots\beta_\ell$, $\eta_1,\ldots
\eta_N$, $c,C_1,\ldots C_N$ satisfying the conditions in \textup{(i)--(iv)}
above and such that~\eqref{kpmetric} is positive for
$\boldsymbol\sigma(\xi_1,\ldots\xi_\ell) \in\Delta^0$.

Then there is a complex manifold $M$ obtained by a blow-down of a projective
bundle ${\hat M}=P({\cL_0} \oplus \cL_1\oplus \cdots \oplus {\cL}_{\ell})\to
S$ which gives rise to these data. Further, for any smooth function $\Theta$
on $[\beta_0,\beta_{\ell}]$ satisfying~\eqref{Theta-pos}--\eqref{Theta-bound},
a K\"ahler metric of the form \eqref{metric}, with $F_j(t) = \Mpc(t)
\Theta(t)$, is globally defined on $M$ and admits a hamiltonian $2$-form of
order $\ell$.
\end{thm}

\begin{proof} Consider the explicit form~\eqref{metric} of the metric on the
open subset $M^0$ of $M$. By Lemma~\ref{tricky}, the map
$\boldsymbol\sigma=(\sigma_1,\ldots \sigma_\ell)\colon M\to \R^\ell$ generates
a rigid hamiltonian torus action, and the K\"ahler quotient metric
(i.e.,~\eqref{kpmetric}) is clearly semisimple, so that by
Theorem~\ref{generalizedcalabi}, there is a cover of $M$ which is given by the
generalized Calabi construction. The covering is straightforward: there is a
discrete group $\Gamma$ of holomorphic isometries of $S$ which lifts to the
bundle $M_0\times_{\T^c} \cV$ and $\hat M$ is the quotient. We shall therefore
suppose $\Gamma$ is trivial in the following.

The K\"ahler metrics $\pm\omega_a$ are determined by~\eqref{kpmetric} and the
constants $c_{a0}$ and ${\boldsymbol c}_a = (c_{a1}, \ldots c_{a\ell})$
appearing in the generalized Calabi data are
\begin{equation}\label{constants}
c_{a0} = \eta_{a}^{\ell}, \quad c_{ar}= (-1)^r\eta_a^{\ell-r}, \quad
r=1,\ldots\ell.
\end{equation}
Since the roots $\xi_1,\ldots\xi_\ell$ of the momentum polynomial $\Mpn(t)$
are smooth, functionally independent and pairwise distinct on $M^0$, with
orthogonal gradients, the toric K\"ahler manifold $\cV$ appearing in the
generalized Calabi data is orthotoric.

\begin{numlist}
\item $\boldsymbol\sigma$ is a momentum map by definition, and by
Proposition~\ref{orthopoly}, its image $\Delta$ is a simplex as stated.  In
particular, the codimension one faces $\Fa_0,\Fa_1,\ldots \Fa_\ell$ correspond
to the boundary points $\beta_0<\beta_1<\cdots <\beta_\ell$ of $D$,
and $\cV$ is biholomorphic to $\C P^\ell$.

\item $M^0$ is a holomorphic principal $\T^c$-bundle over $S$ and the blow up
of $\hat M$ along the inverse image of the codimension one faces is
equivariantly biholomorphic to a projective bundle $M^0\times_{\T^c} \C
P^\ell$ with a global fibre preserving $\T^c$ action.  This action identifies
$\hat M$ with $P(\cL_0\oplus\cL_1\oplus\cdots \oplus \cL_\ell)$ for
holomorphic line bundles $\cL_0,\cL_1\ldots \cL_\ell$ over $S$ (uniquely
determined as stated) in such a way that the $2(\ell-1)$-dimensional orbits of
$\T^c$ in each fibre are orbits of elements of $P(\cL_0\oplus\cdots \oplus
\cL_\ell)$ with one homogeneous coordinate vanishing.  We label the line
bundles so that the codimension one face $\Fa_j$ corresponds to the orbit of
$\T^c$ with $\cL_j$ component vanishing.

\item We only need to construct the line bundles $\cL_{j,a}$ and establish the
formula~\eqref{chernclasses} for their first Chern classes. The explicit
form~\eqref{metric} of the metric on the principal bundle $M^0$ shows that the
connection $1$-forms $\theta_r$, with $\theta_r(K_s)=\delta_{rs}$, satisfy
\begin{equation}
d\theta_r = \sum_{a=1}^N (-1)^r \eta_a^{\ell-r} \omega_a,
\end{equation}
where $\pm\omega_a$ are the K\"ahler forms of the globally defined metrics
$\pm g_a$ on $S_a$. Note that the $\theta_r$ are not necessarily integral.
The integral principal connection forms are those which evaluate to integers
on the Euler fields $X_0,\ldots X_\ell$, which, according to
Proposition~\ref{orbi-simplex}, are given by~\eqref{duallat}:
\begin{equation*}\notag
\theta_{p,q} =  \sum_{r=0}^\ell \tfrac 1 2 \cc\bigl(\sigma_{r}^\beta(\hat
\beta_q)-\sigma_{r}^\beta(\hat \beta_p)\bigr) \theta_r,
\end{equation*}
More specifically, this is the connection form of the line bundle
$\cL_p^{-1}\otimes\cL_q$ (up to a sign convention). The curvature form of
$\cL_p^{-1}\otimes\cL_q$ is therefore
\begin{equation}\label{curvdiff}
\begin{split}
d\theta_{p,q}&= \sum_{a=1}^N
\sum_{r=0}^\ell \tfrac 1 2 \cc(-1)^r\bigl(\sigma_{r}^\beta(\hat\beta_q)
-\sigma_{r}^\beta(\hat \beta_p)\bigr)\eta_a^{\ell-r}\omega_a\\
&=\sum_{a=1}^N \tfrac 1 2 \cc\biggl(\prod_{k\neq q} (\eta_a-\beta_k)
- \prod_{k\neq p} (\eta_a-\beta_k) \biggr)\omega_a.
\end{split}\end{equation}
It follows that for each $a=1,\ldots N$, the corresponding $2$-form in this
sum is \emph{integral} in the sense that the cohomology class
\begin{equation}\label{chernclass}
\tfrac 1 2 \cc\biggl(\prod_{k\neq q} (\eta_a-\beta_k)
- \prod_{k\neq p} (\eta_a-\beta_k) \biggr) [\omega_a/2\pi]
\end{equation}
is in the image of $H^2(S_a,\Z)$ in $H^2(S_a,\R)$.  If $\eta_a=\beta_j$ for
some $j$, we deduce (by taking $p=j$, $q\neq j$) that $\frac 1 2 \cc
\bigl(\prod_{k\neq j} (\eta_a-\beta_k)\bigr)\omega_a$ is integral. Otherwise,
this will differ from an integral class by a constant.  Hence there are
constants $C_1,C_2,\ldots C_N$ such that for each $j=0,\ldots\ell$ and
$a=1,\ldots N$, the $2$-form
\begin{equation}\label{curvform}
\tfrac 1 2 \cc \biggl(\prod_{k\neq j} (\eta_a-\beta_k) \biggr)
\bigl(C_a(\eta_a-\beta_j)+1\bigr) \omega_a
\end{equation}
is also integral. Now the Lefschetz Theorem for $(1,1)$-classes implies that
there are holomorphic line bundles $\cL_{j,a}$ with connection over $S_a$
whose curvature forms are given by~\eqref{curvform}: the first Chern classes
are then as stated in~\eqref{chernclasses}. It follows that $\cL_j$ is the
tensor product of $\bigotimes_{a=1}^N \pi_a^* \cL_{j,a}$ by a flat line bundle
$\cF_j$. Since any flat line bundle on $S$ is a tensor product of flat line
bundles pulled back from the factors $S_a$, we may use the freedom in the
choice of $\cL_{j,a}$ to make $\cF_j$ trivial.

Finally, note that for each $a$, the Chern classes $c_1(\cL_{j,a})$ cannot
vanish for all $j$. It follows that the manifold $S_a$ is Hodge, i.e., admits
a K\"ahler metric whose K\"ahler class is integral in cohomology.

\item For any $\boldsymbol\sigma$ in $\Delta^0$, the K\"ahler quotient
metric~\eqref{kpmetric} is global on $S$, so that $\Mpn(\eta_a)=
\prod_{j=1}^{\ell} (\xi_j-\eta_a)$ does not vanish on $\Delta^0$. Hence no
$\eta_a$ can belong to any of the open intervals
$(\beta_{j-1},\beta_j)$. Clearly when $\eta_a=\beta_j$ for some
$j=0,\ldots\ell$, $\Mpn(\eta_a)$ vanishes on the codimension one face $\Fa_j$
of $\Delta$ and this is precisely the condition that a blow-down occurs (over
the factor $S_a$).  The rest is immediate from the definition of generalized
Calabi data apart from the normalization of the Fubini--Study metric on
$S_a$. For this we note that the formula~\eqref{chernclasses} gives
\begin{equation}\label{FSnormalization}
c_1(\cL_{j,a}) = \tfrac 12 \cc\biggl(\prod_{k\neq j}(\beta_j -
\beta_k)\biggr)[\omega_a/2\pi]
\end{equation}
(and $c_1(\cL_{k,a})=0$ for $k\neq j$).  Since $\cL_{j,a}$ has to be
$\cO(-1)$, we must have
\begin{equation*}
[\rho_a/2\pi] = (d_a+1) c_1(\cL_{j,a}) = \tfrac12 \cc (d_a+1)
\biggl(\prod_{k\neq j} (\beta_j-\beta_k)\biggr)[\omega_a/2\pi],
\end{equation*}
where $\rho_a = \frac{\Scal_a}{2d_a} \omega_a$ is the Ricci form of the
Fubini--Study metric $\pm g_a$. The holomorphic sectional curvature $\pm\frac1
{d_a(d_a+1)}\Scal_a$ is therefore as stated. (Note that this can be always
achieved by rescaling $\omega_a$.)

\item This is immediate from the explicit form of the metric, the generalized
Calabi construction, and the necessity of the conditions of
Theorem~\ref{CPorthotoric} for the compactification of orthotoric K\"ahler
metrics on $\C P^\ell$.
\end{numlist}

For the converse, observe first that the integrality conditions ensure (by the
Lefschetz Theorem for (1,1)-classes) that there are holomorphic line bundles
$\cL_{j,a}$ over $S_a$ with first Chern classes given by \eqref{chernclasses},
equipped with compatible connections whose curvatures are given
by~\eqref{curvform}, and we define $\cL_j=\prod_{a=1}^N \pi_a^*\cL_{j,a}$.

Let $\hat M=P(\cL_0\oplus\cL_1\oplus\cdots\oplus\cL_\ell)\cong
M^0\times_{\T^c}\C P^\ell$, where $\T^c$ acts by scalar multiplication on each
line bundle $\cL_j$ modulo overall scalar multiplication on the direct sum
(which acts trivially on the projective bundle), $M^0$ is the union of the
open $\T^c$ orbits in each fibre of $\hat M\to S$, and $\C P^\ell$ is toric
under $\T^c$.

By \S\ref{s:restricttoric} (see also Theorem~\ref{generalizedcalabi}) $\hat M$
has a blow-down $M$, which collapses a family of divisors (which are closures
of complex codimension one $\T^c$-orbits) corresponding to $\eta_a\in\cC$
along the $\C P^{d_a}$ fibrations induced by the connection on $\hat M$.

Because of the sufficiency of the conditions of Theorem~\ref{CPorthotoric} for
the compactification of orthotoric K\"ahler metrics on $\C P^\ell$, we have
generalized Calabi data for the construction of a K\"ahler metric on $M$ using
Theorem~\ref{generalizedcalabi}, where $\cV=\C P^\ell$ equipped with this
orthotoric structure, the connection has curvature~\eqref{curvform}, and the
constants are given by~\eqref{constants}. On $M^0$ the K\"ahler structure is
given by~\eqref{metric}.

The hamiltonian $2$-form $\phi = \sum_{r=1}^\ell (\sigma_r d\sigma_1
-d\sigma_{r+1})\wedge dt_r$ (defined on $M^0$) also extends on $M$. Indeed, it
follows from \cite[\S 2.2]{ACG1} that the $2$-jet of $\phi$ is a parallel
section (over $M^0$) of a vector bundle with linear connection globally
defined on $M$. Since $M\smallsetminus M^0$ has codimension at least two in
$M$, $\phi$ extends to the whole of $M$.
\end{proof}

\begin{rem} It follows from the proof of Theorem~\ref{generalizedcalabi} that
$M$ is a bundle of restricted toric K\"ahler manifolds over
$\prod_{a\notin\cC} S_a$. The typical fibre $\cX$ is a toric K\"ahler manifold
of dimension $2k$, $k=\ell+\sum_{a\in\cC} d_a$, obtained as a blow-down of a
$\C P^\ell$ bundle over a product of $\#\cC<\ell+1$ projective spaces as in
\S\ref{s:restricttoric}, and admits a hamiltonian $2$-form of order $\ell$.
However, by~\cite[\S 2.4,\S 5.4]{ACG1}, the Fubini--Study metric on $\C P^k$
admits a hamiltonian $2$-form of order $\ell$ with any number of distinct
constant roots between $0$ and $\ell+1$, with all factors in the K\"ahler
quotient being blown down over some face of the Delzant polytope of $\C
P^\ell$. It follows that $\cX$ is biholomorphic to $\C P^k$, and $M$ is a
bundle of projective spaces over $\prod_{a\notin\cC} S_a$. However, the metric
on the fibres need not be orthotoric unless $k=\ell$. In fact it is not hard
to see directly that the blow-down of
$P(\cL_0\otimes\cO\oplus\cL_1\otimes\cO\oplus\cdots\oplus
\cL_j\otimes\cO(-1)\oplus\cdots \oplus\cL_\ell\otimes\cO)\to S'\times\C P^d$ is
biholomorphic to $P(\cL_0\oplus\cL_1\oplus\cdots
\cL_j\otimes\C^{d+1}\oplus\cdots \oplus \cL_\ell)\to S'$.
\end{rem}

\subsection{K\"ahler surfaces with hamiltonian $2$-forms}

In this subsection we specialize to the case that $(M,g,J,\omega)$ is a smooth
compact \ka surface with a {nontrivial} hamiltonian $2$-form; here nontrivial
means that $\phi$ is not a constant multiple of $\omega$.  We obtain a
complete classification, overcoming the issue of coverings raised in the
previous subsection. We first recall that if $\phi$ is a nontrivial
hamiltonian $2$-form, then for any real numbers $a,b$ ($a \neq 0$), the affine
deformation $a\phi + b\omega$ is again a nontrivial hamiltonian $2$-form (of
the same order as $\phi$).

\begin{prop}\label{prop1} Let $(M,g,J)$ be a connected \ka surface not
of constant holomorphic sectional curvature. Then $(M,g,J)$ admits at most one
\textup(up to an affine deformation\textup) nontrivial hamiltonian $2$-form,
even locally.
\end{prop}
\begin{proof} According to~\cite[Lemmas 2 and 6]{ACG0}, the primitive
part $\phi_0$ of a nontrivial hamiltonian $2$-form $\phi$ defines (on the open
dense subset $U$ where $\phi_0 \neq 0$) a conformally \ka hermitian structure
$(g,I)$, such that $I$ and $J$ induce opposite orientations on $U$; then the
antiselfdual tensor $W^-$ of $g$, with respect to orientation induced by $J$,
has degenerate spectrum on $U$, hence on $M$; moreover, $\phi_0$ is an
eigenform of $W^-$, whose eigenvalue, at each point where $W^-$ does not
vanish, is the (unique) simple eigenvalue of $W^-$.  We also know that $\phi$
commutes with the Ricci form $\rho$---see \cite[\S 2.2]{ACG1}---so on an open
subset where $\rho_0 \neq 0$, $\phi_0$ is proportional to $\rho_0$.

It follows that on any open subset where $W^- \neq 0$ or $\rho_0\neq 0$, the
primitive parts of two nontrivial hamiltonian $2$-forms, $\phi$ and $\phi'$,
are related by $\phi_0 = f \phi_0'$ for a smooth function $f$; since the
primitive part of any hamiltonian $2$-form satisfies $d(\phi_0/|\phi_0|^3) =
0$ (see \cite[Lemma 2]{ACG0}), $f$ must be a constant, i.e., $\phi_0 = a
\phi_0'$. By unique continuation \cite[\S 2.2]{ACG1}, this equality holds
everywhere on $M$ and so $\phi - a\phi'$ is a hamiltonian $2$-form with
vanishing primitive part, hence a multiple of $\omega$. Thus $\phi$ and
$\phi'$ are affinely equivalent unless $W^-$ and $\rho_0$ are identically
zero.
\end{proof}
\begin{rem} The above result is optimal: according to \cite[\S 2.3]{ACG1},
each of the manifolds $\C P^2$, $\C^2$ and $\C\cH^2$ endowed with its
canonical K\"ahler structure admits a $9$-dimensional family of nontrivial
hamiltonian $2$-forms.
\end{rem}

\begin{thm}\label{surfaces} Let $(M,J)$ be a compact complex surface which
supports a \ka metric $g$ with a nontrivial hamiltonian $2$-form $\phi$.
Then the following cases occur.
\begin{numlist}
\item $\phi$ is of order zero\textup; then $(M,J)$ is biholomorphic to a
compact locally symmetric \ka surface of reducible type.

\item $\phi$ is of order one\textup; then $(M,J)$ is biholomorphic to either
${\C P}^2$ or to a ruled surface of the form $P({\mathcal O}\oplus {\mathcal
L}) \to S$ where $S$ is a compact complex curve
and ${\mathcal L}$ is a holomorphic line bundle over $S$ of
positive degree.

\item $\phi$ is of order two\textup; then $(M,J)$ is biholomorphic to ${\C
P}^2$.
\end{numlist}
Each complex surface listed in {\rm (i)--(iii)} above admits {\rm (}infinitely
many{\rm )} K\"ahler metrics with nontrivial hamiltonian $2$-forms of the
corresponding order.
\end{thm}
\begin{numlproof}
\item If the order of $\phi$ is zero, i.e., if $\phi$ is parallel, then, by
the deRham decomposition theorem, the universal cover $(\tilde M,\tilde g)$
of $(M,g)$ is a K\"ahler product $({\mathbb U}_1\times{\mathbb U}_2,
g_1\times g_2)$ where each ${\mathbb U}_i$ biholomorphic to $\C P^1$,
$\C\cH^1$ or $\C$, equipped with a K\"ahler metric $g_i$.  Taking the
conjugate complex structure on one of the factors defines a \ka structure
$(g,I)$ on $M$, with the opposite orientation to $(g,J)$.  By a result of
Kotschick \cite{kotschick}, $(M,J)$ is either a geometric complex surface~\cite{thurston} or is a minimal ruled surface.

If $(M,J)$ is geometric complex surface, the fundamental group acts
biholomorphically and isometrically with respect to the product of constant
curvature metrics on ${\mathbb U}_i$, i.e., $(M,J)$ carries a reducible
locally symmetric \ka structure.

If $(M,J)$ is a minimal ruled surface, then it is biholomorphic to the total
space of the projectivization $P(E)$ of a rank $2$ holomorphic vector bundle
$E$ over a compact complex curve $S$ (see for instance~\cite{BPV}) and so,
without loss, $\mathbb U_1$ is the universal cover of $S$ and $\mathbb U_2=\C
P^1$: the K\"ahler product metric $\tilde g=g_1\times g_2$ must be compatible
with the holomorphic splitting. If ${\mathbb U}_1= \C P^1$ as well, then $M =
{\C P}^1 \times {\C P}^1$ so it admits a product symmetric structure. Suppose
${\mathbb U}_1 = \C$ or $\C\cH^1$; by Liouville's Theorem, any holomorphic
isometry of $(\tilde M,\tilde g)$ has the form $\Psi(z,w) = (\psi_1(z),
\psi_2(z,w))$, where $\psi_1$ is a holomorphic isometry of $({\mathbb U}_1,
g_1)$ and (for any fixed $z$) $w\mapsto\psi_2(z,w)$ is a holomorphic isometry
of $(\C P^1, g_2)$. Since $\psi_1$ is a holomorphic automorphism of ${\mathbb
U}_1$, it preserves a constant curvature metric on ${\mathbb U}_1$; similarly,
since Isom$(g_2)$ is a compact subgroup of $PSL(2,\C)$, it lies in a conjugate
of $PSU(2)$ and hence preserves a constant curvature metric on $\C P^1$.  Thus
the fundamental group preserves the product of constant curvature metrics on
${\mathbb U}_1$ and $\C P^1$, so $(M,J)$ is again a geometric complex surface
supporting a reducible locally symmetric \ka structure.

\item Suppose now that $\phi$ has order $1$. By Theorem~\ref{thm:main}, after
blowing up $M$ at most once, we get a compact complex surface ${\hat M}$ which
is a holomorphic ${\C P}^1$-bundle over a compact complex curve ${S}$, i.e.,
$M$ is a ruled complex surface \cite{BPV}. If ${S} \cong {\C P}^1$, then $M$
is either ${\C P}^1\times{\C P}^1$ or a Hirzebruch surface $F_k = P({\mathcal
O}\oplus {\mathcal O}(k)) \to {\C P}^1$.  Of these surfaces, only $F_1$ is not
minimal: it is the blow-up ${\C P}^2$ at one point. We conclude that $M$ is
either ${\C P}^2$ or can be written as $P({\mathcal O}\oplus {\mathcal O}(k))
\to {\C P}^1$, $k\in\Z$.  If $S$ has genus {\bf g}$(S) \ge 1$, by using again
Theorem~\ref{thm:main}, we have $M = {\hat M}$ and therefore $M$ is a
(minimal) ruled surface $P(E)$ over a {compact} complex curve $S$, with the
induced ${\C}^{\times}$-action tangent to the projective fibers. Clearly in
the latter case $E$ must be split, and so without loss, $E =\cO\oplus \cL$.

As a final point, we have to show that we can assume $\deg \cL >0$ (or $k>0$
in the case of $F_k$). But this is an immediate consequence of
Theorem~\ref{thm:main}: the formula \eqref{chernclasses} specializes to give
(see also \eqref{curvform}) $c_1({\mathcal L}) = \frac12
c(\beta_1-\beta_0)[\omega_{S}/2\pi]$, where $\beta_0 < \beta_1$ and $c\neq 0$,
while $\pm \omega_{S}$ is the \ka structure induced on stable quotient
$S$. Thus $\deg\cL\neq 0$ and since $P(E) \cong P(E\otimes {\mathcal L}^*)$,
we can assume that $\deg\cL >0$.

\item This is an immediate consequence of Proposition~\ref{orthopoly}.

\smallbreak\noindent It follows from Theorem~\ref{thm:main} that each complex
surface listed in Theorem \ref{surfaces} does admit infinitely many
(non-isometric) \ka metrics with nontrivial hamiltonian $2$-forms of the
corresponding order.
\end{numlproof}

\begin{rem} The complex surfaces in Theorem \ref{surfaces} also admit
{\it extremal} \ka metrics with nontrivial hamiltonian $2$-forms, see
\cite{calabi1,christina1}.
\end{rem}

\subsection{Examples: extremal and weakly Bochner-flat K\"ahler metrics}
\label{s:ewbf}

We turn now to the construction of particular types of K\"ahler metrics with
hamiltonian $2$-forms. From this point of view, the notion of a hamiltonian
$2$-form is simply a device which provides constructions of interesting
K\"ahler manifolds, and this uses very little of the theory that we have
developed: the converse part of Theorem~\ref{thm:main}, which essentially
amounts to the sufficiency of the conditions for the compactification of a
toric K\"ahler metric and for the construction of K\"ahler metrics on
blow-downs.  In fact, we shall mainly restrict attention here to metrics on
projective line bundles (with no blow-downs) where these issues are trivial.

We recall from~\cite{ACG1} how Bochner-flat, weakly Bochner-flat and extremal
K\"ahler metrics with hamiltonian $2$-forms arise. A K\"ahler manifold $M$ is
\emph{Bochner-flat} if the Bochner tensor (a component of the K\"ahler
curvature) vanishes, \emph{weakly Bochner-flat} (WBF) if the Bochner tensor is
co-closed, and \emph{extremal} if the scalar curvature is a Killing potential
(i.e., its symplectic gradient is a Killing vector field).

By the differential Bianchi identity, a K\"ahler metric is WBF if and only if
its normalized Ricci form $\tilde\rho = \rho + \frac{\Scal}{2 (m + 1)} \,
\omega$ is a hamiltonian $2$-form. It follows that a WBF K\"ahler manifold is
extremal. Any K\"ahler--Einstein manifold is WBF, since $\tilde\rho$ a
constant multiple of $\omega$; however, the hamiltonian $2$-form in this case
is trivial. To deal with this, and the case of extremal K\"ahler metrics, we
shall suppose that there is a nontrivial hamiltonian $2$-form $\phi$ on $M$
such that $\tilde\rho=a\phi+b\omega$ in the case of WBF K\"ahler metrics, and
such that the scalar curvature $\Scal=a\trace_\omega\phi+b$ in the case of
extremal K\"ahler metrics (for constants $a,b$).

Suppose that we have a K\"ahler manifold $(M,g,J,\omega)$ with a hamiltonian
$2$-form $\phi$ of order $\ell$ where the K\"ahler quotient is a product of
$N$ K\"ahler manifolds $S_a$ of dimension $2d_a$, corresponding to the
constant roots $\eta_a$ of $\phi$. The K\"ahler metric then has the explicit
form~\eqref{metric} and there is the following local classification
result~\cite{ACG1}.
\begin{numlist}
\item $g$ is extremal, with $\Scal$ as above, if and only if
\begin{itemize}
\item for all $j$, $\pF_{j}''(t)=\check \Mpc(t) q(t)$, where $\check
\Mpc(t)=\prod_{a=1}^N (t-\eta_a)^{d_a-1}$ and $q$ is a polynomial of degree
$\ell+N$ independent of $j$;
\item for all $a$, $g_a$ has constant scalar curvature $q(\eta_a)/\prod_{b\neq
a} (\eta_a-\eta_b)$.
\end{itemize}
$g$ then has constant scalar curvature if and only if $q$ has degree
$\ell+N-1$.
\item $g$ is weakly Bochner-flat, with $\tilde\rho$ as above, if and only if
\begin{itemize}
\item for all $j$, $\pF_j'(t)= \Mpc(t)q(t)$, where $\Mpc(t)=\prod_{a=1}^N
(t-\eta_a)^{d_a}$ and $q$ is a polynomial of degree $\ell+1$ independent of
$j$;
\item for all $a$, $S_a$ is K\"ahler--Einstein with scalar curvature $d_a
q(\eta_a)$.
\end{itemize}
$g$ is then K\"ahler--Einstein if and only if $q$ has degree $\ell$.
\item \cite{bryant} $g$ is Bochner-flat, with $\tilde\rho$ as above, if and
only if
\begin{itemize}
\item for all $j$, $\pF_j(t)= \hat \Mpc(t)q(t)$ where $\hat\Mpc(t)=
\prod_{a=1}^N (t-\eta_a)^{d_a+1}$ and $q$ is a polynomial of degree $\ell+2-N$
independent of $j$;
\item for all $a$, $S_a$ has constant holomorphic sectional curvature and
scalar curvature $d_a(d_a+1)q(\eta_a)\prod_{b\neq a} (\eta_a-\eta_b)$.
\end{itemize}
$g$ has constant holomorphic sectional curvature if and only if $q$ has degree
$\ell+1-N$.
\end{numlist}

We want to combine this local classification with the global construction of
Theorem~\ref{thm:main}. To do this, we have to satisfy the \emph{boundary
conditions} of~\eqref{Theta-bound}, and the \emph{integrality conditions} for
the first Chern classes $c_1(\cL_{j,a})$ given by~\eqref{chernclasses}.

\begin{rem}\label{counting}  In practice, we need enough freedom in the
choice of $F(t)$ and the constants both to satisfy these boundary conditions
and to prescribe the first Chern classes freely (up to some open conditions),
since otherwise we face potentially nontrivial diophantine problems on our
data. Let us analyse the implications of this in the case that there are no
blow-downs, i.e., $M = P(\cL_{0} \oplus \cdots \oplus \cL_{\ell}) \to S$ where
$S$ is the product of compact Hodge complex manifolds $S_1, \ldots S_{N}$.
Thus we have $N(\ell+1)$ integrality conditions, together with $2(\ell+1)$
boundary conditions for the function $F(t)=\Mpc(t) \Theta(t)$, giving
$(N+2)(\ell+1)$ constraints on $F(t)$ and the constants $\beta_0,\ldots
\beta_\ell$, $\eta_1,\ldots \eta_N$ and $\cc,C_1,\ldots C_N$. Three of these
constants, say $\cc,\beta_0,\beta_\ell$ are useless for satisfying the
constraints, since there is a homethety freedom $g\mapsto k g$ in the K\"ahler
metric and an affine freedom $\xi_j\mapsto a\xi_j+b$ in the orthotoric
coordinates.  We therefore have $2N+\ell-1$ effective constants. This leaves
$N(\ell+1)+2(\ell+1)-2N-\ell+1=(N+1)(\ell-1)+4$ constraints on $F(t)$.
Subtracting this from the number of coefficients defining $F(t)$ gives the
expected dimension of the moduli space of solutions, which we require to be
nonnegative.
\begin{numlist}
\item In the extremal case, $F(t)$ is determined by $\ell+3+N$ constants,
giving $N(2-\ell)$ dimensional moduli and forcing $\ell\leq 2$.
\item In the WBF case, $F(t)$ is determined by $\ell+3$ constants, giving
$N(1-\ell)$ dimensional moduli and forcing $\ell\leq 1$
\item In the Bochner-flat case, $F(t)$ is determined by $\ell+3-N$ constants,
giving $-N\ell$ dimensional moduli and forcing $\ell=0$, $N\leq 2$. (This
paremeter count agrees with the classification of Bryant: the only compact
Bochner-flat K\"ahler manifolds are products of at most two constant
holomorphic sectional curvature manifolds.)
\end{numlist}
\end{rem}

We concentrate here on WBF K\"ahler metrics on projective line bundles, by
assuming the existence of a hamiltonian $2$-form $\phi$ of order $1$. In this
case, once we fix the base manifolds $S_a$ and the line bundles $\cL_{j,a}$,
the moduli are zero dimensional.

\subsubsection{The general setting}
In order to render our discussion as self-contained as possible, we first
recall our notations.  Let $(S_a,\pm g_a,\pm \omega_a)$, $a=1,\ldots N$, be
compact connected K\"ahler manifolds of real dimension $2d_a$, associated to
the distinct constant roots $\eta_a$ of the hamiltonian $2$-form. A K\"ahler
metric with a hamiltonian $2$-form of order $1$ is defined on a projective
line bundle over $S=S_1\times\cdots\times S_N$, using a metric of the form
\begin{equation} \label{lineartype} \begin{split}
g &= \sum_{a=1}^N (z-\eta_a) g_a
+ \frac{\prod_{a=1}^N (z-\eta_a)^{d_a}}{F(z)} \, dz^2
+ \frac{F(z)}{\prod_{a=1}^N (z-\eta_a)^{d_a}} \, \theta^2,\\
\omega &= \sum_{a=1}^N (z-\eta_a)\omega_a + dz\wedge\theta,
\qquad\qquad\qquad d\theta = \sum_{a=1}^N \omega_a,
\end{split} \end{equation}
where we normalize the momentum interval for $z$ to $[-1,1]$ and require
$|\eta_a|>1$.  Note that each K\"ahler metric $g_a$ can be positive or
negative definite, depending on the sign of $\eta_a$, and for convenience it
is taken here with the opposite sign to the one used in
equation~\eqref{metric}---observe that $\Mpn(\eta_a)= \eta_a-z$ rather than
$z-\eta_a$). It is convenient to set $\eta_a=-1/x_a$: now the sign of $g_a$ is
the sign of $x_a$.

The projective line bundle is $M=P(\cO \oplus \cL) \cong P(\cO \oplus
\cL^{-1})$, where up to a sign convention, $\theta$ is a connection form on
the principal $S^1$-bundle associated to $\cL$ with curvature $d\theta$. By
Theorem~\ref{thm:main}, $g$ compactifies on $M$ when $F(z)$ satisfies the
following boundary conditions (for the fibrewise compactification on $\C
P^1$):
\begin{equation} \label{boundary}
F(\pm1)=0,\qquad F'(\pm 1)= \mp 2 \Mpc(\pm 1).
\end{equation}
For the existence of $\cL$, we require that $\omega_a$ is integral, i.e.,
$[\omega_a/2\pi]$ is in the image of $H^2(S_a,\Z)$ in $H^2(S_a,\R)$, and we
write $\cL = \bigotimes_a \cL_a$, where $\cL_a$ is (the pullback to $M$ of) a
line bundle on $S_a$ with $c_1(\cL_a)=[\omega_a/2\pi]$.

In order to obtain WBF K\"ahler metrics, the $S_a$ must K\"ahler--Einstein,
i.e., with Ricci form $\rho_a = s_a\omega_a$. Since $[\rho_a/2\pi]$ is an
integral class, the first Chern class of the anti-canonical bundle,
$s_a=p_a/q_a$ for integers $p_a,q_a$. If $s_a\neq 0$, we take $p_a$ maximal so
that the anti-canonical bundle has a $p_a$th root (i.e., $[\rho_a/2\pi p_a]$
is a primitive class); then $\cL_a$ is $\cK_a^{\smash{-q_a/p_a}}$ twisted
by a flat line bundle.

\begin{rem} If $S_a$ is a Riemann surface $\Sigma_{\bf g}$ of genus ${\bf
g}$, then $p_a=2|{\bf g}-1|$, while if $S_a=\C P^{d_a}$, then $p_a=d_a+1$ so
that $\cK^{-1/p_a}=\cO(1)$. More generally, if the scalar curvature of $S_a$
is positive, then $p_a\leq d_a+1$ by Kobayashi--Ochiai~\cite{kob-och}.
\end{rem}

The remaining conditions to obtain a WBF metric (as in \S\ref{s:ewbf}(ii)
above) are
\begin{equation} \label{Fprime}
F'(z)=\Mpc(z)(b_{-1}z^{2}+b_{0}z +b_{1}),
\end{equation}
where $\Mpc(z) = \prod_{a=1}^N(z-\eta_a)^{d_a}$, and
\begin{equation} \label{einsteinconst}
2s_a=  b_{-1}\eta_a^2 + b_{0} \eta_a + b_{1}.
\end{equation}
Using the boundary conditions~\eqref{boundary} and the equation \eqref{Fprime}
for $F'$, we deduce that $b_{0}=-2$ and $b_{1}=-b_{-1}$. So, re-naming
$b_{-1}$ to $B$, equation \eqref{Fprime} becomes
\begin{equation}\label{newFprime}
F'(z)=\Mpc(z)\bigl(B(z^{2}-1)-2z\bigr)
\end{equation}
and~\eqref{einsteinconst} gives
\begin{equation}\label{ca}
B (1-x_a^2) = 2 x_a (x_a s_a - 1).
\end{equation}
$g$ is K\"ahler--Einstein if and only if $B=0$, which holds if and only if
$s_a=1/x_a$ for all $a$. (This implies in particular that the base factors
have positive scalar curvature.)

On the other hand, given the above, then \eqref{boundary} is satisfied if
and only if we set $F(z)=\int_{-1}^z \Mpc(t)\bigl(B (t^2-1) - 2 t\bigr) dt$
and
\begin{equation} \label{integralbc}
\int_{-1}^1 \Mpc(t)\bigl(B (t^2-1) - 2 t\bigr)  dt = 0.
\end{equation}
Since $F'(z)$ only changes sign once on the interval $(-1,1)$, $F(z)$ as
defined above will not have any zeroes between $z=-1$ and $z=1$.  Therefore,
as the sign of $F(z)$ equals the sign of $\Mpc(z)$ between $-1$ and $1$, the
metric $g$ will be positive definite.

So, in conclusion, the problem of constructing a WBF K\"ahler metric on $M$
(for given K\"ahler--Einstein manifolds $S_a$ with $s_a=p_a/q_a$) reduces to
finding solutions $B,x_1,\ldots x_N$ to \eqref{ca} and \eqref{integralbc}.
However, $\Mpc(t)(1-t^2)$ has constant sign on $(-1,1)$, so $B$ is uniquely
determined by~\eqref{integralbc}: substituting for $B$ from~\eqref{ca} (for
each $a$) it suffices to show that there exist distinct $(x_1,\ldots x_N)$
with $0<|x_a|<1$ such that
\begin{equation}\label{hdef}
h_a(x_1,\ldots x_N):= \int_{-1}^1 \widetilde\Mpc(t) H_a(t) dt
\end{equation}
vanishes for $a=1,\ldots N$, where $\widetilde\Mpc(t)=\prod_{b=1}^N (x_b
t+1)^{d_b}$ and
\begin{equation*}
H_a(t)= x_a(x_as_a-1)(1-t^2)+t (1-x_a^2)=x_a^2 s_a (1-t^2) + (t-x_a)
(x_a t+1)
\end{equation*}
\begin{rem} If $s_b\neq s_a$, $x_b$ cannot equal $x_a$, again since
$\Mpc(t)(1-t^2)$ has constant sign on $(-1,1)$. Hence if $x_a=x_b$, then
$s_a=s_b$ and $S_a\times S_b$ is K\"ahler--Einstein. Thus we do not actually
need to check that $x_1,\ldots x_N$ are distinct: if $x_a=x_b$, we still get a
WBF K\"ahler metric, but the hamiltonian $2$-form has fewer distinct constant
roots.
\end{rem}

\subsubsection{WBF K\"ahler metrics over K\"ahler--Einstein manifolds}

We consider the simplest case $N=1$, when $S=S_1$ is a K\"ahler--Einstein
manifold.  Replacing the momentum coordinate $z$ by $-z$ if necessary (and
dropping the $1$ subscripts) we may suppose that we have to find $0<x<1$ such
that $h(x)=0$, where
\begin{equation*}
h(x):=
\int_{-1}^1 (x t+1)^{d} \bigl( x (x s-1) (1-t^2) + t (1-x^2) \bigr) dz.
\end{equation*}
Since $h(0)=0$, $h'(0)=2(d-2)/3$ and the sign of $h(1)$ is equal to the sign
of $s-1$, we certainly have a solution $0<x<1$ to $h(x)=0$ if $d>2$ and
$s<1$.

For the case $d=2$ we calculate directly that
\begin{equation}
h(x)=\frac{4x^2}{15} \bigl( s (x^2+5) - 6x \bigr)
\end{equation}
and there is a solution $0<x<1$ to $h(x)=0$ if and only if $0<s<1$.

\begin{thm} \label{t:wbf1}
There are WBF K\"ahler metrics of the form \eqref{lineartype} on\textup:
\begin{bulletlist}
\item $P(\cO \oplus \cL) \to S$, where $S$ is a compact Ricci-flat K\"ahler
manifold of complex dimension $\geq 3$ whose K\"ahler form $\omega_S$ is
integral, and $\cL$ is a holomorphic line bundle with
$c_1(\cL)=[\omega_S/2\pi]$\textup;

\item $P(\cO \oplus \cK^{-q/p}{\otimes}\cL_0 ) \to S$, where $S$ is a
compact negative K\"ahler--Einstein manifold of complex dimension $\geq3$,
$q\in \Z$ with $q<0$, $\cK$ is the canonical bundle on $S$, and $\cL_0$ is a
flat line bundle on $S$\textup;

\item $P(\cO \oplus \cK^{-q/p}{\otimes}\cL_0 ) \to S$, where $S$ is a
compact positive K\"ahler--Einstein manifold of complex dimension $\geq 2$,
$q\in \Z$ with $q>p$, $\cK$ is the canonical bundle on $S$, and $\cL_0$ is a
flat line bundle on $S$.
\end{bulletlist}
\end{thm}

For the case $d=1$, we compute that
\begin{equation}\label{m1h}
h(x)=-\frac{2x}{3} \bigl( x^2+1 - 2sx \bigr)
\end{equation}
and there is a solution $0<x<1$ to $h(x)=0$ if and only if $s>1$. Since $S$ in
this case is $\C P^1$, $\cK=\cO(2)$ and the only possibility is $s=2$,
$\cL=\cO(1)$, in accordance with the classification of~\cite{ACG0}.

\subsubsection{WBF K\"ahler metrics over products of two K\"ahler--Einstein
manifolds}

In this section, we give a taste of the case $N=2$, but we postpone a more
thorough analysis to a subsequent paper. In this case we are looking for
common zeros of the functions
\begin{equation*}\notag
h_a(x_1, x_2):=
\int_{-1}^1 (x_1 t+1)^{d_1}(x_2 t+1)^{d_2}
\bigl( x_a (x_a s_a -1)(1-t^2) + t (1-x_a^2) \bigr) dt
\end{equation*}
(for $a=1,2$) with $0<|x_a|<1$. Analysing this problem in general involves
some delicate calculus arguments, but there are some special cases which are
straightforward. One of the simplest is the case that $d_1=d_2$, and
$s_1=-s_2$, when symmetry solves the problem for us, and we recover some of
the K\"ahler--Einstein metrics of Koiso and Sakane.

\begin{thm} \cite{koi-sak1,koi-sak2}
On the total space of $P(\cO\oplus \cO(k,-k)) \to \C P^d \times \C P^{d}$,
with $1\leq k\leq d$, there is a K\"ahler--Einstein metric, given
\textup(on a dense open set\textup) by
\begin{equation*}
g=\Bigl(\frac{d+1}{k}+z\Bigr)g_1 + \Bigl(\frac{d+1}{k}-z\Bigr)g_2
+ \frac{z^{2}-\frac{(d+1)^2}{k^2}}{F(z)} \,dz^{2} +
\frac{F(z)}{z^{2}-\frac{(d+1)^2}{k^2}} \,\theta^{2},
\end{equation*}
where $(g_1,\omega_1)$ and $(g_2,\omega_2)$ are Fubini--Study metrics on the
$\C P^d$ factors with holomorphic sectional curvature $2/k$, $d\theta=
\omega_1-\omega_2$ and $F(z)= \int_{-1}^{z}
2t\bigl(\frac{(d+1)^2}{k^2}-t^{2}\bigr)\,dt = -\frac{(d+1)^2}{k^2}(1-z^2) +
\frac12 (1-z^4)$.
\end{thm}
\begin{proof}
Let $s_1=-s_2=\frac{d+1}{k}$ and $x_1=-x_2=\frac{k}{d+1}$.  Then clearly
$0<|x_a|<1$ and $h_a(x_1,x_2)=0$ for $a=1,2$. Further, $x_a=1/s_a$ so the WBF
metric is K\"ahler--Einstein.
\end{proof}

In a subsequent paper, we generalize these metrics by proving the following.

\begin{thm} There is a WBF K\"ahler metric on the total space
$P(\cO\oplus\cO(k_1,k_2))\to\C P^{d_1}\times\C P^{d_2}$ in the following
cases\textup:
\begin{bulletlist}
\item $k_1>d_1+1$ and $k_2>d_2+1$\textup;
\item $1\leq k_1\leq d_1$ and $1\leq -k_2\leq d_2$.
\end{bulletlist}
\end{thm}

We illustrate this with the case $P(\cO\oplus\cO(1,-2))\to\C P^{2}\times\C
P^{3}$, where $d_1=2$, $d_2=3$, $s_1=3$, $s_2=4/(-2)=-2$.  The graphs of
$h_{1}=0$ (solid) and $h_{2}=0$ (dashed) for $0<x_1<1$ and $-1<x_2<0$ are
plotted below. Proving that the graphs do cross as shown is a tedious calculus
exercise.

\begin{figure}[ht]
\begin{center}
\includegraphics[width=.3\textwidth]{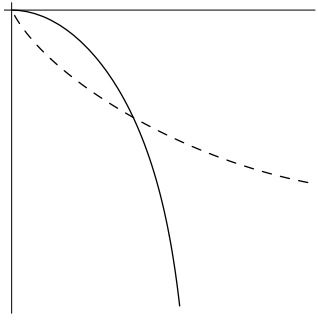}
\caption{$d_{1}=2$, $d_{2}=3$, $s_{1}=3$, $s_2=-2$}
\end{center}
\end{figure}

We end this section by giving an example with a blow-down.  Consider again
$P(\cO\oplus\cO(1,-1))\to\C P^1\times \C P^1$. This carries a Koiso--Sakane
K\"ahler--Einstein metric by setting $x_1=-x_2=1/2$, but it also admits two
blow-downs in which a $\C P^1$ factor collapses at an endpoint of the momentum
interval $[-1,1]$.  Such a collapse corresponds to setting $x_1=1$ and/or
$x_2=-1$. If we carry out both blow-downs, the resulting manifold is $\C P^3$,
which admits a WBF metric, namely the Fubini--Study metric, so let us consider
the case of a single blow-down. The two complex manifolds we obtain are both
isomorphic to $P(\cO\oplus\cO(1)\otimes\C^2)\to\C P^1$, so without loss, we
suppose $x_1=1$, $-1<x_2<0$.

This changes the boundary condition for $F'(z)$ at $z=-1$, since
$\Mpc(z)=(z+1)(z+1/x_2)$ vanishes at one of the endpoints. By a
straightforward application of L'H\^opital's rule we obtain
\begin{equation*}
(F'/\Mpc)(-1) = 4
\end{equation*}
Thus~\eqref{newFprime} is replaced by:
\begin{equation*}
F'(z)= (z+1)(z+1/x_2)\bigl(B(z^{2}-1)-2z+z(z-1)\bigr).
\end{equation*}
Setting $x_1=1$ automatically solves one of the integrality constraint with
$s_1=2$, so it remains to show that we can find $x_2$ to satisfy the second
constraint, with $s_2=-2$. Proceeding as in the case of no blow-downs, this
reduces to showing there is $-1<x<0$ with $f(x)=0$ where
\begin{equation*}
f(x)= \int_{-1}^1 (t+1)(x t+1) \bigl( -x (2x+1) (1-t^2) + t (1-x^2) +
\tfrac12 (x-1)(t-1)(xt+1)\bigr)dt.
\end{equation*}
This holds because $f(-1)$ is negative, while $f(0)=0$ and $f'(0)$ is
negative. Further $f(-1/2)$ is nonzero so the solution $x=x_2$ does not equal
$1/s_2$, and the metric is not K\"ahler--Einstein. 

\begin{thm} There is a WBF K\"ahler metric on $P(\cO\oplus\cO(1)\otimes\C^2)
\to \C P^1$ whose normalized Ricci form is a hamiltonian $2$-form of order
one. In particular, this is an extremal K\"ahler metric with non-constant
scalar curvature.
\end{thm}

\subsubsection{Further extremal K\"ahler metrics}

Since any WBF metric is extremal, the results presented so far provide new
examples of extremal K\"ahler metrics on projective line bundles and their
blow-downs. Furthermore, to obtain an extremal K\"ahler metric, it suffices
that the base manifolds $S_a$ are Hodge K\"ahler manifolds of constant scalar
curvature, giving examples which are not WBF in general.

On the other hand, such an approach is not very satisfactory, since it
produces only one extremal K\"ahler metric in each case, whereas the parameter
count of Remark~\ref{counting} suggests that these metrics should come in $N$
dimensional families (parameterized by admissible K\"ahler classes on $M$).
When the base manifolds have non-negative scalar curvatures, we can obtain
such $N$ dimensional families.

\begin{thm} For $a=1,\ldots N$, let $(S_a,\pm\omega_a)$ be Hodge K\"ahler
manifolds of constant nonnegative scalar curvature, let $\cL_a$ be a
holomorphic line bundles on each $S_a$ with $c_1(\cL_a)=[\omega_a/2\pi]$ and
let $\cL=\bigotimes_{a=1}^N \cL_a$. Then $M= P(\cO\oplus\cL)$ admits an $N$
parameter family of extremal K\"ahler metrics. Furthermore, if the K\"ahler
forms $\pm\omega_a$ do not all have the same sign \textup(i.e., if $c_1(\cL)$
is strictly indefinite\textup) there is an $N-1$ dimensional subfamily of
constant scalar curvature K\"ahler metrics on $M$.
\end{thm}
This Theorem generalizes results of Hwang~\cite{hwang} and
Hwang--Singer~\cite{hwang-singer}, and the proof is not materially different.
The first of these two papers considers the case that the base manifold has
constant eigenvalues of the Ricci tensor (e.g., a product of
K\"ahler--Einstein manifolds) and the idea to weaken this condition is
explored in the second paper. However, it is the notion of a hamiltonian
$2$-form that has selected for us a more general hypothesis for the base.  We
shall discuss this, and further results, in more detail in a subsequent paper.

Finally, we remark that the parameter count of Remark~\ref{counting} suggests
that one should be able to construct examples of extremal K\"ahler metrics on
projective plane bundles (and their blow downs) over products of constant
scalar curvature manifolds. Unfortunately, the existence problem here is
considerably less tractible than in the case of WBF metrics on projective line
bundles. Nevertheless, we hope to be able obtain examples in subsequent work.

\end{document}